\documentclass[pdflatex,sn-mathphys-ay]{sn-jnl}

\usepackage[utf8]{inputenc}
\usepackage[T1]{fontenc}
\usepackage{amsmath,amssymb,amsfonts,bm,amsthm}
\usepackage{graphicx}
\usepackage{subfigure}
\usepackage{manyfoot}
\usepackage{url}
\usepackage{xcolor}

\graphicspath{{./art/}{./figures/}{./}}

\newcommand{\IfLegacyFigureExists}[3]{%
  \IfFileExists{#1}{#2}{%
    \IfFileExists{art/#1}{#2}{%
      \IfFileExists{figures/#1}{#2}{#3}}}}
\newcommand{\IfTwoLegacyFiguresExist}[4]{%
  \IfLegacyFigureExists{#1}{\IfLegacyFigureExists{#2}{#3}{#4}}{#4}}
\newcommand{\IfThreeLegacyFiguresExist}[5]{%
  \IfLegacyFigureExists{#1}{\IfLegacyFigureExists{#2}{\IfLegacyFigureExists{#3}{#4}{#5}}{#5}}{#5}}
\newcommand{\IfFourLegacyFiguresExist}[6]{%
  \IfLegacyFigureExists{#1}{\IfLegacyFigureExists{#2}{\IfLegacyFigureExists{#3}{\IfLegacyFigureExists{#4}{#5}{#6}}{#6}}{#6}}{#6}}
\newif\iflegacysamplepaths
\legacysamplepathstrue
\IfLegacyFigureExists{1.5,0.4,a=-0.93.pdf}{}{\legacysamplepathsfalse}
\IfLegacyFigureExists{a=-0.93.pdf}{}{\legacysamplepathsfalse}
\IfLegacyFigureExists{1.5,1.3,a=-0.93.pdf}{}{\legacysamplepathsfalse}
\IfLegacyFigureExists{1.5,0.4,a=0.pdf}{}{\legacysamplepathsfalse}
\IfLegacyFigureExists{a=0.pdf}{}{\legacysamplepathsfalse}
\IfLegacyFigureExists{1.5,1.3,a=0.pdf}{}{\legacysamplepathsfalse}
\IfLegacyFigureExists{1.5,0.4,beta=-0.8.pdf}{}{\legacysamplepathsfalse}
\IfLegacyFigureExists{b=-0.8.pdf}{}{\legacysamplepathsfalse}
\IfLegacyFigureExists{1.5,1.3,beta=-0.8.pdf}{}{\legacysamplepathsfalse}
\IfLegacyFigureExists{1.5,0.4,10.6.pdf}{}{\legacysamplepathsfalse}
\IfLegacyFigureExists{a=10.6.pdf}{}{\legacysamplepathsfalse}
\IfLegacyFigureExists{1.5,1.3,a=10.6.pdf}{}{\legacysamplepathsfalse}

\newtheorem{proposition}{Proposition}
\newtheorem{example}{Example}
\newtheorem{thm}{Theorem}
\newtheorem{definition}{Definition}
\newtheorem{lem}{Lemma}
\newtheorem{remark}{Remark}

\begin{document}

\title[Weighted Sub-fractional Brownian Motion]{Weighted Sub-fractional Brownian Motion: Covariance Structure, Path Properties and Euler Approximation}

\author*[1]{\fnm{Jose H.} \sur{Ramirez-Gonzalez}}\email{hermenegildo.ramirez@usp.br}

\affil*[1]{\orgdiv{Institute of Mathematics and Statistics}, \orgname{University of S\~ao Paulo (USP)}, \city{S\~ao Paulo}, \country{Brazil}}

\abstract{
Weighted sub-fractional Brownian motion was introduced as a centered Gaussian process with covariance
\[
Q_{a,b}(s,t)=\frac{1}{1-b}\int_0^{s\wedge t}u^a\big\{(s-u)^b+(t-u)^b-(s+t-2u)^b\big\}\,du.
\]
This kernel was proved to be positive definite for \(a>-1\) and \(b\in[0,2]\setminus\{1\}\), and also for \(a>-1\), \(-1<b<0\), and \(a+b+1\ge0\). In this paper, we replace the power weight \(u^a\) by a measurable non-negative function \(f:\mathbb R_+\to\mathbb R_+\) satisfying
\[
\int_0^t f(u)(t-u)^b\,du<\infty,\qquad t>0,
\]
and consider the kernel
\[
R_{f,b}(s,t)=\frac{1}{1-b}\int_0^{s\wedge t}f(u)\big\{(s-u)^b+(t-u)^b-(s+t-2u)^b\big\}\,du.
\]
We prove that, for \(b>-1\), \(b\ne1\), the kernel \(R_{f,b}\) is positive definite for every such function \(f\) if and only if \(b\in[0,1)\cup(1,2]\). For every \(b\in(-1,0)\cup(2,\infty)\), we construct a non-negative function \(f\) satisfying the integrability condition for which \(R_{f,b}\) is not positive definite. The case \(b=1\) is obtained as a logarithmic limit. For the associated Gaussian process, we study H\"older regularity, total and quadratic variation, non-stationarity, and long-range dependence. For \(b\in(0,1)\cup(1,2]\), we also define differential equations driven by this process and establish pathwise convergence of the corresponding Euler approximation, together with strong \(L^p\)-rates under the additional regularity assumptions stated below. The logarithmic boundary \(b=1\) is included in the covariance theory but is not part of the pathwise-equation and Euler results.
}

\keywords{weighted sub-fractional Brownian motion, Gaussian processes, Young integration, Euler approximation, strong \(L^p\) convergence, long-range dependence}

\pacs[MSC 2020]{60G15, 60G17, 60G22, 60H10, 65C30}

\maketitle

\section{Introduction}\label{sec:introduction}
Weighted sub-fractional Brownian motion arose from the study of occupation-time fluctuations of age-dependent critical branching particle systems; see \cite{VW,FVW,LMR2024}. Its covariance is
\begin{equation}\label{intro-power-kernel}
Q_{a,b}(s,t)
=
\frac{1}{1-b}
\int_0^{s\wedge t}
u^a\bigl\{(s-u)^b+(t-u)^b-(s+t-2u)^b\bigr\}\,du,
\qquad s,t\ge0,
\end{equation}
where \(a>-1\), \(b>-1\), \(b\ne1\), and \(u^a\) is the temporal weight. In the constant-weight case \(a=0\), direct integration gives
\[
Q_{0,b}(s,t)
=
\frac{1}{(1-b)(1+b)}
\left[
s^{b+1}+t^{b+1}
-\frac12\bigl\{(s+t)^{b+1}+|s-t|^{b+1}\bigr\}
\right].
\]
Thus, up to a multiplicative constant, \(Q_{0,b}\) coincides with the sub-fractional Brownian covariance with parameter \(h=b+1\). This covariance appears in occupation-time fluctuation limits \cite{BGT-SPL2004} and is positive definite exactly for \(0<h\le4\). The range \(0<h\le2\) corresponds to sub-fractional Brownian motion, whereas \(2<h\le4\) corresponds to negative sub-fractional Brownian motion \cite{BGT2007}. Related occupation-time fluctuation limits led to weighted fractional Brownian motion and to the study of its path properties \cite{BGT-gamma,Propier}.

At \(b=1\), equivalently \(h=2\), the preceding proportionality factor
is singular and the usual sub-fractional covariance in square brackets
vanishes identically. The normalized limit considered here is instead the
non-degenerate logarithmic covariance introduced below.

The positive definiteness of \(Q_{a,b}\) was studied in \cite{LMR2024}. For \(a>-1\), it was proved that \(Q_{a,b}\) is positive definite when \(b\in[0,1)\cup(1,2]\), and also when \(-1<b<0\) and \(a+b+1\ge0\). It was further shown that positive definiteness fails when \(-1<b<0\) and \(a+b+1<0\), and when \(b>a+3\). Therefore, the positive-definiteness problem remained unresolved in general in the region
\[
a>-1,\qquad 2<b\le a+3.
\]
For \(a=0\), the answer follows from the known covariance range of sub-fractional and negative sub-fractional Brownian motion \cite{BGT2007}.

In this paper, we replace the power weight \(u^a\) by a measurable non-negative function \(f:\mathbb R_+\to\mathbb R_+\). For \(b>-1\), \(b\ne1\), we define
\begin{equation}\label{intro-general-kernel}
R_{f,b}(s,t)=\frac{1}{1-b}\int_0^{s\wedge t}f(u)\bigl\{(s-u)^b+(t-u)^b-(s+t-2u)^b\bigr\}\,du,\qquad s,t\ge0,
\end{equation}
and assume that
\[
\int_0^t f(u)(t-u)^b\,du<\infty,\qquad t>0.
\]
A measurable non-negative function satisfying this condition is called an admissible weight for the parameter \(b\).

We determine the values of \(b\) for which \(R_{f,b}\) is positive definite for every admissible weight. We prove that, among \(b>-1\), \(b\ne1\), this holds exactly when
\(
b\in[0,1)\cup(1,2].
\)

For each \(b\in(-1,0)\cup(2,\infty)\), we construct an admissible non-negative weight for which \(R_{f,b}\) is not positive definite. Hence \([0,1)\cup(1,2]\) is the maximal range in which \eqref{intro-general-kernel} defines a covariance for every admissible weight. This result does not settle the positive definiteness of every individual power weight in the remaining region \(2<b\le a+3\).

At \(b=1\), the normalized kernels converge to a logarithmic covariance. A general form of this covariance was studied in \cite{gonzalez2024} for non-negative locally integrable weights. In the notation of the present paper, the covariance used there is \(2K_f\). The corresponding Gaussian processes were also fitted separately to the latitude and longitude components of bat trajectories obtained from telemetry data. This application motivates the use of weights beyond the power family \(f(u)=u^a\).

The family \eqref{intro-general-kernel} includes several known cases. If \(f\equiv c>0\), then \(R_{f,b}\) is proportional to the sub-fractional covariance with exponent \(b+1\). If \(b=0\), the associated process is a time-changed Brownian motion with clock
\[
t\longmapsto\int_0^t f(u)\,du.
\]
Let \(\zeta_{f,b}=(\zeta_{t,f,b})_{t\ge0}\) denote the centered Gaussian process with covariance \(R_{f,b}\). For \(0<b<1\), we study H\"older regularity, quadratic variation, and total variation. For \(1<b\le2\), under a local power bound on \(f\), we prove that the process admits an absolutely continuous realization. We also study non-stationarity and long-range dependence.

Finally, we define differential equations driven by \(\zeta_{f,b}\) and study the corresponding Euler approximation. For \(0<b<1\), the integral is defined in the Young sense, and we obtain pathwise and strong \(L^p\)-convergence. For \(1<b\le2\), the equation is treated pathwise through the absolutely continuous realization of the process; uniform pathwise convergence holds under a local power bound on \(f\), while the order-one strong \(L^p\)-estimate is obtained when the Gaussian velocity has a continuous modification.

Section~\ref{sec:covariance-construction} establishes the covariance range. Section~\ref{sec:path-properties} states the path and dependence results. Section~\ref{sec:pathwise-euler} treats differential equations driven by the process and their Euler approximation. Detailed proofs of the sample-path and Euler results, related covariance constructions on Euclidean index sets, closed covariance formulas, and numerical experiments are provided in Sections~\ref{sec:supp-path-proofs}--\ref{sec:supp-paths} of the Supplementary Material included below.
\medskip

\noindent\textit{Notation.}
All random variables are defined on a complete probability space
\((\Omega,\mathcal F,\mathbb P)\). Throughout,
\(\mathbb R_+=[0,\infty)\),
\(a\wedge b=\min\{a,b\}\),
\(a\vee b=\max\{a,b\}\), and
\(\mathbf 1_D\) denotes the indicator of a set \(D\).
A symmetric real-valued kernel \(K\) on an index set is positive definite if
\[
\sum_{i,j=1}^m c_i c_j K(t_i,t_j)\ge0
\]
for every \(m\ge1\), every \(t_1,\ldots,t_m\), and every
\(c_1,\ldots,c_m\in\mathbb R\).

For a compact interval \(I\subset\mathbb R\), a function
\(x:I\to\mathbb R\), and \(0<\gamma<1\), define
\[
\|x\|_{\infty;I}:=\sup_{t\in I}|x_t|,
\qquad
[x]_{\gamma;I}:=
\sup_{\substack{s,t\in I\\s\ne t}}
\frac{|x_t-x_s|}{|t-s|^\gamma},
\qquad
\|x\|_{\gamma;I}:=\|x\|_{\infty;I}+[x]_{\gamma;I}.
\]
The space \(C^\gamma(I)\) consists of the functions with finite
\(\|\cdot\|_{\gamma;I}\). When \(I=[0,T]\), or when the interval is clear
from the context, it is omitted from the notation.

For a partition
\(\pi=\{a=t_0<\cdots<t_n=d\}\), let
\[
|\pi|:=\max_{0\le k<n}(t_{k+1}-t_k),
\]
and write \([u,v]\in\pi\) when \(u\) and \(v\) are consecutive partition
points. The total variation of \(x\) on \(I=[a,d]\) is
\[
V_I(x):=
\sup_\pi\sum_{k=0}^{n-1}|x_{t_{k+1}}-x_{t_k}|,
\]
where the supremum is taken over all finite partitions of \(I\).

We write
\[
Z_{s,t}:=Z_t-Z_s,
\qquad
\operatorname{Lip}(F):=
\sup_{x\ne y}\frac{|F(x)-F(y)|}{|x-y|},
\qquad
\|Y\|_{L^p}:=(\mathbb E|Y|^p)^{1/p}.
\]
Unless otherwise stated, semimartingale assertions refer to the completed
natural filtration of the chosen continuous modification.

Symbols such as \(C\) and \(C_\theta\) denote finite positive constants that
may change from line to line and depend only on the indicated parameters.
The notation \(O(\cdot)\), \(o(\cdot)\), and \(\sim\) has its usual asymptotic
meaning. Equality of finite-dimensional distributions is denoted by
\(\stackrel{\mathrm{f.d.d.}}{=}\). Finally,
\[
B(p,q):=\int_0^1 r^{p-1}(1-r)^{q-1}\,dr,
\qquad p,q>0,
\]
denotes the Beta function, and \(\Gamma\) denotes the Gamma function.
\section{Weighted sub-fractional Brownian motion}

\subsection{Covariance construction and maximal universal range}
\label{sec:covariance-construction}

Let \(b>-1\), \(b\ne1\), and let
\(f:\mathbb R_+\to\mathbb R_+\) be measurable. Assume that
\begin{equation}\label{eq:weight-finiteness}
\int_0^t f(u)(t-u)^b\,du<\infty,
\qquad t>0.
\end{equation}
A non-negative measurable function satisfying
\eqref{eq:weight-finiteness} is called an admissible weight for the
parameter \(b\).

For \(s,t\ge0\), define
\begin{equation}\label{weighted-covariance}
R_{f,b}(s,t)
:=
\frac{1}{1-b}
\int_0^{s\wedge t}
f(r)
\bigl[
(s-r)^b+(t-r)^b-(s+t-2r)^b
\bigr]\,dr.
\end{equation}
When \(b<0\), the integral is understood in the improper sense at the
endpoints where a power vanishes. The value assigned to the integrand at
a single endpoint does not affect the integral.

For the boundary value \(b=1\), we assume instead that \(f\) is locally
integrable:
\begin{equation}\label{eq:log-weight-finiteness}
\int_0^t f(u)\,du<\infty,
\qquad t>0.
\end{equation}

The following theorem determines the maximal range of \(b\) for which
\eqref{weighted-covariance} is positive definite for every admissible
weight.

\begin{thm}\label{thm:maximal-range}
Let \(b>-1\), \(b\ne1\).  The following statements hold.
\begin{itemize}
\item[(i)] If \(b\in[0,1)\cup(1,2]\), then for every measurable \(f\ge0\) satisfying \eqref{eq:weight-finiteness}, the kernel \(R_{f,b}\) is finite and positive definite.  Hence there exists a centered Gaussian process \(\zeta_{f,b}=\{\zeta_{t,f,b}:t\ge0\}\) with covariance \(R_{f,b}\).
\item[(ii)] If \(b\in(-1,0)\cup(2,\infty)\), then there exists a non-negative locally integrable function \(f\) satisfying \eqref{eq:weight-finiteness} for which \(R_{f,b}\) is not positive definite.
\end{itemize}
Consequently, among \(b>-1\) with \(b\ne1\), the set \([0,1)\cup(1,2]\) is the largest range for which the construction is a covariance for every admissible non-negative weight.
\end{thm}

\begin{proof}
We first prove part~\textup{(i)}. For
\(b\in(0,1)\cup(1,2]\), define
\[
Q_b(x,y)
:=
\frac{x^b+y^b-(x+y)^b}{1-b},
\qquad x,y\ge0.
\]
For \(b=0\), set
\[
Q_0(x,y)
:=
\mathbf 1_{(0,\infty)}(x)
\mathbf 1_{(0,\infty)}(y),
\qquad x,y\ge0.
\]
The kernel \(Q_0\) is positive definite because it has rank one. Moreover,
up to changes at the endpoints of the integral, which do not affect its
value,
\[
R_{f,0}(s,t)
=
\int_0^{s\wedge t}f(r)\,dr.
\]
Thus, if
\[
A(t):=\int_0^t f(r)\,dr,
\]
then \(A\) is finite, continuous, and non-decreasing, and
\[
R_{f,0}(s,t)=A(s\wedge t)
\]
is the covariance of the time-changed Brownian motion
\((B(A(t)))_{t\ge0}\).

We now consider \(0<b<2\), \(b\ne1\). We claim that
\begin{equation}\label{eq:Qb-laplace-feature}
Q_b(x,y)
=
\frac{b}{\Gamma(2-b)}
\int_0^\infty
(1-e^{-\lambda x})(1-e^{-\lambda y})
\lambda^{-b-1}\,d\lambda,
\qquad x,y\ge0.
\end{equation}
The integral is finite. Indeed, as \(\lambda\downarrow0\),
\[
(1-e^{-\lambda x})(1-e^{-\lambda y})
\lambda^{-b-1}
=
O(\lambda^{1-b}),
\]
which is integrable at zero because \(b<2\). As
\(\lambda\to\infty\), the integrand is
\(O(\lambda^{-b-1})\), which is integrable because \(b>0\).

Let \(F_b(x,y)\) denote the right-hand side of
\eqref{eq:Qb-laplace-feature}. For \(x,y>0\), differentiation under the
integral sign is justified locally by dominated convergence and gives
\[
\frac{\partial^2 F_b}{\partial x\,\partial y}(x,y)
=
\frac{b}{\Gamma(2-b)}
\int_0^\infty
e^{-\lambda(x+y)}\lambda^{1-b}\,d\lambda
=
b(x+y)^{b-2}.
\]
On the other hand,
\[
\frac{\partial^2 Q_b}{\partial x\,\partial y}(x,y)
=
b(x+y)^{b-2}.
\]
Both \(F_b\) and \(Q_b\) vanish when either argument is zero. Hence their
difference has zero mixed derivative and vanishes on both coordinate
axes, which implies that \(F_b=Q_b\) on \(\mathbb R_+^2\). This proves
\eqref{eq:Qb-laplace-feature}.

For each \(\lambda>0\), the kernel
\[
(x,y)\longmapsto
(1-e^{-\lambda x})(1-e^{-\lambda y})
\]
has rank one and is therefore positive definite. Since
\[
\frac{b}{\Gamma(2-b)}\lambda^{-b-1}\,d\lambda
\]
is a positive measure, representation
\eqref{eq:Qb-laplace-feature} shows that \(Q_b\) is positive definite for
\(0<b<2\), \(b\ne1\). At \(b=2\),
\[
Q_2(x,y)
=
\frac{x^2+y^2-(x+y)^2}{1-2}
=
2xy,
\]
which is again a rank-one positive-definite kernel.

Fix \(r\ge0\) and define
\[
K_r(s,t)
:=
\mathbf 1_{[r,\infty)}(s)
\mathbf 1_{[r,\infty)}(t)
Q_b(s-r,t-r).
\]
The kernel \(K_r\) is positive definite. Indeed, for arbitrary
\(t_1,\ldots,t_m\ge0\) and \(a_1,\ldots,a_m\in\mathbb R\),
\[
\sum_{i,j=1}^m a_i a_j K_r(t_i,t_j)
=
\sum_{\substack{1\le i,j\le m\\ t_i,t_j\ge r}}
a_i a_j Q_b(t_i-r,t_j-r)
\ge0,
\]
because \(Q_b\) is positive definite.

For \(b\in[0,1)\cup(1,2]\), one has
\begin{equation}\label{eq:R-mixture-K}
R_{f,b}(s,t)
=
\int_0^\infty f(r)K_r(s,t)\,dr
=
\int_0^{s\wedge t}
f(r)Q_b(s-r,t-r)\,dr.
\end{equation}
We verify that this integral is absolutely convergent. The
Cauchy--Schwarz inequality for the positive-definite kernel \(Q_b\)
gives
\[
|Q_b(x,y)|^2
\le Q_b(x,x)Q_b(y,y).
\]
Consequently, applying the ordinary Cauchy--Schwarz inequality with
respect to the measure \(f(r)\,dr\),
\begin{align*}
&\int_0^{s\wedge t}
f(r)|Q_b(s-r,t-r)|\,dr\\
&\quad\le
\left(
\int_0^{s\wedge t}
f(r)Q_b(s-r,s-r)\,dr
\right)^{1/2}
\left(
\int_0^{s\wedge t}
f(r)Q_b(t-r,t-r)\,dr
\right)^{1/2}\\
&\quad\le
\left(
\int_0^s
f(r)Q_b(s-r,s-r)\,dr
\right)^{1/2}
\left(
\int_0^t
f(r)Q_b(t-r,t-r)\,dr
\right)^{1/2}.
\end{align*}
For \(b\ne0\),
\[
Q_b(x,x)
=
\frac{2-2^b}{1-b}x^b,
\qquad x\ge0.
\]
The constant
\[
c_b:=\frac{2-2^b}{1-b}
\]
is strictly positive for every \(b\ne1\). Therefore
\begin{equation}\label{diag-covariance-general}
R_{f,b}(t,t)
=
c_b\int_0^t f(u)(t-u)^b\,du,
\end{equation}
or equivalently,
\[
\frac{1-b}{2-2^b}R_{f,b}(t,t)
=
\int_0^t f(u)(t-u)^b\,du.
\]
The right-hand side is finite by
\eqref{eq:weight-finiteness}. Hence the two diagonal integrals in the
preceding Cauchy--Schwarz bound are finite. For \(b=0\), finiteness
follows directly from
\[
R_{f,0}(t,t)=\int_0^t f(r)\,dr<\infty.
\]
Thus \(R_{f,b}(s,t)\) is finite for every \(s,t\ge0\).

Since the sum below is finite and each entry is absolutely integrable,
we may interchange the finite sum and the integral in
\eqref{eq:R-mixture-K}. For arbitrary
\(t_1,\ldots,t_m\ge0\) and \(a_1,\ldots,a_m\in\mathbb R\),
\begin{align*}
\sum_{i,j=1}^m a_i a_jR_{f,b}(t_i,t_j)
&=
\int_0^\infty
f(r)
\sum_{i,j=1}^m a_i a_jK_r(t_i,t_j)\,dr\ge0,
\end{align*}
because \(f\ge0\) and \(K_r\) is positive definite for every \(r\ge0\).
Therefore \(R_{f,b}\) is finite and positive definite.

It follows from the standard existence theorem for Gaussian processes
that there exists a centered Gaussian process
\[
\zeta_{f,b}=\{\zeta_{t,f,b}:t\ge0\}
\]
with covariance \(R_{f,b}\). This proves part~\textup{(i)}.

We now prove part~\textup{(ii)}. For
\(b\in(-1,0)\cup(2,\infty)\), define \(Q_b\) on
\((0,\infty)^2\) by
\[
Q_b(x,y)
=
\frac{x^b+y^b-(x+y)^b}{1-b}.
\]
We first show that this kernel is not positive definite. Set
\[
c_b:=\frac{2-2^b}{1-b}>0.
\]
For \(M>1\), consider the determinant
\[
D_b(M)
:=
\det
\begin{pmatrix}
Q_b(1,1)&Q_b(1,M)\\
Q_b(M,1)&Q_b(M,M)
\end{pmatrix}.
\]
Since
\[
Q_b(1,1)=c_b,
\qquad
Q_b(M,M)=c_bM^b,
\]
we have
\[
D_b(M)
=
c_b^2M^b-Q_b(1,M)^2.
\]

Suppose first that \(-1<b<0\). Then \(M^b\to0\) and
\((M+1)^b\to0\) as \(M\to\infty\). Therefore
\[
Q_b(M,M)=c_bM^b\longrightarrow0
\]
and
\[
Q_b(1,M)
=
\frac{1+M^b-(1+M)^b}{1-b}
\longrightarrow
\frac{1}{1-b}.
\]
Consequently,
\[
D_b(M)\longrightarrow-\frac{1}{(1-b)^2}<0.
\]
Thus \(D_b(M)<0\) for all sufficiently large \(M\).

Suppose now that \(b>2\). As \(M\to\infty\),
\[
(M+1)^b
=
M^b+bM^{b-1}+O(M^{b-2}),
\]
and hence
\[
1+M^b-(M+1)^b
=
-bM^{b-1}(1+o(1)).
\]
Since \(1-b=-(b-1)\),
\[
Q_b(1,M)
=
\frac{b}{b-1}M^{b-1}(1+o(1)).
\]
It follows that
\[
D_b(M)
=
c_b^2M^b
-
\left(\frac{b}{b-1}\right)^2
M^{2b-2}(1+o(1)).
\]
Because \(2b-2>b\), the second term dominates the first one, and hence
\[
D_b(M)<0
\]
for every sufficiently large \(M\).

We have therefore shown that, for each
\(b\in(-1,0)\cup(2,\infty)\), there exists \(M>1\) such that
\begin{equation}\label{eq:negative-Q-determinant}
D_b(M)<0.
\end{equation}

Fix such an \(M\). For \(0<\varepsilon<1\), define
\[
f_\varepsilon(r)
:=
\frac{1}{\varepsilon}
\mathbf 1_{[0,\varepsilon]}(r).
\]
This function is non-negative and locally integrable. It also satisfies
\eqref{eq:weight-finiteness} for every \(b>-1\). Indeed, if
\(0<t\le\varepsilon\), then
\[
\int_0^t
f_\varepsilon(u)(t-u)^b\,du
=
\frac{1}{\varepsilon}
\int_0^t(t-u)^b\,du
=
\frac{t^{b+1}}{\varepsilon(b+1)}
<\infty.
\]
If \(t>\varepsilon\), then
\[
\int_0^t
f_\varepsilon(u)(t-u)^b\,du
=
\frac{1}{\varepsilon}
\int_0^\varepsilon(t-u)^b\,du
=
\frac{t^{b+1}-(t-\varepsilon)^{b+1}}
{\varepsilon(b+1)}
<\infty.
\]

Let
\[
x_1=1,
\qquad
x_2=M.
\]
Since \(0<\varepsilon<1\le x_i\wedge x_j\), for
\(i,j\in\{1,2\}\),
\begin{align*}
R_{f_\varepsilon,b}(x_i,x_j)
&=
\frac{1}{\varepsilon}
\int_0^\varepsilon
Q_b(x_i-r,x_j-r)\,dr\\
&=
\int_0^1
Q_b(x_i-\varepsilon v,x_j-\varepsilon v)\,dv.
\end{align*}
The function \(Q_b\) is continuous on \((0,\infty)^2\), and the
arguments in the last integral remain in a compact subset of
\((0,\infty)^2\) for all sufficiently small \(\varepsilon\).
Consequently,
\[
R_{f_\varepsilon,b}(x_i,x_j)
\longrightarrow
Q_b(x_i,x_j),
\qquad \varepsilon\downarrow0.
\]

Let
\[
D_{b,\varepsilon}(M)
:=
\det
\begin{pmatrix}
R_{f_\varepsilon,b}(1,1)&R_{f_\varepsilon,b}(1,M)\\
R_{f_\varepsilon,b}(M,1)&R_{f_\varepsilon,b}(M,M)
\end{pmatrix}.
\]
Since the determinant is a continuous polynomial in the matrix entries,
\[
D_{b,\varepsilon}(M)
\longrightarrow
D_b(M)<0.
\]
Hence \(D_{b,\varepsilon}(M)<0\) for every sufficiently small
\(\varepsilon>0\). The corresponding \(2\times2\) matrix is therefore
not positive semidefinite, and \(R_{f_\varepsilon,b}\) is not positive
definite.

This proves part~\textup{(ii)} and establishes the maximality of
\([0,1)\cup(1,2]\).
\end{proof}
\begin{definition}
A centered Gaussian process \(\{\zeta_{t,f,b}:t\ge0\}\) with covariance function \eqref{weighted-covariance}, where \(b\in(0,1)\cup(1,2]\) and \(f\) satisfies \eqref{eq:weight-finiteness}, is called a \emph{weighted sub-fractional Brownian motion}. For \(b=0\) the same name is used with the covariance convention
\[
        R_{f,0}(s,t)=\int_0^{s\wedge t}f(r)\,dr.
\]
\end{definition}

Two examples used below are \(f(u)=u^a\), \(a>-1\), and \(f(u)=e^{au}\), \(a\in\mathbb R\), both of which are admissible for every \(b>-1\).

\begin{remark}\label{rem:constant-weight-representations}
Let \(f\equiv c>0\). Direct integration of
\eqref{weighted-covariance} gives
\begin{equation}\label{eq:constant-weight-covariance}
R_{c,b}(s,t)
=
\frac{c}{(1-b)(1+b)}
\left[
s^{b+1}+t^{b+1}
-\frac12\bigl((s+t)^{b+1}+|t-s|^{b+1}\bigr)
\right],
\end{equation}
for every \(b>-1\), \(b\ne1\). In the positive-definite
constant-weight range \(-1<b<1\) or \(1<b\le3\), let
\(\zeta_{c,b}\) denote a centered Gaussian process with this covariance,
and set \(H=(b+1)/2\).

If \(-1<b<1\) and \(B^H\) is a two-sided fractional Brownian motion with covariance
\[
\mathbb E(B_t^HB_s^H)
=
\frac12
\bigl(
|t|^{2H}+|s|^{2H}-|t-s|^{2H}
\bigr),
\]
then
\[
\zeta_{t,c,b}
\stackrel{\mathrm{f.d.d.}}{=}
\left(
\frac{c}{2(1-b)(1+b)}
\right)^{1/2}
\bigl(B_t^H+B_{-t}^H\bigr),
\qquad t\ge0.
\]

If \(1<b\le3\) and \(\nu_b\) is a centered Gaussian process with covariance
\[
\mathbb E(\nu_{s,b}\nu_{t,b})
=
(s+t)^{b-1}-|s-t|^{b-1},
\]
then
\[
\zeta_{t,c,b}
\stackrel{\mathrm{f.d.d.}}{=}
\left(
\frac{bc}{2(b-1)}
\right)^{1/2}
\int_0^t \nu_{r,b}\,dr,
\qquad t\ge0.
\]
These representations agree, after normalization, with those in
\cite{BGT2007}.
\end{remark}

\subsection{Sample path properties}\label{sec:path-properties}
We consider separately the ranges \(0<b<1\) and \(1<b\le2\). For
\(0<b<1\), under the comparison assumptions stated below, we prove that
the process has zero quadratic variation and infinite total variation.
For \(1<b\le2\), if
\[
0\le f(u)\le cu^a,\qquad 0<u\le T,
\]
for some \(c>0\) and \(a>-1\), we prove that the process admits an
absolutely continuous, and therefore finite-variation, realization on
\([0,T]\).

We first consider \(0<b<1\). For non-negative functions \(f\) and \(g\),
we write \(f\asymp g\) on a set if there exist constants
\(0<c_1\le c_2<\infty\) such that
\[
c_1f\le g\le c_2f
\]
on that set. For the constant weight \(f\equiv1\),
\eqref{eq:constant-weight-covariance} and the increment bounds for
sub-fractional Brownian motion in \cite{BGT-SPL2004} give
\begin{equation}\label{eq:constant-increment-bounds}
\frac{2-2^b}{(1-b)(1+b)}(t-s)^{1+b}
\le
\mathbb E(\zeta_{t,1,b}-\zeta_{s,1,b})^2
\le
\frac{1}{(1-b)(1+b)}(t-s)^{1+b},
\qquad 0\le s<t.
\end{equation}

Throughout this subsection, zero quadratic variation means that, for
every deterministic sequence of partitions \(\pi_n\) whose mesh tends
to zero,
\[
\sum_{[u,v]\in\pi_n}
|\zeta_{v,f,b}-\zeta_{u,f,b}|^2
\longrightarrow0
\qquad\text{almost surely}.
\]
All pathwise statements refer to the continuous modifications specified
in the corresponding propositions. The detailed arguments are organized in
Section~\ref{sec:supp-path-proofs}:
Section~\ref{sec:supp-infinite-variation-lemma} proves the auxiliary
infinite-variation lemma, and
Sections~\ref{sec:supp-proof-prop1}--\ref{sec:supp-proof-prop8} prove
Propositions~\ref{prop:rough-bounded-weight}--\ref{prop:logarithmic-limit}, respectively.

\begin{proposition}\label{prop:rough-bounded-weight}
Let \(T>0\), let \(0<b<1\), let \(f\) be an admissible weight, and suppose that
\(1\asymp f\) on \([0,T]\). Then \(\zeta_{f,b}\) admits a continuous
modification which has zero quadratic variation and infinite total
variation almost surely on \([0,T]\). In particular, this modification
is not a semimartingale.
\end{proposition}

\begin{proposition}\label{prop:rough-power-comparable}
Let \(T>0\), let \(0<b<1\), and let \(f\) be an admissible weight.
Suppose that
\[
u^\alpha\asymp f(u)
\qquad\text{on }(0,T]
\]
for some \(\alpha>-1\). Then \(\zeta_{f,b}\) admits a modification
continuous on \((0,T]\) such that, almost surely, for every
\(\varepsilon\in(0,T)\), its restriction to \([\varepsilon,T]\) has
zero quadratic variation and infinite total variation. In particular,
for every \(\varepsilon\in(0,T)\), this restriction is not a
semimartingale with respect to its completed natural filtration.
\end{proposition}

\begin{proposition}\label{prop:gaussian-velocity}
Let \(1<b\le2\), fix \(T>0\), and assume that
\[
0\le f(u)\le cu^a,\qquad 0<u\le T,
\]
for some \(c>0\) and \(a>-1\). Then there exists a jointly measurable
centered Gaussian process
\[
\nu_{f,b}=\{\nu_{t,f,b}:0\le t\le T\}
\]
with covariance
\begin{equation}\label{eq:velocity-covariance}
C_{f,b}(s,t)
:=
\int_0^{s\wedge t}
f(u)(s+t-2u)^{b-2}\,du.
\end{equation}
When \(1<b<2\), the value of the integrand at the single endpoint
\(u=s=t\) is set equal to zero; this convention does not alter the
integral. Moreover,
\begin{equation}\label{eq:velocity-L1}
\int_0^T|\nu_{s,f,b}|\,ds<\infty
\qquad\text{almost surely},
\end{equation}
and the process
\begin{equation}\label{eq:ac-realization}
\widetilde\zeta_{t,f,b}
:=
\sqrt b\int_0^t\nu_{s,f,b}\,ds,
\qquad 0\le t\le T,
\end{equation}
is centered Gaussian with covariance \(R_{f,b}\). Thus the Gaussian
process with covariance \(R_{f,b}\) has an absolutely continuous, and
therefore continuous finite-variation, realization on \([0,T]\). This
realization is used below.

If \(a\ge0\), then, after redefining \(f(0)\), if necessary,
\(\nu_{f,b}\) has a continuous modification with every H\"older exponent
\[
\delta<\frac{b-1}{2}.
\]
After replacing \(\nu_{f,b}\) by this modification, the realization
\(\widetilde\zeta_{f,b}\) may be chosen in \(C^1([0,T])\), with
\[
\frac{d\widetilde\zeta_{t,f,b}}{dt}
=
\sqrt b\,\nu_{t,f,b},
\qquad 0\le t\le T.
\]
\end{proposition}

We consider covariance at large temporal separation. Since the process
is not stationary, the usual autocovariance terminology for stationary
time series cannot be applied without qualification. We use two related
notions. First, the level-covariance asymptotic describes the behavior
of \(R_{f,b}(s,T+t)\) as the future time \(T+t\) tends to infinity,
while \(s\) and \(t\) remain fixed. Secondly, for a centered
second-order process \(X\), fixed \(0<r\leq v\) and \(0<s\leq t\), we
say that \(X\) satisfies the increment long-range-dependence criterion
with exponent \(\kappa\) if
\[
\lim_{T\to\infty}T^\kappa
\mathbb E\bigl[(X_v-X_r)(X_{T+t}-X_{T+s})\bigr]
\]
exists and is finite and non-zero. This is the increment-covariance
criterion used for Gaussian processes arising from occupation-time
fluctuation limits; see, for instance,
\cite{BGT2007,Samorodnitsky2007,PipirasTaqqu2017}. When the remote
increment covariance is asymptotic to \(CT^{-\kappa}\) with \(C\ne0\),
the summability convention gives long memory if \(0<\kappa\leq1\) and
short memory if \(\kappa>1\); the case \(\kappa=1\) is the usual
borderline. If \(\kappa=0\), the covariance does not decay and is also
non-summable, hence it has long memory under the same convention. Thus, in
the present family, the distinction between long and short memory is
made through the decay rate of remote increment covariances, not
through stationarity of the process itself.

For \(x\geq0\), set
\[
A_f(x):=\int_0^x f(u)(x-u)\,du.
\]
In the parameter ranges used below, this quantity is finite under
\eqref{eq:weight-finiteness}. Indeed, choose \(y>x\). Since \(b\geq0\)
in these ranges,
\[
(y-u)^b\geq(y-x)^b,
\qquad 0\leq u\leq x,
\]
and therefore \eqref{eq:weight-finiteness}, applied at the terminal
time \(y\), implies that \(f\) is integrable on \([0,x]\). Consequently,
\[
0\leq A_f(x)\leq x\int_0^x f(u)\,du<\infty.
\]

\begin{proposition}\label{prop:remote-covariance}
Suppose that \(f\) satisfies \eqref{eq:weight-finiteness}. Then the
following assertions hold.
\begin{itemize}
\item[(i)] For fixed \(s,t\geq0\), if \(b\in[0,1)\), then
\begin{equation}\label{eq:remote-level-rough}
R_{f,b}(s,T+t)
\longrightarrow
\frac{1}{1-b}\int_0^s f(u)(s-u)^b\,du,
\qquad T\to\infty.
\end{equation}
If \(b\in(1,2]\), then
\begin{equation}\label{eq:remote-level-smooth}
\frac{R_{f,b}(s,T+t)}{T^{b-1}}
\longrightarrow
\frac{b}{b-1}\int_0^s f(u)(s-u)\,du,
\qquad T\to\infty.
\end{equation}

\item[(ii)] Let \(b\in(0,1)\cup(1,2]\). For \(0<r\leq v\) and
\(0<s\leq t\),
\begin{equation}\label{eq:remote-increment-limit}
\begin{split}
\lim_{T\to\infty}T^{2-b}
&\mathbb E\left[
(\zeta_{v,f,b}-\zeta_{r,f,b})
(\zeta_{T+t,f,b}-\zeta_{T+s,f,b})
\right]\\
&=
b(t-s)\{A_f(v)-A_f(r)\}\\
&=
b(t-s)\left[
\int_r^v f(u)(v-u)\,du
+
(v-r)\int_0^r f(u)\,du
\right].
\end{split}
\end{equation}
\end{itemize}
\end{proposition}

For \(b=0\), disjoint increments are uncorrelated because
\(\zeta_{f,0}\) is a time-changed Brownian motion. For
\(b\in(0,1)\cup(1,2]\), assume in addition that \(t>s\) and
\(A_f(v)>A_f(r)\); equivalently, the constant on the right-hand side
of \eqref{eq:remote-increment-limit} is non-zero. Under this
non-degeneracy condition,
Proposition~\ref{prop:remote-covariance}(ii) gives
\[
\kappa=2-b,
\]
and the remote increment covariance is asymptotic to a non-zero
constant times \(T^{b-2}\). Under the summability convention described
above, \(0<b<1\) corresponds to short memory because \(2-b>1\),
whereas \(1<b<2\) corresponds to long memory because
\(0<2-b<1\). At \(b=2\), one has \(\kappa=0\), so the remote increment
covariance does not decay. Consequently, it is non-summable and has long
memory under the preceding convention; it is the non-decaying endpoint
of the long-memory regime. The case \(b=1\) is the logarithmic
borderline treated by the limiting covariance. We use
``long-range dependence'' for the normalized increment-covariance
limit and ``long memory'' or ``short memory'' for the summability
classification.

A second-order process is weakly stationary if its mean is constant and
its covariance depends only on the time lag. For a centered second-order
process, we use the covariance-level notion of intrinsic stationarity:
\[
\mathbb E|X_{t+h}-X_t|^2
\]
must depend only on the lag \(h\), whenever both times belong to the
index set.

\begin{proposition}\label{prop:nonstationary}
Assume that \(b\in(0,1)\cup(1,2]\), that \(f\ge0\) satisfies
\eqref{eq:weight-finiteness}, and that \(f>0\) on a set of positive
Lebesgue measure. Then \(\zeta_{f,b}\) is not weakly stationary.

Moreover, it is not intrinsically stationary in either of the following
situations:
\begin{itemize}
\item[(a)] \(1<b\le2\), \(f\) is locally bounded on
\([0,\infty)\), and for some \(s>0\),
\[
\int_0^s f(u)(s-u)^{b-2}\,du>0;
\]

\item[(b)] \(0<b<1\), \(f\) is continuous at zero and \(f(0)>0\).
\end{itemize}
In particular, intrinsic stationarity fails for constant weights; in
the rough regime it also fails for weights continuous and strictly
positive at the origin, and in the smooth regime under condition
\textup{(a)}.
\end{proposition}

\begin{proposition}\label{prop:local-increments}
The following assertions hold.
\begin{itemize}
\item[(i)] Let \(b=0\), let \(f\ge0\) be locally integrable, and suppose
that \(f\) is continuous at \(s\). Then
\begin{equation}\label{eq:local-b0}
\lim_{\varepsilon\downarrow0}
\frac{
\mathbb E|\zeta_{s+\varepsilon,f,0}-\zeta_{s,f,0}|^2
}{\varepsilon}
=
f(s).
\end{equation}

\item[(ii)] Let \(T>0\), let \(0<b<1\), and suppose that \(f\) satisfies
\eqref{eq:weight-finiteness}. If \(1\asymp f\) on \([0,T]\), then for
every fixed \(s<T\),
\begin{equation}\label{eq:local-rough-bounds}
0<
\liminf_{\varepsilon\downarrow0}
\frac{
\mathbb E|\zeta_{s+\varepsilon,f,b}-\zeta_{s,f,b}|^2
}{\varepsilon^{1+b}}
\le
\limsup_{\varepsilon\downarrow0}
\frac{
\mathbb E|\zeta_{s+\varepsilon,f,b}-\zeta_{s,f,b}|^2
}{\varepsilon^{1+b}}
<\infty.
\end{equation}
If \(u^a\asymp f(u)\) on \((0,T]\) for some \(a>-1\), then the same
conclusion holds for every fixed \(s\in(0,T)\).

\item[(iii)] Let \(1<b\le2\), suppose that \(f\) satisfies
\eqref{eq:weight-finiteness}, and assume that \(f\) is locally bounded
in a neighborhood of a fixed \(s>0\). Then
\begin{equation}\label{eq:local-smooth-limit}
\lim_{\varepsilon\downarrow0}
\frac{
\mathbb E|\zeta_{s+\varepsilon,f,b}-\zeta_{s,f,b}|^2
}{\varepsilon^2}
=
2^{b-2}b
\int_0^s f(u)(s-u)^{b-2}\,du.
\end{equation}
\end{itemize}
\end{proposition}

\begin{proposition}\label{prop:small-time-diagonal}
Let \(b\in[0,1)\cup(1,2]\). Suppose that \(f\ge0\) satisfies
\eqref{eq:weight-finiteness} and is continuous at zero. Then
\begin{equation}\label{small-time-diagonal-general}
\lim_{\varepsilon\downarrow0}
\frac{\mathbb E(\zeta_{\varepsilon,f,b}^2)}
     {\varepsilon^{1+b}}
=
\frac{2-2^b}{(1-b)(1+b)}f(0).
\end{equation}

If \(a>-1\) and \(f(u)=u^a\), then
\begin{equation}\label{small-time-diagonal-power}
\lim_{\varepsilon\downarrow0}
\frac{\mathbb E(\zeta_{\varepsilon,f,b}^2)}
     {\varepsilon^{1+b+a}}
=
\frac{2-2^b}{1-b}B(a+1,b+1).
\end{equation}
\end{proposition}

The covariance studied in \cite[Theorem~1]{gonzalez2024} is \(2K_f\)
in the notation below. We include the covariance limit because it
identifies the present normalization as the boundary \(b=1\) of
\(R_{f,b}\).

\begin{proposition}\label{prop:logarithmic-limit}
Let \(f\ge0\) be locally integrable and define, with \(0\log0=0\),
\begin{equation}\label{eq:logarithmic-covariance}
K_f(s,t)
=
\int_0^{s\wedge t}f(u)
\Bigl[
(s+t-2u)\log(s+t-2u)
-(s-u)\log(s-u)
-(t-u)\log(t-u)
\Bigr]\,du.
\end{equation}
Then \(K_f\) is finite and positive definite. For every fixed
\(s,t\ge0\),
\[
R_{f,b}(s,t)\longrightarrow K_f(s,t),
\qquad b\to1,
\]
where \(b\) approaches one through \((0,1)\cup(1,2]\).

Moreover, for either of the families
\[
f_a(u)=u^a,\qquad a>-1,
\]
and
\[
f_a(u)=e^{au},\qquad a\in\mathbb R,
\]
one has
\[
R_{f_a,b}(s,t)\longrightarrow K_1(s,t),
\qquad (a,b)\to(0,1),
\]
with \(b\in(0,1)\cup(1,2]\) and, in the power-weight case,
\(a>-1\) along the approach.
\end{proposition}

Let \(\zeta_f=(\zeta_{t,f})_{t\ge0}\) denote a centered Gaussian
process with covariance \(K_f\). The logarithmic family with covariance
\(2K_f\) is studied in \cite{gonzalez2024}; that work treats variation,
quadratic variation, memory properties, and covariance dependence under
the hypotheses stated there. Thus \(\zeta_f\) is the
\(1/\sqrt2\)-rescaling of the process in that reference, and the
qualitative properties that are invariant under non-zero deterministic
rescaling transfer unchanged.

Sample-path simulations for the processes \(\zeta_{f,b}\) and
\(\zeta_f\) with power and exponential weights are provided in
Section~\ref{sec:supp-paths}; see Figure~\ref{fig:sample-paths}.

\section{Pathwise integration, equations and Euler approximation}
\label{sec:pathwise-euler}

Throughout this section, \(T>0\) is fixed, and \(\mathfrak b\) denotes
the drift coefficient. The covariance exponent continues to be denoted by
\(b\). The equation and Euler results below are restricted to
\(b\in(0,1)\cup(1,2]\); the logarithmic boundary \(b=1\) is not treated
in this section. We also exclude \(b=0\), since the driver is then a
time-changed Brownian motion and the corresponding equation belongs to the
semimartingale SDE setting.

If \(0<b<1\) and \(f\) is bounded on \([0,T]\), then the
covariance identity and comparison with the constant-weight process give
\begin{equation}\label{eq:rough-driver-increment-bound}
\mathbb E|\zeta_{t,f,b}-\zeta_{s,f,b}|^2
\le C_{T,f,b}|t-s|^{1+b},
\qquad 0\le s<t\le T.
\end{equation}
Consequently, Gaussian moment estimates and Kolmogorov's continuity
theorem yield a modification with H\"older continuous paths of every
order
\[
\rho<\frac{1+b}{2}.
\]
In particular, one may choose \(\rho>1/2\), and integration with respect
to \(\zeta_{f,b}\) is understood in the Young sense
\cite{Young1936,FrizVictoir2010,FrizHairer2014}.

If \(1<b\le2\), we use the absolutely continuous realization provided
by Proposition~\ref{prop:gaussian-velocity}. Integration is then
interpreted through
\[
\int_0^t G_s\,d\widetilde\zeta_{s,f,b}
=
\sqrt b\int_0^t G_s\nu_{s,f,b}\,ds.
\]
When \(f\) is bounded on \([0,T]\), this realization may be chosen in
\(C^1([0,T])\).

The detailed proofs for this section are collected in
Section~\ref{sec:supp-euler-proofs}.
Section~\ref{sec:supp-young-aux} contains the auxiliary Young estimates,
Section~\ref{sec:supp-proof-thm2} proves
Theorem~\ref{thm:young-equation-wsfbm},
Section~\ref{sec:supp-proof-thm3} proves
Theorem~\ref{thm:finite-variation-equation-wsfbm},
Section~\ref{sec:supp-euler-aux} contains the auxiliary Euler estimates, and
Section~\ref{sec:supp-proof-thm4} proves Theorem~\ref{thm:euler-wsfbm}.

\subsection{Existence, uniqueness and stability}

For a prescribed initial value \(x_0\in\mathbb R\), we consider
\begin{equation}\label{eq:general-wsfbm-equation}
X_t
=
x_0
+
\int_0^t\mathfrak b(X_r,r)\,dr
+
\int_0^t\sigma(X_r,r)\,dZ_r.
\end{equation}
Assume the following regularity condition.

\begin{description}
\item[(H)]
The functions \(\mathfrak b,\sigma\) are bounded and continuous on
\(\mathbb R\times[0,T]\). The derivatives
\[
\partial_x\mathfrak b,\quad
\partial_t\mathfrak b,\quad
\partial_x\sigma,\quad
\partial_{xx}\sigma,\quad
\partial_t\sigma,\quad
\partial_{xt}\sigma
\]
exist, are continuous, and are bounded on
\(\mathbb R\times[0,T]\).
\end{description}

The boundedness of \(\mathfrak b\) and \(\sigma\) is used in the
uniform one-step estimates underlying the Euler analysis.

\begin{thm}\label{thm:young-equation-wsfbm}
Let \(0<b<1\), let \(f\) be an admissible weight that is bounded on
\([0,T]\), and choose
\[
\frac12<\rho<\frac{1+b}{2}.
\]
For the H\"older realization
\[
Z=\zeta_{f,b}\in C^\rho([0,T]),
\]
condition \textup{(H)} implies that
\eqref{eq:general-wsfbm-equation} has a unique solution
\(X\in C^\rho([0,T])\), almost surely.

More precisely, for every \(M>0\), there exists a finite constant
\(C_{T,M}\) such that, whenever \(Z,\widetilde Z\in C^\rho([0,T])\)
satisfy
\[
Z_0=\widetilde Z_0=0
\]
and
\[
|x_0|\vee |y_0|
\vee [Z]_{\rho;[0,T]}
\vee [\widetilde Z]_{\rho;[0,T]}
\le M,
\]
the corresponding solutions \(X\) and \(Y\) satisfy
\begin{equation}\label{eq:young-solution-stability}
\|X-Y\|_{\rho;[0,T]}
\le
C_{T,M}
\left(
|x_0-y_0|
+
\|Z-\widetilde Z\|_{\rho;[0,T]}
\right).
\end{equation}
\end{thm}

\begin{thm}\label{thm:finite-variation-equation-wsfbm}
Let \(1<b\le2\), and suppose that
\[
0\le f(u)\le cu^a,
\qquad 0<u\le T,
\]
for some \(c>0\) and \(a>-1\). Select the absolutely continuous
realization
\[
Z_t=\sqrt b\int_0^t\nu_r\,dr,
\qquad 0\le t\le T,
\]
provided by Proposition~\ref{prop:gaussian-velocity}. Under
\textup{(H)}, for almost every sample path,
\eqref{eq:general-wsfbm-equation} has a unique absolutely continuous
solution. Equivalently, this solution is the unique Carath\'eodory
solution of
\begin{equation}\label{eq:random-ode-fv}
\dot X_t
=
\mathfrak b(X_t,t)
+
\sqrt b\,\sigma(X_t,t)\nu_t
\quad\text{for a.e. }t\in[0,T],
\qquad
X_0=x_0.
\end{equation}
\end{thm}

\subsection{Euler approximation}

For a deterministic partition
\[
\pi=\{0=t_0<\cdots<t_n=T\},
\]
put
\[
h:=|\pi|,
\qquad
\Delta_k:=t_{k+1}-t_k,
\qquad
Z_{t_k,t_{k+1}}:=Z_{t_{k+1}}-Z_{t_k}.
\]
Define
\begin{equation}\label{eq:euler-wsfbm}
X_{t_{k+1}}^\pi
=
X_{t_k}^\pi
+
\mathfrak b(X_{t_k}^\pi,t_k)\Delta_k
+
\sigma(X_{t_k}^\pi,t_k)Z_{t_k,t_{k+1}},
\qquad
X_0^\pi=x_0,
\end{equation}
and, for \(t\in[t_k,t_{k+1}]\),
\begin{equation}\label{eq:euler-interpolation-wsfbm}
X_t^\pi
=
X_{t_k}^\pi
+
\mathfrak b(X_{t_k}^\pi,t_k)(t-t_k)
+
\sigma(X_{t_k}^\pi,t_k)Z_{t_k,t}.
\end{equation}

\begin{thm}\label{thm:euler-wsfbm}
Assume \textup{(H)} and let \(p\ge1\).

\begin{itemize}
\item[(i)] Let \(0<b<1\), let \(f\) be an admissible weight that is bounded on \([0,T]\), set
\(Z=\zeta_{f,b}\), and choose
\[
\frac12<\rho<\frac{1+b}{2}.
\]
Then, almost surely,
\begin{equation}\label{eq:euler-young-pathwise-final}
\sup_{0\le t\le T}|X_t-X_t^\pi|
\le
C_T
\exp\!\left\{
C_T\bigl(1+[Z]_{\rho;[0,T]}\bigr)^{1/\rho}
\right\}
\bigl(1+[Z]_{\rho;[0,T]}\bigr)^2
h^{2\rho-1}.
\end{equation}
Moreover,
\begin{equation}\label{eq:euler-lp-young}
\left\|
\sup_{0\le t\le T}|X_t-X_t^\pi|
\right\|_{L^p}
\le
C_{p,T,\rho}h^{2\rho-1}.
\end{equation}

\item[(ii)] Let \(1<b\le2\), and suppose that
\[
0\le f(u)\le cu^a,
\qquad 0<u\le T,
\]
for some \(c>0\) and \(a>-1\). Use the realization
\[
Z_t=\sqrt b\int_0^t\nu_r\,dr
\]
provided by Proposition~\ref{prop:gaussian-velocity}, and define
\[
A_t
:=
t+\sqrt b\int_0^t|\nu_r|\,dr,
\qquad
\omega_A(h)
:=
\sup_{\substack{s,t\in[0,T]\\|t-s|\le h}}
|A_t-A_s|.
\]
Then, almost surely,
\begin{equation}\label{eq:euler-fv-modulus}
\sup_{0\le t\le T}|X_t-X_t^\pi|
\le
C_Te^{C_TA_T}A_T\omega_A(h).
\end{equation}
Consequently, the Euler interpolation converges uniformly to \(X\)
as \(h\downarrow0\).

If \(a\ge0\), so that \(\nu\) admits a continuous Gaussian
modification, then
\begin{equation}\label{eq:euler-lp-fv}
\left\|
\sup_{0\le t\le T}|X_t-X_t^\pi|
\right\|_{L^p}
\le
C_{p,T}h.
\end{equation}
\end{itemize}

Here and below, deterministic constants may depend on
\(T,p,\rho,b,f\) and the coefficient bounds in \textup{(H)}, but not
on the partition \(\pi\).
\end{thm}

\clearpage
\section*{Supplementary Material}
\addcontentsline{toc}{section}{Supplementary Material}

\setcounter{secnumdepth}{3}
\setcounter{section}{0}
\setcounter{subsection}{0}
\setcounter{subsubsection}{0}
\setcounter{equation}{0}
\setcounter{figure}{0}
\setcounter{table}{0}
\setcounter{lem}{0}
\setcounter{proposition}{0}
\setcounter{example}{0}

\renewcommand{\thesubsection}{S.\arabic{subsection}}
\renewcommand{\theequation}{S.\arabic{equation}}
\renewcommand{\thefigure}{S.\arabic{figure}}
\renewcommand{\thetable}{S.\arabic{table}}
\renewcommand{\thelem}{S.\arabic{lem}}
\renewcommand{\theproposition}{S.\arabic{proposition}}
\renewcommand{\theexample}{S.\arabic{example}}

\providecommand*{\theHsubsection}{}
\providecommand*{\theHsubsubsection}{}
\providecommand*{\theHequation}{}
\providecommand*{\theHfigure}{}
\providecommand*{\theHtable}{}
\providecommand*{\theHlem}{}
\providecommand*{\theHproposition}{}
\providecommand*{\theHexample}{}
\renewcommand*{\theHsubsection}{supp.\arabic{subsection}}
\renewcommand*{\theHsubsubsection}{supp.\arabic{subsection}.\arabic{subsubsection}}
\renewcommand*{\theHequation}{supp.\arabic{equation}}
\renewcommand*{\theHfigure}{supp.\arabic{figure}}
\renewcommand*{\theHtable}{supp.\arabic{table}}
\renewcommand*{\theHlem}{supp.\arabic{lem}}
\renewcommand*{\theHproposition}{supp.\arabic{proposition}}
\renewcommand*{\theHexample}{supp.\arabic{example}}

\section*{Supplementary contents}

Section~\ref{sec:supp-path-proofs} contains the detailed proofs of the sample-path results stated in Propositions~\ref{prop:rough-bounded-weight}--\ref{prop:logarithmic-limit}. More precisely, Section~\ref{sec:supp-infinite-variation-lemma} gives the auxiliary infinite-variation lemma, and Sections~\ref{sec:supp-proof-prop1}--\ref{sec:supp-proof-prop8} contain the eight proofs. Section~\ref{sec:supp-euler-proofs} contains the pathwise-equation and Euler arguments: Sections~\ref{sec:supp-young-aux} and \ref{sec:supp-euler-aux} collect the auxiliary estimates, while Sections~\ref{sec:supp-proof-thm2}, \ref{sec:supp-proof-thm3}, and \ref{sec:supp-proof-thm4} prove Theorems~\ref{thm:young-equation-wsfbm}, \ref{thm:finite-variation-equation-wsfbm}, and \ref{thm:euler-wsfbm}, respectively. Section~\ref{sec:euclidean-covariances} gives related covariance constructions on Euclidean index sets. Section~\ref{sec:supp-numerical} gives closed covariance formulas and numerical comparisons, Section~\ref{sec:supp-inference} reports likelihood and prediction calculations for simulated data, and Section~\ref{sec:supp-paths} presents sample paths.

Throughout, \(\zeta_{f,b}\) denotes the centered Gaussian process with covariance
\begin{equation}\label{supp-weighted-covariance}
R_{f,b}(s,t)
=
\frac{1}{1-b}
\int_0^{s\wedge t}
f(r)
\bigl\{
(s-r)^b+(t-r)^b-(s+t-2r)^b
\bigr\}\,dr,
\end{equation}
where \(b\in[0,1)\cup(1,2]\) and \(f\) satisfies the assumptions stated in Section~\ref{sec:covariance-construction}. At \(b=1\), we use the centered Gaussian process \(\zeta_f\) with covariance \(K_f\), where the corresponding covariance \(2K_f\) was studied in \cite{gonzalez2024}:
\begin{equation}\label{supp-log-covariance}
K_f(s,t)
=
\int_0^{s\wedge t}
f(r)
\bigl[
(s+t-2r)\log(s+t-2r)
-(s-r)\log(s-r)
-(t-r)\log(t-r)
\bigr]\,dr,
\end{equation}
with \(0\log0=0\).

The calculations use the two weight families, defined for \(u>0\),
\[
\mathcal C_1
=
\bigl\{f_a^{(1)}(u)=u^a:\ a>-1\bigr\},
\qquad
\mathcal C_2
=
\bigl\{f_a^{(2)}(u)=e^{au}:\ a\in\mathbb R\bigr\}.
\]
Both families are admissible for every \(b>-1\). Any non-negative value may be assigned at \(u=0\), without changing the covariance integrals. We write
\[
R_{a,b}^{(i)}:=R_{f_a^{(i)},b},
\qquad
K_a^{(i)}:=K_{f_a^{(i)}},
\qquad i=1,2.
\]
\subsection{Detailed proofs of the sample-path results}\label{sec:supp-path-proofs}

All references below point directly to the corresponding propositions, theorems, and equations in the first part of this document.

\subsubsection{Auxiliary infinite-variation lemma}\label{sec:supp-infinite-variation-lemma}

\begin{lem}\label{lem:infinite-variation-gaussian}
Let \(a<d\), and let \(X=(X_t)_{t\in[a,d]}\) be a centered continuous
Gaussian process. Assume that there exist constants \(c>0\) and
\(\gamma\in(0,2)\) such that
\[
\mathbb E|X_t-X_s|^2
\ge c|t-s|^\gamma,
\qquad a\le s<t\le d.
\]
Then
\[
V_{[a,d]}(X)=\infty
\qquad\text{almost surely}.
\]
\end{lem}

\begin{proof}
For \(n\ge1\), let
\[
t_k^{(n)}=a+\frac{k(d-a)}{n},
\qquad k=0,\ldots,n,
\]
and set
\[
V_n
=
\sum_{k=0}^{n-1}
\left|X_{t_{k+1}^{(n)}}-X_{t_k^{(n)}}\right|.
\]
Each increment is a centered Gaussian random variable. Since a centered
normal random variable with variance \(\sigma^2\) has mean absolute value
\(\sqrt{2/\pi}\,\sigma\), the assumed lower bound gives
\[
\mathbb EV_n
\ge
\sqrt{\frac{2c}{\pi}}\,
n\left(\frac{d-a}{n}\right)^{\gamma/2}
=
\sqrt{\frac{2c}{\pi}}\,
(d-a)^{\gamma/2}n^{1-\gamma/2}.
\]
Since \(\gamma<2\),
\begin{equation}\label{eq:expected-discrete-variation}
\mathbb EV_n\longrightarrow\infty.
\end{equation}

Regard \(X\) as a centered Gaussian random element of the separable
Banach space \(E=C([a,d])\), equipped with the uniform norm. Define
\[
V(x):=V_{[a,d]}(x),
\qquad x\in E.
\]
Let \(\mathcal P_{\mathbb Q}\) denote the countable collection of finite
partitions of \([a,d]\) whose interior points belong to
\(\mathbb Q\cap(a,d)\). Uniform continuity gives
\[
V(x)
=
\sup_{\pi\in\mathcal P_{\mathbb Q}}
\sum_{[u,v]\in\pi}|x(v)-x(u)|.
\]
For each fixed partition, the corresponding variation sum is continuous
on \(E\). Hence \(V\) is Borel measurable. Furthermore,
\[
V(x+y)\le V(x)+V(y),
\qquad
V(\lambda x)=|\lambda|V(x),
\]
and the set \(\{x\in E:V(x)<\infty\}\) is a vector subspace of \(E\).
Thus \(V\) is a measurable pseudo-seminorm in the sense of
\cite[Section~1.2]{Fernique1975}.

Suppose that
\[
\mathbb P\{V(X)<\infty\}>0.
\]
The set \(\{x\in E:V(x)<\infty\}\) is a measurable vector subspace.
The Gaussian zero--one law for measurable linear subspaces therefore
implies that it has probability one; see, for example,
\cite{Bogachev1998}. Fernique's integrability theorem for measurable
pseudo-seminorms \cite[Theorem~1.3.2, p.~11]{Fernique1975} then yields
an \(\alpha>0\) such that
\[
\mathbb E\exp\{\alpha V(X)^2\}<\infty.
\]
In particular, \(\mathbb EV(X)<\infty\). Since \(V_n\le V(X)\) almost
surely for every \(n\),
\[
\mathbb EV_n\le\mathbb EV(X)<\infty,
\]
contradicting \eqref{eq:expected-discrete-variation}. Therefore
\[
\mathbb P\{V(X)<\infty\}=0,
\]
which proves the assertion.
\end{proof}

\subsubsection{Proof of Proposition~\ref{prop:rough-bounded-weight}}\label{sec:supp-proof-prop1}
\begin{proof}
Since \(1\asymp f\) on \([0,T]\), there exist constants
\(0<m\le M<\infty\) such that
\[
m\le f(u)\le M,
\qquad 0\le u\le T.
\]

For \(0\le s<t\le T\), using the covariance formula and the diagonal
identity
\[
R_{f,b}(t,t)
=
\frac{2-2^b}{1-b}
\int_0^t f(u)(t-u)^b\,du,
\]
we obtain
\begin{equation}\label{eq:increment-variance-decomposition}
\begin{split}
\mathbb E(\zeta_{t,f,b}-\zeta_{s,f,b})^2
&=
\frac{2-2^b}{1-b}
\int_s^t f(u)(t-u)^b\,du\\
&\quad+
\frac{1}{1-b}
\int_0^s f(u)
\Bigl[
2(t+s-2u)^b
-2^b(t-u)^b
-2^b(s-u)^b
\Bigr]\,du.
\end{split}
\end{equation}

Both integrands in \eqref{eq:increment-variance-decomposition} are
non-negative. Indeed, the first one is non-negative because
\(0<b<1\). For the second one, put
\[
x=t-u,\qquad y=s-u.
\]
Since \(x,y\ge0\) and \(r\mapsto r^b\) is concave,
\[
\left(\frac{x+y}{2}\right)^b
\ge
\frac{x^b+y^b}{2}.
\]
Equivalently,
\[
2(x+y)^b-2^b x^b-2^b y^b\ge0,
\]
which is precisely the bracket in
\eqref{eq:increment-variance-decomposition}.

Because the right-hand side of
\eqref{eq:increment-variance-decomposition} is linear in \(f\) and has
non-negative integrands, comparison of the weights gives
\[
m\,
\mathbb E(\zeta_{t,1,b}-\zeta_{s,1,b})^2
\le
\mathbb E(\zeta_{t,f,b}-\zeta_{s,f,b})^2
\le
M\,
\mathbb E(\zeta_{t,1,b}-\zeta_{s,1,b})^2.
\]
Combining this inequality with
\eqref{eq:constant-increment-bounds}, we obtain
\begin{equation}\label{eq:bounded-weight-increment-bounds}
C_1|t-s|^{1+b}
\le
\mathbb E(\zeta_{t,f,b}-\zeta_{s,f,b})^2
\le
C_2|t-s|^{1+b},
\qquad 0\le s<t\le T,
\end{equation}
where one may take
\[
C_1
=
m\,\frac{2-2^b}{(1-b)(1+b)},
\qquad
C_2
=
M\,\frac{1}{(1-b)(1+b)}.
\]

Let \(p\ge2\). Since the increments are centered Gaussian random
variables, there exists a constant \(c_p>0\) such that
\[
\mathbb E
|\zeta_{t,f,b}-\zeta_{s,f,b}|^p
=
c_p
\left(
\mathbb E
|\zeta_{t,f,b}-\zeta_{s,f,b}|^2
\right)^{p/2}.
\]
Consequently, by the upper bound in
\eqref{eq:bounded-weight-increment-bounds},
\[
\mathbb E
|\zeta_{t,f,b}-\zeta_{s,f,b}|^p
\le
c_p C_2^{p/2}|t-s|^{p(1+b)/2}.
\]
Kolmogorov's continuity theorem therefore yields a modification whose
sample paths are Hölder continuous of every order
\[
\delta<\frac{1+b}{2}.
\]
We work with this modification and retain the notation
\(\zeta_{f,b}\).

Fix
\[
\delta\in\left(\frac12,\frac{1+b}{2}\right).
\]
On an event of probability one,
\([\zeta_{f,b}]_{\delta;[0,T]}<\infty\). For every finite partition
\(\pi\) of \([0,T]\),
\begin{align*}
\sum_{[u,v]\in\pi}
|\zeta_{v,f,b}-\zeta_{u,f,b}|^2
&\le
[\zeta_{f,b}]_{\delta;[0,T]}^2
\sum_{[u,v]\in\pi}|v-u|^{2\delta}\\
&\le
[\zeta_{f,b}]_{\delta;[0,T]}^2
|\pi|^{2\delta-1}
\sum_{[u,v]\in\pi}|v-u|=
T[\zeta_{f,b}]_{\delta;[0,T]}^2
|\pi|^{2\delta-1}.
\end{align*}
Since \(2\delta-1>0\), the right-hand side converges to zero whenever
\(|\pi|\to0\). Thus, on the same probability-one event, the process has
zero quadratic variation along every deterministic sequence of
partitions whose mesh tends to zero.

The lower bound in
\eqref{eq:bounded-weight-increment-bounds} and
Lemma~\ref{lem:infinite-variation-gaussian}, applied with
\(
\gamma=1+b\in(1,2),
\)
give
\(
V_{[0,T]}(\zeta_{f,b})=\infty
,\,\text{almost surely}.
\)

Finally, suppose that this continuous modification were a
semimartingale. Its quadratic variation along deterministic partitions
would converge in probability to its semimartingale quadratic
variation. Since the preceding pathwise estimate gives convergence to
zero, its quadratic variation would be identically zero. A continuous
semimartingale with zero quadratic variation has a constant local
martingale part and therefore has finite variation. This contradicts
the almost sure infinite total variation established above. Hence
\(\zeta_{f,b}\) is not a semimartingale.
\end{proof}

\subsubsection{Proof of Proposition~\ref{prop:rough-power-comparable}}\label{sec:supp-proof-prop2}
\begin{proof}
Let
\[
f_\alpha(u):=u^\alpha,\qquad u>0.
\]
The value assigned to \(f_\alpha\) at \(u=0\) is immaterial. Since
\(\alpha>-1\) and \(b>-1\), the reference weight is admissible:
\[
\int_0^t u^\alpha(t-u)^b\,du
=
t^{\alpha+b+1}B(\alpha+1,b+1)
<\infty,
\qquad t>0.
\]
We first establish uniform increment bounds for the Gaussian process
associated with \(f_\alpha\).

For \(0\le s<t\le T\), the increment-variance decomposition
\eqref{eq:increment-variance-decomposition} gives
\begin{equation}\label{eq:power-reference-variance}
\begin{split}
\mathbb E
\bigl(\zeta_{t,f_\alpha,b}-\zeta_{s,f_\alpha,b}\bigr)^2
&=
\frac{2-2^b}{1-b}
\int_s^t u^\alpha(t-u)^b\,du\\
&\quad+
\frac{1}{1-b}
\int_0^s u^\alpha
\Bigl[
2(t+s-2u)^b
-2^b(t-u)^b
-2^b(s-u)^b
\Bigr]\,du.
\end{split}
\end{equation}

We rewrite the second term. For fixed \(0\le u<s\), put
\[
F_u(r,r'):=(r+r'-2u)^b.
\]
Since
\[
\frac{\partial^2 F_u}{\partial r\,\partial r'}(r,r')
=
b(b-1)(r+r'-2u)^{b-2},
\]
rectangular integration over \([s,t]^2\) yields
\begin{align*}
&2^b(t-u)^b
-2(t+s-2u)^b
+2^b(s-u)^b\\
&\qquad=
b(b-1)
\int_s^t\int_s^t
(r+r'-2u)^{b-2}\,dr'\,dr.
\end{align*}
Therefore
\[
\frac{
2(t+s-2u)^b
-2^b(t-u)^b
-2^b(s-u)^b
}{1-b}
=
b\int_s^t\int_s^t
(r+r'-2u)^{b-2}\,dr'\,dr.
\]
Multiplying by \(u^\alpha\), integrating over \(u\in[0,s)\), and using
Tonelli's theorem gives
\begin{equation}\label{eq:power-reference-triple-integral}
\begin{split}
&\frac{1}{1-b}
\int_0^s u^\alpha
\Bigl[
2(t+s-2u)^b
-2^b(t-u)^b
-2^b(s-u)^b
\Bigr]\,du\\
&\qquad=
b\int_s^t\int_s^t\int_0^s
u^\alpha(r+r'-2u)^{b-2}\,du\,dr'\,dr.
\end{split}
\end{equation}
The endpoint \(u=s\) does not affect either side. Notice also that the
integrand on the right-hand side is non-negative.

Write
\[
h:=t-s,\qquad r=s+x,\qquad r'=s+y,\qquad u=s-v.
\]
Then the right-hand side of
\eqref{eq:power-reference-triple-integral} becomes
\begin{equation}\label{eq:power-singular-integral}
b\int_0^h\int_0^h\int_0^s
(s-v)^\alpha(x+y+2v)^{b-2}\,dv\,dy\,dx.
\end{equation}

Fix \(\varepsilon\in(0,T)\) and assume that
\[
\varepsilon\le s<t\le T.
\]
We split the integral with respect to \(v\) into the intervals
\([0,s/2]\) and \([s/2,s]\).

On \(0\le v\le s/2\), one has
\[
\frac{\varepsilon}{2}\le s-v\le T.
\]
Hence
\[
(s-v)^\alpha\le C_{\varepsilon,T,\alpha}.
\]
Consequently,
\begin{align}
&\int_0^h\int_0^h\int_0^{s/2}
(s-v)^\alpha(x+y+2v)^{b-2}\,dv\,dy\,dx
\notag\\
&\quad\le
C_{\varepsilon,T,\alpha}
\int_0^h\int_0^h\int_0^\infty
(x+y+2v)^{b-2}\,dv\,dy\,dx
\notag\\
&\quad=
\frac{C_{\varepsilon,T,\alpha}}{2(1-b)}
\int_0^h\int_0^h
(x+y)^{b-1}\,dy\,dx
\notag\\
&\quad=
\frac{C_{\varepsilon,T,\alpha}}{2b(1-b)}
\int_0^h
\bigl[(x+h)^b-x^b\bigr]\,dx
\notag\\
&\quad=
\frac{C_{\varepsilon,T,\alpha}(2^{b+1}-2)}
     {2b(1-b)(b+1)}
h^{b+1}.
\label{eq:power-near-diagonal-bound}
\end{align}
The identity obtained after integrating with respect to \(v\) holds
for \(x+y>0\), and therefore almost everywhere on \([0,h]^2\). The
resulting function \((x+y)^{b-1}\) is integrable on this square because
\(b-1>-1\).

On \(s/2\le v\le s\), one has
\[
x+y+2v\ge s\ge\varepsilon.
\]
Since \(b-2<0\),
\[
(x+y+2v)^{b-2}\le\varepsilon^{b-2}.
\]
Moreover, since \(\alpha>-1\),
\[
\int_{s/2}^s(s-v)^\alpha\,dv
=
\int_0^{s/2}w^\alpha\,dw
=
\frac{(s/2)^{\alpha+1}}{\alpha+1}
\le C_{T,\alpha}.
\]
It follows that
\begin{align}
&\int_0^h\int_0^h\int_{s/2}^s
(s-v)^\alpha(x+y+2v)^{b-2}\,dv\,dy\,dx
\notag\\
&\qquad\le
C_{\varepsilon,T,\alpha}h^2
\le
C_{\varepsilon,T,\alpha}T^{1-b}h^{b+1},
\label{eq:power-origin-bound}
\end{align}
because \(0<h\le T\) and \(1-b>0\).

We now estimate the first term in
\eqref{eq:power-reference-variance}. Since \(u^\alpha\) is bounded
above and bounded away from zero on \([\varepsilon,T]\), there exist
constants
\[
0<m_{\varepsilon,T,\alpha}
\le M_{\varepsilon,T,\alpha}<\infty
\]
such that
\[
m_{\varepsilon,T,\alpha}
\le u^\alpha
\le M_{\varepsilon,T,\alpha},
\qquad \varepsilon\le u\le T.
\]
Therefore
\begin{equation}\label{eq:first-term-power-bounds}
m_{\varepsilon,T,\alpha}
\int_s^t(t-u)^b\,du
\le
\int_s^t u^\alpha(t-u)^b\,du
\le
M_{\varepsilon,T,\alpha}
\int_s^t(t-u)^b\,du,
\end{equation}
where
\begin{equation}\label{eq:first-term-power-exact}
\int_s^t(t-u)^b\,du
=
\frac{(t-s)^{b+1}}{b+1}.
\end{equation}

The second term in
\eqref{eq:power-reference-variance} is non-negative by
\eqref{eq:power-reference-triple-integral}. Hence the lower bound in
\eqref{eq:first-term-power-bounds} gives a lower bound of order
\(h^{1+b}\). The upper bound in
\eqref{eq:first-term-power-bounds}, together with
\eqref{eq:power-near-diagonal-bound} and
\eqref{eq:power-origin-bound}, gives the corresponding upper bound.
Thus there exist constants
\[
0<C_{1,\varepsilon,T,\alpha,b}
\le C_{2,\varepsilon,T,\alpha,b}<\infty
\]
such that
\begin{equation}\label{eq:power-reference-increment-bounds}
C_{1,\varepsilon,T,\alpha,b}|t-s|^{1+b}
\le
\mathbb E
\bigl|\zeta_{t,f_\alpha,b}-\zeta_{s,f_\alpha,b}\bigr|^2
\le
C_{2,\varepsilon,T,\alpha,b}|t-s|^{1+b},
\qquad
\varepsilon\le s<t\le T.
\end{equation}

By assumption, there exist constants \(0<c_-\le c_+<\infty\) such that
\[
c_-u^\alpha\le f(u)\le c_+u^\alpha,
\qquad 0<u\le T.
\]
In the increment-variance decomposition
\eqref{eq:increment-variance-decomposition}, the coefficient and
integrand in the first term are non-negative, and the kernel in the
second term is non-negative by
\eqref{eq:power-reference-triple-integral}. Since the decomposition is
linear in the weight, comparison of the weights gives
\[
c_-\,
\mathbb E
\bigl|\zeta_{t,f_\alpha,b}-\zeta_{s,f_\alpha,b}\bigr|^2
\le
\mathbb E
\bigl|\zeta_{t,f,b}-\zeta_{s,f,b}\bigr|^2
\le
c_+\,
\mathbb E
\bigl|\zeta_{t,f_\alpha,b}-\zeta_{s,f_\alpha,b}\bigr|^2.
\]
It follows that, for every \(\varepsilon\in(0,T)\), there are constants
\(0<C_{1,\varepsilon}\le C_{2,\varepsilon}<\infty\) such that
\begin{equation}\label{eq:power-comparable-increment-bounds}
C_{1,\varepsilon}|t-s|^{1+b}
\le
\mathbb E
\bigl|\zeta_{t,f,b}-\zeta_{s,f,b}\bigr|^2
\le
C_{2,\varepsilon}|t-s|^{1+b},
\qquad
\varepsilon\le s<t\le T.
\end{equation}

Fix
\[
\delta\in\left(\frac12,\frac{1+b}{2}\right).
\]
For \(n\ge1\), set
\[
I_n:=\left[\frac{T}{n+1},T\right].
\]
Since the increments are centered Gaussian random variables, for every
\(p\ge2\) there is a constant \(c_p>0\) such that
\[
\mathbb E
\bigl|\zeta_{t,f,b}-\zeta_{s,f,b}\bigr|^p
=
c_p
\left(
\mathbb E
\bigl|\zeta_{t,f,b}-\zeta_{s,f,b}\bigr|^2
\right)^{p/2}.
\]
Choose \(p\) sufficiently large that
\[
\delta<
\frac{1+b}{2}-\frac1p.
\]
Using the upper bound in
\eqref{eq:power-comparable-increment-bounds} with
\(\varepsilon=T/(n+1)\), Kolmogorov's continuity theorem gives a
modification \(\zeta^{(n)}_{f,b}\) on \(I_n\) whose paths are
\(\delta\)-Hölder continuous.

These modifications can be chosen consistently. Indeed, if \(m>n\),
then \(I_n\subset I_m\), and both
\(\zeta^{(n)}_{f,b}\) and \(\zeta^{(m)}_{f,b}\) are modifications of
the same process on \(I_n\). They therefore agree almost surely at
every point of a fixed countable dense subset of \(I_n\). Since both
versions are continuous, they are indistinguishable on \(I_n\).
Intersecting the corresponding probability-one events over the
countable collection of pairs \(m>n\), we may patch the versions to
obtain one modification, still denoted by \(\zeta_{f,b}\), which is
continuous on \((0,T]\) and \(\delta\)-Hölder continuous on every
compact subinterval of \((0,T]\).

Fix a rational number \(q\in(0,T)\). On an event of probability one,
\[
[\zeta_{f,b}]_{\delta;[q,T]}<\infty.
\]
For every finite partition \(\pi\) of \([q,T]\),
\begin{align*}
\sum_{[u,v]\in\pi}
\bigl|\zeta_{v,f,b}-\zeta_{u,f,b}\bigr|^2
&\le
[\zeta_{f,b}]_{\delta;[q,T]}^2
\sum_{[u,v]\in\pi}|v-u|^{2\delta}\\
&\le
(T-q)
[\zeta_{f,b}]_{\delta;[q,T]}^2
|\pi|^{2\delta-1}.
\end{align*}
Since \(2\delta-1>0\), the right-hand side tends to zero whenever
\(|\pi|\to0\). Thus the quadratic variation on \([q,T]\) is zero
almost surely.

On the other hand, the lower bound in
\eqref{eq:power-comparable-increment-bounds}, applied with
\(\varepsilon=q\), and
Lemma~\ref{lem:infinite-variation-gaussian}, with
\[
\gamma=1+b\in(1,2),
\]
give
\[
V_{[q,T]}(\zeta_{f,b})=\infty
\qquad\text{almost surely}.
\]

Intersecting these probability-one events over
\(q\in\mathbb Q\cap(0,T)\), we obtain a single event of probability one
on which both conclusions hold for every rational \(q\in(0,T)\).

Let now \(\varepsilon\in(0,T)\) be arbitrary. Choose rational numbers
\(q_-,q_+\) such that
\[
0<q_-<\varepsilon<q_+<T.
\]
Since the path is \(\delta\)-Hölder continuous on \([q_-,T]\), the same
partition estimate shows that its quadratic variation on
\([\varepsilon,T]\) is zero. Moreover,
\[
V_{[\varepsilon,T]}(\zeta_{f,b})
\ge
V_{[q_+,T]}(\zeta_{f,b})
=
\infty.
\]
Therefore, on one event of probability one, both conclusions hold for
every \(\varepsilon\in(0,T)\).

Finally, fix \(\varepsilon\in(0,T)\) and suppose that the restriction of
this continuous modification to \([\varepsilon,T]\) were a
semimartingale with respect to its completed natural filtration. Its
quadratic variation along deterministic partitions would converge in
probability to its semimartingale quadratic variation. Since the
preceding pathwise estimate gives convergence to zero, its quadratic
variation would be identically zero.

Writing the continuous semimartingale as
\[
\zeta_{f,b}=M+A,
\]
where \(M\) is a continuous local martingale and \(A\) has finite
variation, one has
\[
[\zeta_{f,b}]=[M].
\]
Thus \([M]=0\), so \(M\) is constant and \(\zeta_{f,b}\) has finite
variation on \([\varepsilon,T]\). This contradicts
\[
V_{[\varepsilon,T]}(\zeta_{f,b})=\infty
\qquad\text{almost surely}.
\]
Hence the restriction is not a semimartingale. No continuity at the
origin is asserted.
\end{proof}

\subsubsection{Proof of Proposition~\ref{prop:gaussian-velocity}}\label{sec:supp-proof-prop3}
\begin{proof}
For \(1<b<2\),
\begin{equation}\label{eq:laplace-power}
z^{b-2}
=
\frac{1}{\Gamma(2-b)}
\int_0^\infty e^{-\lambda z}\lambda^{1-b}\,d\lambda,
\qquad z>0,
\end{equation}
whereas for \(b=2\) the function \(z^{b-2}\) is identically one.

When \(1<b<2\), set the value at \((s,t)=(u,u)\) equal to zero.
Formula \eqref{eq:laplace-power} shows that, for every \(u\ge0\),
\[
(s,t)\longmapsto
\mathbf 1_{(u,\infty)}(s)
\mathbf 1_{(u,\infty)}(t)
(s+t-2u)^{b-2}
\]
is positive definite. Indeed, it is a positive mixture of the rank-one
kernels
\[
(s,t)\longmapsto
\mathbf 1_{(u,\infty)}(s)e^{-\lambda(s-u)}
\mathbf 1_{(u,\infty)}(t)e^{-\lambda(t-u)}.
\]
For \(b=2\), the corresponding elementary kernel is
\[
(s,t)\longmapsto
\mathbf 1_{[u,\infty)}(s)\mathbf 1_{[u,\infty)}(t),
\]
which has rank one. These endpoint conventions do not change any
Lebesgue integral. Integration with respect to the non-negative measure
\(f(u)\,du\) therefore proves that \(C_{f,b}\) is positive definite.

Its diagonal is finite for every \(0<t\le T\), because
\begin{equation}\label{eq:velocity-diagonal-bound}
\begin{split}
C_{f,b}(t,t)
&=
2^{b-2}\int_0^t f(u)(t-u)^{b-2}\,du\\
&\le
2^{b-2}c
\int_0^t u^a(t-u)^{b-2}\,du\\
&=
2^{b-2}cB(a+1,b-1)t^{a+b-1}.
\end{split}
\end{equation}
If \(0\le s<t\le T\), then
\[
s+t-2u\ge t-s,
\qquad 0\le u\le s.
\]
Since \(b-2\le0\),
\[
C_{f,b}(s,t)
\le
(t-s)^{b-2}\int_0^s f(u)\,du
<\infty,
\]
where the final integral is finite because
\(f(u)\le cu^a\) and \(a>-1\). Thus \(C_{f,b}\) is finite on
\([0,T]^2\).

We next construct a jointly measurable realization. First suppose that
\(1<b<2\), and set
\[
\mathcal H
=
L^2\!\left(
(0,T)\times(0,\infty),
\frac{f(u)}{\Gamma(2-b)}
\lambda^{1-b}\,du\,d\lambda
\right),
\]
with
\[
k_t(u,\lambda)
=
\mathbf 1_{[0,t)}(u)e^{-\lambda(t-u)}.
\]
For \(b=2\), instead set
\[
\mathcal H=L^2((0,T),f(u)\,du),
\qquad
k_t(u)=\mathbf 1_{[0,t)}(u).
\]
In either case,
\[
\langle k_s,k_t\rangle_{\mathcal H}
=
C_{f,b}(s,t),
\]
and \(\mathcal H\) is separable.

For every fixed \(h\in\mathcal H\), the map
\[
t\longmapsto\langle k_t,h\rangle_{\mathcal H}
\]
is measurable by Fubini's theorem. Since \(\mathcal H\) is separable,
Pettis' measurability theorem implies that
\(t\mapsto k_t\) is strongly measurable. Moreover,
\eqref{eq:velocity-diagonal-bound} gives
\[
\int_0^T\|k_t\|_{\mathcal H}^2\,dt
=
\int_0^T C_{f,b}(t,t)\,dt
<\infty.
\]

If \(\mathcal H=\{0\}\), set \(\nu_t=0\) for every \(t\). Otherwise,
choose an orthonormal basis \((e_j)_{j\ge1}\) of \(\mathcal H\), let
\((\xi_j)_{j\ge1}\) be independent standard normal random variables,
and define
\[
\nu_t^{(N)}
=
\sum_{j=1}^N
\langle k_t,e_j\rangle_{\mathcal H}\xi_j.
\]
For each fixed \(t\),
\[
\sum_{j=1}^\infty
|\langle k_t,e_j\rangle_{\mathcal H}|^2
=
\|k_t\|_{\mathcal H}^2
<\infty.
\]
Hence the Gaussian series converges almost surely and in
\(L^2(\Omega)\). Define \(\nu_t\) as its limit on the measurable set
where the series converges, and set it equal to zero elsewhere.

Each partial sum is jointly measurable, so
\[
(t,\omega)\longmapsto\nu_t(\omega)
\]
is jointly measurable. Moreover,
\(t\mapsto\nu_t^{(N)}\) is strongly measurable as an
\(L^2(\Omega)\)-valued map, and
\[
\nu_t^{(N)}\longrightarrow\nu_t
\qquad\text{in }L^2(\Omega)
\]
for every \(t\). Thus \(t\mapsto\nu_t\) is strongly measurable as an
\(L^2(\Omega)\)-valued map.

In both the zero and non-zero Hilbert-space cases, the
finite-dimensional distributions are centered Gaussian and
\[
\mathbb E(\nu_s\nu_t)
=
\langle k_s,k_t\rangle_{\mathcal H}
=
C_{f,b}(s,t).
\]

By \eqref{eq:velocity-diagonal-bound},
\begin{align}
\int_0^T\|\nu_s\|_{L^2(\Omega)}\,ds
&=
\int_0^T C_{f,b}(s,s)^{1/2}\,ds \notag\\
&\le
2^{(b-2)/2}
\sqrt{cB(a+1,b-1)}
\int_0^T s^{(a+b-1)/2}\,ds
<\infty.
\label{eq:velocity-bochner-integrability}
\end{align}
The last integral is finite because
\[
a+b-1>-1.
\]
Therefore \(s\mapsto\nu_s\) is Bochner integrable as an
\(L^2(\Omega)\)-valued map. The same bound and Tonelli's theorem give
\[
\mathbb E\int_0^T|\nu_s|\,ds
\le
\int_0^T\|\nu_s\|_{L^2(\Omega)}\,ds
<\infty,
\]
which proves \eqref{eq:velocity-L1}.

The pathwise Lebesgue integral and the Bochner integral agree almost
surely. To verify that the integral process is Gaussian, let
\[
I_N(t)
:=
\int_0^t\nu_s^{(N)}\,ds
=
\sum_{j=1}^N
\left(
\int_0^t\langle k_s,e_j\rangle_{\mathcal H}\,ds
\right)\xi_j.
\]
Each \(I_N(t)\) is centered Gaussian. Moreover, Minkowski's inequality
gives
\[
\left\|
I_N(t)-\int_0^t\nu_s\,ds
\right\|_{L^2(\Omega)}
\le
\int_0^t
\|\nu_s^{(N)}-\nu_s\|_{L^2(\Omega)}\,ds
\longrightarrow0.
\]
This follows by dominated convergence, since
\[
\|\nu_s^{(N)}-\nu_s\|_{L^2(\Omega)}
\le
\|k_s\|_{\mathcal H},
\]
and the latter function is integrable by
\eqref{eq:velocity-bochner-integrability}.

More generally, for arbitrary \(t_1,\ldots,t_m\in[0,T]\) and
\(c_1,\ldots,c_m\in\mathbb R\),
\[
\sum_{i=1}^m c_iI_N(t_i)
\longrightarrow
\sum_{i=1}^m c_i\int_0^{t_i}\nu_s\,ds
\qquad\text{in }L^2(\Omega).
\]
Since the left-hand side is Gaussian for every \(N\), its limit is
Gaussian. Hence
\[
\left(\int_0^t\nu_s\,ds\right)_{0\le t\le T}
\]
is a centered Gaussian process.

Furthermore, Cauchy--Schwarz and
\eqref{eq:velocity-bochner-integrability} justify Fubini's theorem:
\[
\mathbb E
\left[
\left(\int_0^s\nu_r\,dr\right)
\left(\int_0^t\nu_v\,dv\right)
\right]
=
\int_0^s\int_0^t
C_{f,b}(r,v)\,dv\,dr.
\]

For \(x,y\ge0\) and \(1<b\le2\), direct integration gives
\[
Q_b(x,y)
=
b\int_0^x\int_0^y
(r+v)^{b-2}\,dv\,dr.
\]
All the integrands below are non-negative, so Tonelli's theorem yields
\begin{align}
R_{f,b}(s,t)
&=
b\int_0^{s\wedge t}
f(u)
\int_u^s\int_u^t
(r+v-2u)^{b-2}\,dv\,dr\,du
\notag\\
&=
b\int_0^s\int_0^t
\left(
\int_0^{r\wedge v}
f(u)(r+v-2u)^{b-2}\,du
\right)dv\,dr
\notag\\
&=
b\int_0^s\int_0^t
C_{f,b}(r,v)\,dv\,dr.
\label{eq:covariance-double-integral}
\end{align}
It follows that the process in \eqref{eq:ac-realization} has covariance
\(R_{f,b}\). By \eqref{eq:velocity-L1}, its trajectories are absolutely
continuous and therefore have finite total variation on \([0,T]\).

Suppose now that \(a\ge0\). Since the value of \(f\) at the single point
\(0\) does not affect any covariance integral, we may redefine \(f(0)\),
if necessary, so that
\[
0\le f(u)\le M,
\qquad 0\le u\le T,
\]
for some \(M<\infty\).

Put
\[
\beta:=b-1\in(0,1].
\]
For the constant weight \(f\equiv1\), direct integration gives
\[
C_{1,b}(s,t)
=
\frac{(s+t)^\beta-|t-s|^\beta}{2\beta}.
\]
For \(0\le s\le t\le T\),
\begin{align*}
\mathbb E|\nu_{t,1,b}-\nu_{s,1,b}|^2
&=
\frac{
(2t)^\beta+(2s)^\beta
-2(s+t)^\beta
+2(t-s)^\beta
}{2\beta}\\
&\le
\frac{|t-s|^\beta}{\beta},
\end{align*}
because the concavity of \(x\mapsto x^\beta\) gives
\[
(2t)^\beta+(2s)^\beta
\le
2(s+t)^\beta.
\]

For every fixed \(u\), the elementary positive-definite kernel used
above has non-negative increment variance. Consequently,
\(0\le f\le M\) implies
\[
\mathbb E|\nu_{t,f,b}-\nu_{s,f,b}|^2
\le
M\,
\mathbb E|\nu_{t,1,b}-\nu_{s,1,b}|^2
\le
\frac{M}{b-1}|t-s|^{b-1}.
\]
Gaussian moment equivalence and Kolmogorov's continuity theorem
therefore give a continuous modification \(\widehat\nu\) whose paths
are H\"older continuous of every order
\[
\delta<\frac{b-1}{2}.
\]

For every fixed \(t\),
\[
\nu_t=\widehat\nu_t
\qquad\text{almost surely}.
\]
Fubini's theorem consequently gives
\[
\nu_t=\widehat\nu_t
\qquad
dt\otimes d\mathbb P\text{-almost everywhere}.
\]
Hence, on an event of probability one,
\[
\int_0^t\nu_s\,ds
=
\int_0^t\widehat\nu_s\,ds
\qquad
\text{for every }t\in[0,T].
\]
Thus the realization in \eqref{eq:ac-realization} may be written using
\(\widehat\nu\). The fundamental theorem of calculus then gives
\[
\frac{d\widetilde\zeta_{t,f,b}}{dt}
=
\sqrt b\,\widehat\nu_{t,f,b},
\qquad 0\le t\le T.
\]
Relabeling \(\widehat\nu\) as \(\nu\) proves the final assertion.
\end{proof}

\subsubsection{Proof of Proposition~\ref{prop:remote-covariance}}\label{sec:supp-proof-prop4}
\begin{proof}
For all sufficiently large \(T\), one has \(s\leq T+t\), and hence
\[
R_{f,b}(s,T+t)
=
\frac{1}{1-b}
\int_0^s f(u)
\Bigl\{
(s-u)^b+(T+t-u)^b-(T+t+s-2u)^b
\Bigr\}\,du.
\]

For \(b=0\), \eqref{eq:remote-level-rough} follows directly from the
time-change covariance. Suppose that \(0<b<1\). Concavity and
subadditivity of \(x\mapsto x^b\) give
\[
0\leq
(T+t+s-2u)^b-(T+t-u)^b
\leq
(s-u)^b,
\qquad 0\leq u\leq s.
\]
The difference on the left converges to zero uniformly for
\(u\in[0,s]\) as \(T\to\infty\). Moreover,
\[
0\leq
(s-u)^b+(T+t-u)^b-(T+t+s-2u)^b
\leq
(s-u)^b.
\]
Since \(f(u)(s-u)^b\) is integrable by \eqref{eq:weight-finiteness}, dominated convergence proves
\eqref{eq:remote-level-rough}.

Suppose now that \(1<b\leq2\). The contribution of the first term
\((s-u)^b\) is \(o(T^{b-1})\). For the remaining difference, the
mean-value theorem gives, for each \(u\in[0,s]\), a point
\(\xi_{T,u}\) between \(T+t-u\) and \(T+t+s-2u\) such that
\[
\frac{
(T+t+s-2u)^b-(T+t-u)^b
}{T^{b-1}}
=
b(s-u)
\left(\frac{\xi_{T,u}}{T}\right)^{b-1}.
\]
Uniformly for \(u\in[0,s]\),
\[
\frac{\xi_{T,u}}{T}\longrightarrow1,
\qquad T\to\infty,
\]
and this ratio is uniformly bounded for all sufficiently large \(T\).
Since \(f\) is locally integrable and
\(f(u)(s-u)\) is integrable on \([0,s]\), dominated convergence gives
\[
\frac{1}{T^{b-1}}
\int_0^s f(u)
\Bigl\{
(T+t-u)^b-(T+t+s-2u)^b
\Bigr\}\,du
\longrightarrow
-b\int_0^s f(u)(s-u)\,du.
\]
Using
\[
\frac{1}{1-b}=-\frac{1}{b-1}
\]
proves \eqref{eq:remote-level-smooth}.

We now prove the increment asymptotic. Fix \(x>0\). For all
sufficiently large \(T\), every \(q\in[s,t]\), and every
\(u\in[0,x]\), the quantities
\[
T+q-u
\quad\text{and}\quad
T+q+x-2u
\]
are bounded away from zero. Since \(f\) is locally integrable,
differentiation under the integral sign is justified and gives
\[
\frac{\partial}{\partial q}R_{f,b}(x,T+q)
=
\frac{b}{1-b}
\int_0^x f(u)
\Bigl\{
(T+q-u)^{b-1}
-
(T+q+x-2u)^{b-1}
\Bigr\}\,du.
\]

For each sufficiently large \(T\), \(q\in[s,t]\), and \(u\in[0,x]\),
the mean-value theorem gives a point \(\eta_{T,q,u}\) between
\(T+q-u\) and \(T+q+x-2u\) such that
\[
\begin{split}
&T^{2-b}\frac{b}{1-b}
\Bigl\{
(T+q-u)^{b-1}
-
(T+q+x-2u)^{b-1}
\Bigr\}\\
&\qquad=
b(x-u)
\left(\frac{\eta_{T,q,u}}{T}\right)^{b-2}.
\end{split}
\]
Uniformly for \(q\in[s,t]\) and \(u\in[0,x]\),
\[
\frac{\eta_{T,q,u}}{T}\longrightarrow1.
\]
Moreover, this ratio is bounded above and bounded away from zero for
all sufficiently large \(T\). It follows that
\[
T^{2-b}
\left|
\frac{b}{1-b}
\Bigl\{
(T+q-u)^{b-1}
-
(T+q+x-2u)^{b-1}
\Bigr\}
\right|
\leq
C_{b,x,s,t}(x-u).
\]
Since \(f(u)(x-u)\) is integrable on \([0,x]\), dominated convergence,
uniformly in \(q\in[s,t]\), yields
\[
T^{2-b}
\frac{\partial}{\partial q}R_{f,b}(x,T+q)
\longrightarrow
b\int_0^x f(u)(x-u)\,du
=
bA_f(x).
\]

The fundamental theorem of calculus gives
\begin{align*}
&\mathbb E\left[
(\zeta_{v,f,b}-\zeta_{r,f,b})
(\zeta_{T+t,f,b}-\zeta_{T+s,f,b})
\right]\\
&\qquad=
\int_s^t
\frac{\partial}{\partial q}
\Bigl\{
R_{f,b}(v,T+q)-R_{f,b}(r,T+q)
\Bigr\}\,dq.
\end{align*}
The preceding convergence is uniform for \(q\in[s,t]\), and therefore
\[
\begin{split}
\lim_{T\to\infty}T^{2-b}
&\mathbb E\left[
(\zeta_{v,f,b}-\zeta_{r,f,b})
(\zeta_{T+t,f,b}-\zeta_{T+s,f,b})
\right]\\
&=
b(t-s)\{A_f(v)-A_f(r)\}.
\end{split}
\]
Finally,
\[
\begin{split}
A_f(v)-A_f(r)
&=
\int_r^v f(u)(v-u)\,du
+
\int_0^r f(u)\bigl[(v-u)-(r-u)\bigr]\,du\\
&=
\int_r^v f(u)(v-u)\,du
+
(v-r)\int_0^r f(u)\,du,
\end{split}
\]
which proves \eqref{eq:remote-increment-limit}.
\end{proof}

\subsubsection{Proof of Proposition~\ref{prop:nonstationary}}\label{sec:supp-proof-prop5}
\begin{proof}
Since
\[
\frac{2-2^b}{1-b}>0
\]
for \(b\in(0,1)\cup(1,2]\), and since \(f>0\) on a set of positive
Lebesgue measure, there exists \(t>0\) such that
\[
R_{f,b}(t,t)
=
\frac{2-2^b}{1-b}
\int_0^t f(u)(t-u)^b\,du
>0.
\]
On the other hand, \(\zeta_{0,f,b}=0\) almost surely. Hence the variance
is not constant, and weak stationarity is impossible.

Suppose first that \textup{(a)} holds. By
Proposition~\ref{prop:local-increments}\textup{(iii)},
\[
\mathbb E
\bigl|\zeta_{s+\varepsilon,f,b}-\zeta_{s,f,b}\bigr|^2
\sim
2^{b-2}b\varepsilon^2
\int_0^s f(u)(s-u)^{b-2}\,du,
\qquad \varepsilon\downarrow0.
\]
The constant in this asymptotic relation is strictly positive by
assumption. On the other hand, local boundedness of \(f\) near zero and
the diagonal formula give
\[
\mathbb E\zeta_{\varepsilon,f,b}^2
=
\frac{2-2^b}{1-b}
\int_0^\varepsilon f(u)(\varepsilon-u)^b\,du
=
O(\varepsilon^{b+1}).
\]
Since \(b+1>2\),
\[
\mathbb E\zeta_{\varepsilon,f,b}^2
=
o(\varepsilon^2).
\]
Consequently,
\[
\frac{
\mathbb E|\zeta_{s+\varepsilon,f,b}-\zeta_{s,f,b}|^2
}{\varepsilon^2}
\longrightarrow
2^{b-2}b
\int_0^s f(u)(s-u)^{b-2}\,du
>0,
\]
whereas
\[
\frac{
\mathbb E|\zeta_{\varepsilon,f,b}-\zeta_{0,f,b}|^2
}{\varepsilon^2}
\longrightarrow0.
\]
Thus, for all sufficiently small \(\varepsilon>0\), increments of the
same length \(\varepsilon\), starting respectively at \(0\) and at
\(s\), have different variances. Hence the process is not intrinsically
stationary.

Suppose now that \textup{(b)} holds. For \(s,h\ge0\), write
\[
V_s(h)
:=
\mathbb E
\bigl|\zeta_{s+h,f,b}-\zeta_{s,f,b}\bigr|^2.
\]
Formula \eqref{eq:increment-variance-decomposition} gives
\begin{equation}\label{eq:nonstationarity-variance-difference}
\begin{split}
V_s(h)-V_0(h)
&=
\frac{2-2^b}{1-b}
\int_0^h
\bigl\{f(s+v)-f(v)\bigr\}(h-v)^b\,dv\\
&\quad+J_s(h),
\end{split}
\end{equation}
where
\[
J_s(h)
=
\frac{1}{1-b}
\int_0^s f(u)
\Bigl\{
2(h+2s-2u)^b
-2^b(h+s-u)^b
-2^b(s-u)^b
\Bigr\}\,du.
\]

Rectangular integration of the mixed derivative of
\((r+r'-2u)^b\), with the endpoint \(u=s\) being immaterial, gives the
non-negative identity
\begin{equation}\label{eq:nonstationarity-positive-triple}
J_s(h)
=
b\int_s^{s+h}\int_s^{s+h}\int_0^s
f(u)(r+r'-2u)^{b-2}\,du\,dr'\,dr.
\end{equation}

Since \(f\) is continuous at zero and \(f(0)>0\), there exist
\(\delta>0\) and \(m>0\) such that
\[
f(u)\ge m,
\qquad 0\le u\le\delta.
\]
Choose \(s\in(0,\delta/2)\) and put
\[
h:=\frac{s}{2}.
\]
In \eqref{eq:nonstationarity-positive-triple}, restrict the integration
in \(u\) to
\[
s-\frac h2\le u\le s.
\]
This interval is contained in \([0,\delta]\), and hence \(f(u)\ge m\)
there. Moreover, for
\[
s\le r,r'\le s+h,
\qquad
s-\frac h2\le u\le s,
\]
one has
\[
r+r'-2u\le3h.
\]
Since \(b-2<0\),
\[
(r+r'-2u)^{b-2}\ge(3h)^{b-2}.
\]
Therefore
\begin{align}
J_s(h)
&\ge
bm
\int_s^{s+h}\int_s^{s+h}
\int_{s-h/2}^s
(3h)^{b-2}\,du\,dr'\,dr \notag\\
&=
bm\,h^2\frac h2(3h)^{b-2}\notag\\
&=
\frac{bm}{2}3^{b-2}h^{b+1}.
\label{eq:nonstationarity-positive-lower-bound}
\end{align}

For the first term on the right-hand side of
\eqref{eq:nonstationarity-variance-difference}, define
\[
\omega_f(s)
:=
\sup_{0\le v\le s/2}
|f(s+v)-f(v)|.
\]
Since both \(v\) and \(s+v\) converge uniformly to zero as
\(s\downarrow0\), continuity of \(f\) at zero implies that
\[
\omega_f(s)\longrightarrow0.
\]
Because \(h=s/2\),
\begin{align}
&\left|
\frac{2-2^b}{1-b}
\int_0^h
\bigl\{f(s+v)-f(v)\bigr\}(h-v)^b\,dv
\right| \notag\\
&\qquad\le
\frac{2-2^b}{1-b}\,
\omega_f(s)\int_0^h(h-v)^b\,dv \notag\\
&\qquad=
\frac{2-2^b}{(1-b)(b+1)}
\omega_f(s)h^{b+1}.
\label{eq:nonstationarity-remainder-bound}
\end{align}

Combining
\eqref{eq:nonstationarity-variance-difference},
\eqref{eq:nonstationarity-positive-lower-bound}, and
\eqref{eq:nonstationarity-remainder-bound}, we obtain
\[
V_s(h)-V_0(h)
\ge
\left[
\frac{bm}{2}3^{b-2}
-
\frac{2-2^b}{(1-b)(b+1)}
\omega_f(s)
\right]h^{b+1}.
\]
Since \(\omega_f(s)\to0\), the expression in brackets is strictly
positive for all sufficiently small \(s>0\). Hence, with \(h=s/2\),
\[
V_s(s/2)>V_0(s/2)
\]
for all sufficiently small \(s>0\). Thus increments with the same lag
have different variances, and intrinsic stationarity is impossible.
\end{proof}

\subsubsection{Proof of Proposition~\ref{prop:local-increments}}\label{sec:supp-proof-prop6}
\begin{proof}
For \(b=0\), the time-change covariance gives
\[
\mathbb E
|\zeta_{s+\varepsilon,f,0}-\zeta_{s,f,0}|^2
=
\int_s^{s+\varepsilon}f(u)\,du.
\]
Dividing by \(\varepsilon\) and using the continuity of \(f\) at \(s\)
proves \eqref{eq:local-b0}.

Part \textup{(ii)} follows directly from
\eqref{eq:bounded-weight-increment-bounds} in the case
\(1\asymp f\), and from
\eqref{eq:power-comparable-increment-bounds} in the case
\(u^a\asymp f(u)\).

For part \textup{(iii)}, define
\[
G(r,v)
:=
\int_0^{r\wedge v}
f(u)(r+v-2u)^{b-2}\,du.
\]
By \eqref{eq:covariance-double-integral},
\[
R_{f,b}(r,v)
=
b\int_0^r\int_0^v G(p,q)\,dq\,dp.
\]
Taking the rectangular covariance increment over
\([s,s+\varepsilon]^2\) gives
\begin{equation}\label{eq:local-smooth-average}
\frac{
\mathbb E|\zeta_{s+\varepsilon,f,b}-\zeta_{s,f,b}|^2
}{\varepsilon^2}
=
\frac{b}{\varepsilon^2}
\int_s^{s+\varepsilon}
\int_s^{s+\varepsilon}
G(r,v)\,dv\,dr.
\end{equation}

It remains to prove that
\[
G(r,v)\longrightarrow G(s,s)
\qquad\text{as }(r,v)\to(s,s),\quad r,v\ge s.
\]
On \([0,s]\), since \(b-2\le0\) and \(r,v\ge s\),
\[
0\le
f(u)(r+v-2u)^{b-2}
\le
2^{b-2}f(u)(s-u)^{b-2}.
\]
The function on the right is integrable. Near \(u=s\), this follows
from the local boundedness of \(f\) and \(b-2>-1\); away from \(s\), it
follows from the local integrability of \(f\), which is implied by \eqref{eq:weight-finiteness}. Dominated convergence therefore gives
\[
\int_0^s f(u)(r+v-2u)^{b-2}\,du
\longrightarrow
2^{b-2}\int_0^s f(u)(s-u)^{b-2}\,du.
\]

Let \(m=r\wedge v\). For \(r\) and \(v\) sufficiently close to \(s\),
local boundedness gives \(f(u)\le M\) on \([s,m]\). Moreover,
\[
r+v-2u\ge2(m-u),
\qquad s\le u\le m.
\]
Hence
\begin{align*}
\int_s^m f(u)(r+v-2u)^{b-2}\,du
&\le
M2^{b-2}\int_s^m(m-u)^{b-2}\,du\\
&=
\frac{M2^{b-2}}{b-1}(m-s)^{b-1}
\longrightarrow0.
\end{align*}
Consequently,
\[
G(r,v)\longrightarrow
2^{b-2}\int_0^s f(u)(s-u)^{b-2}\,du.
\]
Applying this convergence to the average over the shrinking square in
\eqref{eq:local-smooth-average} proves \eqref{eq:local-smooth-limit}.
\end{proof}

\subsubsection{Proof of Proposition~\ref{prop:small-time-diagonal}}\label{sec:supp-proof-prop7}
\begin{proof}
For \(b=0\),
\[
R_{f,0}(\varepsilon,\varepsilon)
=
\int_0^\varepsilon f(u)\,du,
\]
and therefore continuity at zero gives
\[
\frac{R_{f,0}(\varepsilon,\varepsilon)}{\varepsilon}
\longrightarrow f(0).
\]

Let \(b\ne0\). By the diagonal identity \eqref{diag-covariance-general} and the change of variables
\(u=\varepsilon v\),
\begin{align*}
\frac{\mathbb E(\zeta_{\varepsilon,f,b}^2)}
     {\varepsilon^{1+b}}
&=
\frac{2-2^b}{1-b}
\int_0^1 f(\varepsilon v)(1-v)^b\,dv.
\end{align*}
Since \(f\) is continuous at zero, it is bounded in a neighborhood of
zero, and dominated convergence gives
\[
\frac{\mathbb E(\zeta_{\varepsilon,f,b}^2)}
     {\varepsilon^{1+b}}
\longrightarrow
\frac{2-2^b}{1-b}
f(0)\int_0^1(1-v)^b\,dv
=
\frac{2-2^b}{(1-b)(1+b)}f(0).
\]
This proves \eqref{small-time-diagonal-general}.

If \(f(u)=u^a\), the same change of variables gives the exact identity
\[
\mathbb E(\zeta_{\varepsilon,f,b}^2)
=
\frac{2-2^b}{1-b}
\varepsilon^{a+b+1}
\int_0^1 v^a(1-v)^b\,dv
=
\frac{2-2^b}{1-b}
B(a+1,b+1)\varepsilon^{a+b+1}.
\]
For \(b=0\), this agrees with
\[
\int_0^\varepsilon u^a\,du
=
\frac{\varepsilon^{a+1}}{a+1},
\]
because \(B(a+1,1)=1/(a+1)\). This proves \eqref{small-time-diagonal-power}.
\end{proof}

\subsubsection{Proof of Proposition~\ref{prop:logarithmic-limit}}\label{sec:supp-proof-prop8}
\begin{proof}
Put
\[
q_b(x,y)
:=
\frac{x^b+y^b-(x+y)^b}{1-b},
\qquad x,y\ge0.
\]
For fixed \(x,y\ge0\), l'H\^opital's rule, with the convention
\(0\log0=0\), gives
\[
q_b(x,y)
\longrightarrow
(x+y)\log(x+y)-x\log x-y\log y,
\qquad b\to1.
\]

Fix \(s,t\ge0\), and set
\[
L:=1+s+t.
\]
Choose \(\varepsilon_0\in(0,1)\) and restrict
\[
b\in[1-\varepsilon_0,1+\varepsilon_0].
\]
For \(z>0\), the mean-value theorem applied to the map
\(\theta\mapsto z^\theta\) gives a point \(\theta\) between \(1\)
and \(b\) such that
\[
\left|
\frac{z^b-z}{1-b}
\right|
=
z^\theta|\log z|.
\]
Consequently, uniformly for \(0<z\le L\),
\[
\left|
\frac{z^b-z}{1-b}
\right|
\le
C_{L,\varepsilon_0}
z^{1-\varepsilon_0}(1+|\log z|).
\]
The expression on the right extends continuously at \(z=0\) with
value zero and is bounded on \([0,L]\).

Since
\[
q_b(x,y)
=
\frac{x^b-x}{1-b}
+
\frac{y^b-y}{1-b}
-
\frac{(x+y)^b-(x+y)}{1-b},
\]
the preceding estimate implies
\[
|q_b(s-u,t-u)|
\le C_{s,t,\varepsilon_0},
\qquad 0\le u\le s\wedge t,
\]
uniformly for \(b\) in the indicated neighborhood of one. Local
integrability of \(f\) therefore permits dominated convergence and
gives
\[
R_{f,b}(s,t)
=
\int_0^{s\wedge t}
f(u)q_b(s-u,t-u)\,du
\longrightarrow
K_f(s,t).
\]
The same bound also shows directly that \(K_f(s,t)\) is finite.

To prove positive definiteness, fix
\(t_1,\ldots,t_m\ge0\) and \(c_1,\ldots,c_m\in\mathbb R\). Let
\(b_n\to1\) with \(b_n\in(0,1)\). Since \(f\) is locally integrable,
it is an admissible weight for every \(b_n>0\). Hence
Theorem~\ref{thm:maximal-range} gives
\[
\sum_{i,j=1}^m
c_ic_jR_{f,b_n}(t_i,t_j)\ge0.
\]
Passing to the limit yields
\[
\sum_{i,j=1}^m
c_ic_jK_f(t_i,t_j)\ge0.
\]
Thus \(K_f\) is positive definite.

It remains to prove the joint limit. Fix \(s,t\ge0\), retain
\(L=1+s+t\), and choose \(\varepsilon_0\in(0,1)\). If
\(|a|\le\varepsilon_0\), then, for \(0<u\le s\wedge t\),
\[
u^a\le L^{\varepsilon_0}+u^{-\varepsilon_0},
\qquad
e^{au}\le e^{\varepsilon_0L}.
\]
These bounds are integrable on \((0,s\wedge t]\), because
\(\varepsilon_0<1\). Moreover,
\[
u^a\longrightarrow1,
\qquad
e^{au}\longrightarrow1,
\qquad a\to0,
\]
for every \(u>0\). Combining these weight bounds with the preceding
uniform domination of \(q_b\), dominated convergence gives
\[
R_{f_a,b}(s,t)\longrightarrow K_1(s,t)
\qquad\text{as }(a,b)\to(0,1)
\]
for both families.
\end{proof}

\subsection{Detailed proofs for the pathwise equations and Euler approximation}\label{sec:supp-euler-proofs}

\subsubsection{Auxiliary estimates for the Young equation}\label{sec:supp-young-aux}

\begin{lem}\label{lem:composition-holder}
Let \(0<\rho<1\), \(0\le s<t\le T\), and let
\(F:\mathbb R\times[0,T]\to\mathbb R\) have bounded derivatives
\[
\partial_xF,\qquad \partial_{xx}F,\qquad \partial_{xt}F.
\]
For \(x,y\in C^\rho([s,t])\), set
\[
h_r:=F(x_r,r)-F(y_r,r).
\]
Then
\begin{align}
\|h\|_{\infty;[s,t]}
&\le
C\|x-y\|_{\infty;[s,t]},
\label{eq:comp-sup}\\
[h]_{\rho;[s,t]}
&\le
C[x-y]_{\rho;[s,t]}
+
C\bigl(
1+[x]_{\rho;[s,t]}+[y]_{\rho;[s,t]}
\bigr)
\|x-y\|_{\infty;[s,t]}.
\label{eq:comp-holder}
\end{align}
The constant depends only on the indicated derivative bounds,
\(T\), and \(\rho\).
\end{lem}

\begin{proof}
Put
\[
d_r:=x_r-y_r.
\]
By the mean-value formula,
\[
h_r=A_rd_r,
\qquad
A_r
=
\int_0^1
\partial_xF\bigl(y_r+\theta d_r,r\bigr)\,d\theta.
\]
Since \(\partial_xF\) is bounded,
\[
\|h\|_{\infty;[s,t]}
\le
\|\partial_xF\|_\infty
\|x-y\|_{\infty;[s,t]},
\]
which proves \eqref{eq:comp-sup}.

For \(r,q\in[s,t]\), the boundedness of
\(\partial_{xx}F\) and \(\partial_{xt}F\) gives
\begin{align*}
|A_r-A_q|
&\le
C\int_0^1
\left(
\left|
y_r+\theta d_r-y_q-\theta d_q
\right|
+
|r-q|
\right)d\theta\\
&\le
C\bigl([x]_{\rho;[s,t]}+[y]_{\rho;[s,t]}\bigr)
|r-q|^\rho
+
C|r-q|.
\end{align*}
Since
\[
|r-q|
\le
(t-s)^{1-\rho}|r-q|^\rho
\le
T^{1-\rho}|r-q|^\rho,
\]
it follows that
\[
[A]_{\rho;[s,t]}
\le
C\bigl(
1+[x]_{\rho;[s,t]}+[y]_{\rho;[s,t]}
\bigr).
\]
Finally, using
\[
[Ad]_\rho
\le
\|A\|_\infty[d]_\rho
+
[A]_\rho\|d\|_\infty
\]
with \(d=x-y\) proves \eqref{eq:comp-holder}.
\end{proof}

\begin{lem}\label{lem:young-local-fixed-radius}
Assume \textup{(H)}, let \(1/2<\rho<1\), and let
\(Z\in C^\rho([0,T])\). Put
\[
K:=[Z]_{\rho;[0,T]}.
\]
There exist deterministic constants \(R_0\ge1\),
\(\eta\in(0,1)\), and \(C<\infty\), depending only on
\textup{(H)}, \(T\), and \(\rho\), with the following property.
If \(I=[s,s+\ell]\subseteq[0,T]\), with \(\ell>0\), satisfies
\begin{equation}\label{eq:young-smallness-condition}
\ell\le\eta,
\qquad
K\ell^\rho\le\eta,
\end{equation}
then, for every prescribed initial value \(\xi\in\mathbb R\), the
Picard map
\[
(\Phi_\xi x)_u
=
\xi
+
\int_s^u\mathfrak b(x_r,r)\,dr
+
\int_s^u\sigma(x_r,r)\,dZ_r
\]
is a strict contraction, with respect to the norm
\[
\|x\|_{I,*}
:=
\|x\|_{\infty;I}
+
\ell^\rho[x]_{\rho;I},
\]
on the closed set
\[
\mathcal B_I(\xi)
:=
\left\{
x\in C^\rho(I):
x_s=\xi,\ 
\ell^\rho[x]_{\rho;I}\le R_0
\right\}.
\]
Since \(\|\cdot\|_{I,*}\) is equivalent to the usual
\(C^\rho(I)\)-norm for fixed \(\ell>0\), the set
\(\mathcal B_I(\xi)\) is complete. Its fixed point \(X\) satisfies
\begin{equation}\label{eq:local-solution-holder-bound}
[X]_{\rho;I}\le C(1+K).
\end{equation}
If \(X\) and \(Y\) are two solutions on \(I\), driven by the same
path \(Z\) but possibly having different initial values, then
\begin{equation}\label{eq:local-flow-lipschitz-two}
\|X-Y\|_{I,*}
\le
2|X_s-Y_s|.
\end{equation}
\end{lem}

\begin{proof}
For \(q\in C^\rho(I)\), the ordinary integral satisfies
\begin{equation}\label{eq:scaled-drift-bound}
\left\|
u\longmapsto\int_s^u q_r\,dr
\right\|_{I,*}
\le
C\ell\|q\|_{\infty;I}.
\end{equation}
Indeed, its supremum norm is bounded by
\(\ell\|q\|_{\infty;I}\), while its \(\rho\)-Hölder seminorm is
bounded by
\(\ell^{1-\rho}\|q\|_{\infty;I}\).

Young--Loeve's inequality gives
\[
\left|
\int_u^v q_r\,dZ_r-q_uZ_{u,v}
\right|
\le
C_\rho[q]_{\rho;I}K|v-u|^{2\rho},
\qquad s\le u<v\le s+\ell.
\]
It follows that
\begin{equation}\label{eq:scaled-young-bound}
\left\|
u\longmapsto\int_s^u q_r\,dZ_r
\right\|_{I,*}
\le
CK\ell^\rho
\left(
\|q\|_{\infty;I}
+
\ell^\rho[q]_{\rho;I}
\right).
\end{equation}

Let \(x\in\mathcal B_I(\xi)\). By the boundedness of
\(\partial_x\sigma\) and \(\partial_t\sigma\),
\[
[\sigma(x,\cdot)]_{\rho;I}
\le
C[x]_{\rho;I}
+
C\ell^{1-\rho}.
\]
Consequently,
\[
\ell^\rho[\sigma(x,\cdot)]_{\rho;I}
\le
C(R_0+\ell).
\]
Since \(\mathfrak b\) and \(\sigma\) are bounded,
\eqref{eq:scaled-drift-bound} and
\eqref{eq:scaled-young-bound} imply
\[
\ell^\rho[\Phi_\xi x]_{\rho;I}
\le
C\ell
+
CK\ell^\rho(1+R_0+\ell).
\]
Choose a deterministic \(R_0\ge1\). By
\eqref{eq:young-smallness-condition},
\[
\ell^\rho[\Phi_\xi x]_{\rho;I}
\le
C\eta+C\eta(1+R_0+\eta).
\]
We may choose \(\eta\in(0,1)\), depending only on the deterministic
coefficient bounds, \(T\), and \(\rho\), so small that
\[
C\eta+C\eta(1+R_0+\eta)\le R_0.
\]
Since \((\Phi_\xi x)_s=\xi\), this proves that
\(\Phi_\xi\) maps \(\mathcal B_I(\xi)\) into itself.

Let now \(x,y\in\mathcal B_I(\xi)\), and put \(d=x-y\). The drift
contribution satisfies
\[
\left\|
u\longmapsto
\int_s^u
\{\mathfrak b(x_r,r)-\mathfrak b(y_r,r)\}\,dr
\right\|_{I,*}
\le
C\ell\|d\|_{I,*}.
\]
Apply Lemma~\ref{lem:composition-holder} with \(F=\sigma\). We have
\[
\|\sigma(x,\cdot)-\sigma(y,\cdot)\|_{\infty;I}
\le
C\|d\|_{\infty;I}
\]
and
\begin{align*}
\ell^\rho
[\sigma(x,\cdot)-\sigma(y,\cdot)]_{\rho;I}
&\le
C\ell^\rho[d]_{\rho;I}\\
&\quad+
C\bigl(
\ell^\rho
+\ell^\rho[x]_{\rho;I}
+\ell^\rho[y]_{\rho;I}
\bigr)
\|d\|_{\infty;I}\\
&\le
C(1+R_0)\|d\|_{I,*}.
\end{align*}
Therefore, by \eqref{eq:scaled-young-bound},
\[
\|\Phi_\xi x-\Phi_\xi y\|_{I,*}
\le
C\left\{
\ell+(1+R_0)K\ell^\rho
\right\}
\|x-y\|_{I,*}.
\]
After decreasing \(\eta\), if necessary, the coefficient on the
right-hand side is at most \(1/2\). Thus \(\Phi_\xi\) is a strict
contraction on \(\mathcal B_I(\xi)\), and it has a unique fixed point
there.

It remains to show that no other \(C^\rho\)-solution exists outside
the fixed-point ball. Let \(X\) be any \(C^\rho\)-solution on \(I\).
Using the equation, boundedness of the coefficients, and
Young--Loeve's inequality, we obtain
\[
[X]_{\rho;I}
\le
C
+
CK
+
CK\ell^\rho[X]_{\rho;I}.
\]
Indeed,
\[
[\sigma(X,\cdot)]_{\rho;I}
\le
C[X]_{\rho;I}
+
C\ell^{1-\rho}.
\]
By decreasing \(\eta\) so that \(CK\ell^\rho\le C\eta\le1/2\), the
last term can be absorbed, yielding
\[
[X]_{\rho;I}
\le
C(1+K),
\]
which proves \eqref{eq:local-solution-holder-bound}. Moreover,
\[
\ell^\rho[X]_{\rho;I}
\le
C\ell^\rho(1+K)
\le
C(\eta^\rho+\eta).
\]
After one final deterministic decrease of \(\eta\), the right-hand
side is at most \(R_0\). Hence every \(C^\rho\)-solution belongs to
the corresponding ball \(\mathcal B_I(X_s)\).

Finally, let \(X\) and \(Y\) be two solutions on \(I\), possibly with
different initial values. Since
\[
\ell^\rho[X]_{\rho;I}
\vee
\ell^\rho[Y]_{\rho;I}
\le R_0,
\]
the preceding difference estimates give
\[
\|X-Y\|_{I,*}
\le
|X_s-Y_s|
+
\frac12\|X-Y\|_{I,*}.
\]
This proves \eqref{eq:local-flow-lipschitz-two} and uniqueness among
all \(C^\rho\)-solutions.

The radius \(R_0\) is independent of \(K\). In particular, a length of
the form
\[
\ell_Z
=
\min\left\{
\eta,
\left(\frac{\eta}{1+K}\right)^{1/\rho}
\right\}
\]
satisfies \eqref{eq:young-smallness-condition}, and
\(\ell_Z\) is of order \((1+K)^{-1/\rho}\).
\end{proof}

\subsubsection{Proof of Theorem~\ref{thm:young-equation-wsfbm}}\label{sec:supp-proof-thm2}
\begin{proof}
Let \(R_0\) and \(\eta\) be the deterministic constants given by
Lemma~\ref{lem:young-local-fixed-radius}. On the event
\[
[Z]_{\rho;[0,T]}<\infty,
\]
which has probability one, set
\begin{equation}\label{eq:young-localization-length}
\ell_Z
:=
\min\left\{
\eta,\,
T,\,
\left(
\frac{\eta}{1+[Z]_{\rho;[0,T]}}
\right)^{1/\rho}
\right\}.
\end{equation}
Every interval \(I\subseteq[0,T]\) of length at most \(\ell_Z\)
satisfies
\[
|I|\le\eta
\]
and
\[
[Z]_{\rho;[0,T]}|I|^\rho
\le
\eta
\frac{[Z]_{\rho;[0,T]}}
     {1+[Z]_{\rho;[0,T]}}
\le\eta.
\]
Hence the smallness conditions of
Lemma~\ref{lem:young-local-fixed-radius} hold on every such interval.

Put
\[
N_Z
:=
\left\lceil\frac{T}{\ell_Z}\right\rceil,
\qquad
s_j:=\frac{jT}{N_Z},
\qquad
j=0,\ldots,N_Z,
\]
and let
\[
I_j:=[s_j,s_{j+1}],
\qquad
\delta_Z:=\frac{T}{N_Z}.
\]
Then \(\delta_Z\le\ell_Z\). Moreover, if \(\ell_Z<T\), then
\[
N_Z<\frac{T}{\ell_Z}+1
\]
and therefore
\[
\delta_Z
=
\frac{T}{N_Z}
>
\frac{T\ell_Z}{T+\ell_Z}
\ge
\frac{\ell_Z}{2}.
\]
If \(\ell_Z=T\), then \(N_Z=1\) and \(\delta_Z=T\).

Starting with \(X_0=x_0\), apply
Lemma~\ref{lem:young-local-fixed-radius} successively on the intervals
\(I_j\), taking the terminal value on \(I_j\) as the initial value on
\(I_{j+1}\). The local fixed points concatenate to a continuous
solution on \([0,T]\). Since the solution is \(\rho\)-H\"older on
each member of a finite partition, it is globally \(\rho\)-H\"older.
Indeed, if \(u<v\) and the interval \([u,v]\) intersects \(n\) members
of the partition, then
\[
|X_v-X_u|
\le
\max_j[X]_{\rho;I_j}
\sum_{k=1}^n |J_k|^\rho,
\]
where the \(J_k\)'s are the corresponding intersections. Since
\[
\sum_{k=1}^n |J_k|^\rho
\le
n^{1-\rho}
\left(\sum_{k=1}^n |J_k|\right)^\rho
=
n^{1-\rho}|v-u|^\rho,
\]
the global H\"older seminorm is finite.

The number of intervals satisfies
\begin{equation}\label{eq:number-local-intervals}
N_Z
\le
C_T
\bigl(1+[Z]_{\rho;[0,T]}\bigr)^{1/\rho},
\end{equation}
where the constant also depends on the fixed quantities
\(\eta\) and \(\rho\). Local uniqueness from
Lemma~\ref{lem:young-local-fixed-radius}, applied successively on the
partition, gives global uniqueness. This proves existence and
uniqueness.

We now prove the stability estimate. Suppose that
\[
[Z]_{\rho;[0,T]}
\vee
[\widetilde Z]_{\rho;[0,T]}
\le M.
\]
Use the common deterministic localization length
\[
\ell_M
:=
\min\left\{
\eta,\,
T,\,
\left(\frac{\eta}{1+M}\right)^{1/\rho}
\right\}.
\]
Define
\[
N_M
:=
\left\lceil\frac{T}{\ell_M}\right\rceil,
\qquad
s_j:=\frac{jT}{N_M},
\qquad
I_j:=[s_j,s_{j+1}],
\qquad
\delta_M:=\frac{T}{N_M}.
\]
The number \(N_M\) and the length \(\delta_M\) are deterministic.
Moreover,
\[
0<c_{T,M}\le\delta_M\le\ell_M,
\]
for a deterministic \(c_{T,M}>0\). Both drivers satisfy the local
smallness conditions on every \(I_j\).

Put
\[
D:=\|Z-\widetilde Z\|_{\rho;[0,T]}
\]
and
\[
d:=X-Y.
\]
On \(I_j\), write
\[
\begin{split}
d_u
&=
d_{s_j}
+
\int_{s_j}^u
\bigl\{
\mathfrak b(X_r,r)-\mathfrak b(Y_r,r)
\bigr\}\,dr\\
&\quad+
\int_{s_j}^u
\bigl\{
\sigma(X_r,r)-\sigma(Y_r,r)
\bigr\}\,d\widetilde Z_r\\
&\quad+
\int_{s_j}^u
\sigma(X_r,r)\,d(Z-\widetilde Z)_r.
\end{split}
\]
By Lemma~\ref{lem:young-local-fixed-radius}, both \(X\) and \(Y\)
belong on \(I_j\) to their corresponding fixed-radius balls. The
drift estimate and Lemma~\ref{lem:composition-holder} therefore give
\[
\left\|
u\mapsto
\int_{s_j}^u
\bigl\{
\mathfrak b(X_r,r)-\mathfrak b(Y_r,r)
\bigr\}\,dr
\right\|_{I_j,*}
\le
C\delta_M\|d\|_{I_j,*},
\]
and
\[
\left\|
u\mapsto
\int_{s_j}^u
\bigl\{
\sigma(X_r,r)-\sigma(Y_r,r)
\bigr\}\,d\widetilde Z_r
\right\|_{I_j,*}
\le
C(1+R_0)M\delta_M^\rho
\|d\|_{I_j,*}.
\]
By the choice of \(\ell_M\), the sum of the two coefficients can be
taken to be at most \(1/2\).

For the remaining term, the scaled Young estimate gives
\[
\begin{split}
&\left\|
u\mapsto
\int_{s_j}^u
\sigma(X_r,r)\,d(Z-\widetilde Z)_r
\right\|_{I_j,*}\\
&\qquad\le
C[Z-\widetilde Z]_{\rho;I_j}\delta_M^\rho
\left(
\|\sigma(X,\cdot)\|_{\infty;I_j}
+
\delta_M^\rho
[\sigma(X,\cdot)]_{\rho;I_j}
\right).
\end{split}
\]
The boundedness of \(\sigma\), together with
\[
\delta_M^\rho[X]_{\rho;I_j}\le R_0,
\]
implies
\[
\left\|
u\mapsto
\int_{s_j}^u
\sigma(X_r,r)\,d(Z-\widetilde Z)_r
\right\|_{I_j,*}
\le
C_M D.
\]
Consequently,
\[
\|X-Y\|_{I_j,*}
\le
|X_{s_j}-Y_{s_j}|
+
\frac12\|X-Y\|_{I_j,*}
+
C_M D,
\]
and hence
\begin{equation}\label{eq:young-local-driver-stability}
\|X-Y\|_{I_j,*}
\le
C_M
\left(
|X_{s_j}-Y_{s_j}|+D
\right).
\end{equation}

Let
\[
e_j:=|X_{s_j}-Y_{s_j}|.
\]
Taking \(u=s_{j+1}\) in
\eqref{eq:young-local-driver-stability} gives
\[
e_{j+1}
\le
C_M(e_j+D).
\]
Since \(e_0=|x_0-y_0|\) and the deterministic partition has
\(N_M<\infty\) intervals, induction yields
\[
\max_{0\le j\le N_M}e_j
+
\max_{0\le j<N_M}
\|X-Y\|_{I_j,*}
\le
C_{T,M}
\left(
|x_0-y_0|+D
\right).
\]

Because \(\delta_M\) is bounded away from zero by a deterministic
constant, the scaled local norms control the ordinary local
\(C^\rho\)-norms:
\[
[X-Y]_{\rho;I_j}
\le
\delta_M^{-\rho}
\|X-Y\|_{I_j,*}
\le
C_{T,M}
\left(
|x_0-y_0|+D
\right).
\]
For arbitrary \(0\le u<v\le T\), split \([u,v]\) at the partition
points. If the resulting pieces have lengths \(a_1,\ldots,a_n\), then
\(n\le N_M\) and
\[
\sum_{k=1}^n a_k^\rho
\le
N_M^{1-\rho}|v-u|^\rho.
\]
Summing the local increment estimates therefore gives
\[
[X-Y]_{\rho;[0,T]}
\le
C_{T,M}
\left(
|x_0-y_0|+D
\right).
\]
The corresponding supremum estimate follows either from the same
partition argument or from
\[
\|X-Y\|_{\infty;[0,T]}
\le
|x_0-y_0|
+
T^\rho[X-Y]_{\rho;[0,T]}.
\]
Thus
\[
\|X-Y\|_{\rho;[0,T]}
\le
C_{T,M}
\left(
|x_0-y_0|
+
\|Z-\widetilde Z\|_{\rho;[0,T]}
\right),
\]
which proves \eqref{eq:young-solution-stability}.
\end{proof}

\subsubsection{Proof of Theorem~\ref{thm:finite-variation-equation-wsfbm}}\label{sec:supp-proof-thm3}
\begin{proof}
By Proposition~\ref{prop:gaussian-velocity}, there exists an event of
probability one on which
\[
\nu(\cdot)\in L^1([0,T])
\]
and
\[
Z_t=\sqrt b\int_0^t\nu_r\,dr,
\qquad 0\le t\le T.
\]
Fix a sample path in this event and define
\[
F(t,x)
:=
\mathfrak b(x,t)
+
\sqrt b\,\sigma(x,t)\nu_t.
\]
The map \(t\mapsto F(t,x)\) is measurable for every \(x\), and
\(x\mapsto F(t,x)\) is continuous for almost every \(t\).

Let
\[
L_{\mathfrak b}:=\|\partial_x\mathfrak b\|_\infty,
\qquad
L_\sigma:=\|\partial_x\sigma\|_\infty.
\]
Then, for all \(x,y\in\mathbb R\) and almost every \(t\),
\[
|F(t,x)-F(t,y)|
\le
\bigl(
L_{\mathfrak b}
+
\sqrt b\,L_\sigma|\nu_t|
\bigr)|x-y|.
\]
The Lipschitz coefficient
\[
t\longmapsto
L_{\mathfrak b}
+
\sqrt b\,L_\sigma|\nu_t|
\]
belongs to \(L^1([0,T])\). Moreover, since
\(\mathfrak b\) and \(\sigma\) are bounded,
\[
|F(t,x)|
\le
\|\mathfrak b\|_\infty
+
\sqrt b\,\|\sigma\|_\infty|\nu_t|,
\]
and the right-hand side is integrable on \([0,T]\), independently of
\(x\).

The Carath\'eodory existence theorem therefore yields an absolutely
continuous solution of \eqref{eq:random-ode-fv} on \([0,T]\).

If \(X\) and \(Y\) are two such solutions with the same initial value,
then
\[
|X_t-Y_t|
\le
\int_0^t
\bigl(
L_{\mathfrak b}
+
\sqrt b\,L_\sigma|\nu_r|
\bigr)
|X_r-Y_r|\,dr.
\]
Gronwall's inequality gives
\[
X_t=Y_t,
\qquad 0\le t\le T,
\]
and hence uniqueness.

Finally, since \(Z\) is absolutely continuous with
\[
\dot Z_t=\sqrt b\,\nu_t
\quad\text{for a.e. }t,
\]
the Lebesgue--Stieltjes integral satisfies
\[
\int_0^t\sigma(X_r,r)\,dZ_r
=
\sqrt b\int_0^t\sigma(X_r,r)\nu_r\,dr.
\]
Thus the integral equation \eqref{eq:general-wsfbm-equation} and the differential equation
\eqref{eq:random-ode-fv} are equivalent.
\end{proof}

\subsubsection{Auxiliary estimates for the Euler approximation}\label{sec:supp-euler-aux}

\begin{lem}\label{lem:young-flow-defect}
Assume \textup{(H)}, let \(1/2<\rho<1\), and let
\(Z\in C^\rho([0,T])\). For \(0\le s\le t\le T\), let
\(\phi_{s,t}(x)\) denote the value at time \(t\) of the solution
started from \(x\) at time \(s\), and define
\[
\Psi_{s,t}(x)
:=
x
+
\mathfrak b(x,s)(t-s)
+
\sigma(x,s)Z_{s,t}.
\]
Put
\[
K:=[Z]_{\rho;[0,T]}.
\]
There exists a deterministic constant \(C<\infty\), depending only
on \textup{(H)}, \(T\), and \(\rho\), such that
\begin{align}
\sup_{0\le s\le t\le T}
\operatorname{Lip}(\phi_{s,t})
&\le
L_T(Z)
:=
1+C\exp\left\{C(1+K)^{1/\rho}\right\},
\label{eq:flow-lip}\\
\sup_{x\in\mathbb R}
|\phi_{s,t}(x)-\Psi_{s,t}(x)|
&\le
C L_T(Z)(1+K)^2|t-s|^{2\rho}.
\label{eq:local-defect-young}
\end{align}
\end{lem}

\begin{proof}
Let \(\ell_Z\) be the localization length defined in
\eqref{eq:young-localization-length}. Every interval of length at most
\(\ell_Z\) satisfies the smallness conditions
\eqref{eq:young-smallness-condition}. Hence
\eqref{eq:local-flow-lipschitz-two} gives
\[
|\phi_{u,v}(x)-\phi_{u,v}(y)|
\le 2|x-y|
\]
whenever \(v-u\le\ell_Z\).

Fix \(0\le s<t\le T\), and divide \([s,t]\) into \(m\) intervals of
length at most \(\ell_Z\). The flow property and multiplication of the
local Lipschitz factors give
\[
\operatorname{Lip}(\phi_{s,t})\le2^m.
\]
The number of intervals may be chosen so that
\[
m
\le
1+\frac{T}{\ell_Z}
\le
C_T(1+K)^{1/\rho},
\]
where the last inequality follows from
\eqref{eq:young-localization-length}. Therefore
\[
\operatorname{Lip}(\phi_{s,t})
\le
\exp\left\{C(1+K)^{1/\rho}\right\},
\]
which proves \eqref{eq:flow-lip}, after increasing \(C\).

We next prove the one-step defect estimate. Put
\[
Y_u:=\phi_{s,u}(x),
\qquad s\le u\le t.
\]
Suppose first that
\[
t-s\le\ell_Z.
\]
By \eqref{eq:local-solution-holder-bound},
\[
[Y]_{\rho;[s,t]}
\le
C(1+K),
\]
uniformly in \(x\). Since \(Y_s=x\), the bounded derivatives in
\textup{(H)} give
\[
|\mathfrak b(Y_r,r)-\mathfrak b(x,s)|
\le
C\bigl(|Y_r-x|+|r-s|\bigr)
\le
C(1+K)|r-s|^\rho.
\]
Consequently,
\[
\left|
\int_s^t
\bigl\{
\mathfrak b(Y_r,r)-\mathfrak b(x,s)
\bigr\}\,dr
\right|
\le
C(1+K)|t-s|^{1+\rho}.
\]

Define
\[
q_r:=\sigma(Y_r,r)-\sigma(x,s).
\]
Then \(q_s=0\), and the boundedness of
\(\partial_x\sigma\) and \(\partial_t\sigma\) implies
\[
[q]_{\rho;[s,t]}
\le
C(1+K).
\]
Young--Loeve's inequality therefore gives
\[
\left|
\int_s^t
\bigl\{
\sigma(Y_r,r)-\sigma(x,s)
\bigr\}\,dZ_r
\right|
\le
C(1+K)K|t-s|^{2\rho}.
\]
Since \(\rho<1\),
\[
|t-s|^{1+\rho}
\le
T^{1-\rho}|t-s|^{2\rho}.
\]
Using
\[
\phi_{s,t}(x)-\Psi_{s,t}(x)
=
\int_s^t
\bigl\{
\mathfrak b(Y_r,r)-\mathfrak b(x,s)
\bigr\}\,dr
+
\int_s^t
\bigl\{
\sigma(Y_r,r)-\sigma(x,s)
\bigr\}\,dZ_r,
\]
we obtain
\begin{equation}\label{eq:young-local-defect-short}
\sup_{x\in\mathbb R}
|\phi_{s,t}(x)-\Psi_{s,t}(x)|
\le
C(1+K)^2|t-s|^{2\rho}
\end{equation}
whenever \(t-s\le\ell_Z\).

Suppose now that
\[
t-s>\ell_Z.
\]
Divide \([s,t]\) into intervals \(J\) of length at most \(\ell_Z\).
On each such interval, \eqref{eq:local-solution-holder-bound} gives
\[
\sup_{u,v\in J}|Y_v-Y_u|
\le
C(1+K)|J|^\rho.
\]
By the definition of \(\ell_Z\),
\[
(1+K)|J|^\rho
\le
(1+K)\ell_Z^\rho
\le
\eta+\eta^\rho,
\]
so the increment on each local interval is bounded by a deterministic
constant. Since the number of intervals is at most
\(C_T(1+K)^{1/\rho}\), it follows that
\[
\sup_{x\in\mathbb R}
|\phi_{s,t}(x)-x|
\le
C_T(1+K)^{1/\rho}.
\]
The boundedness of \(\mathfrak b\) and \(\sigma\) also gives
\[
\begin{split}
\sup_{x\in\mathbb R}
|\Psi_{s,t}(x)-x|
&\le
\|\mathfrak b\|_\infty T
+
\|\sigma\|_\infty K T^\rho\\
&\le
C_T(1+K)
\le
C_T(1+K)^{1/\rho}.
\end{split}
\]
Hence
\begin{equation}\label{eq:young-local-defect-coarse}
\sup_{x\in\mathbb R}
|\phi_{s,t}(x)-\Psi_{s,t}(x)|
\le
C_T(1+K)^{1/\rho}.
\end{equation}

We now use the assumption \(t-s>\ell_Z\). Since \(t-s\le T\), the
minimum defining \(\ell_Z\) cannot be attained at \(T\). If
\(\ell_Z=\eta\), then
\[
(1+K)^2|t-s|^{2\rho}
\ge
\eta^{2\rho}.
\]
If instead
\[
\ell_Z
=
\left(
\frac{\eta}{1+K}
\right)^{1/\rho},
\]
then
\[
(1+K)^2|t-s|^{2\rho}
>
(1+K)^2\ell_Z^{2\rho}
=
\eta^2.
\]
Thus, in either case,
\[
(1+K)^2|t-s|^{2\rho}
\ge
c_{\eta,\rho}>0.
\]
Moreover,
\[
(1+K)^{1/\rho}
\le
C L_T(Z).
\]
Combining these inequalities with
\eqref{eq:young-local-defect-coarse} gives
\[
\sup_{x\in\mathbb R}
|\phi_{s,t}(x)-\Psi_{s,t}(x)|
\le
C L_T(Z)(1+K)^2|t-s|^{2\rho}.
\]
Together with \eqref{eq:young-local-defect-short}, this proves
\eqref{eq:local-defect-young} for every \(0\le s\le t\le T\).
\end{proof}

\begin{lem}\label{lem:gaussian-holder-fernique}
Let \(Z\) be a centered Gaussian process with \(Z_0=0\) and
\[
\mathbb E|Z_t-Z_s|^2
\le
C|t-s|^{2H},
\qquad
0\le s,t\le T,
\]
for some \(0<H\le1\). For every \(0<\rho<H\), \(Z\) admits a
modification such that
\begin{equation}\label{eq:fernique-holder}
\mathbb E\exp\bigl\{\alpha[Z]_{\rho;[0,T]}^2\bigr\}<\infty
\end{equation}
for some \(\alpha>0\).
\end{lem}

\begin{proof}
Fix \(\rho<\eta<H\). For every \(p\ge2\), Gaussian moment
equivalence gives
\[
\mathbb E|Z_t-Z_s|^p
\le
C_p|t-s|^{pH},
\qquad
0\le s,t\le T.
\]
Choose \(p\) sufficiently large that
\[
\eta<H-\frac1p.
\]
Kolmogorov's continuity theorem then yields a modification, still
denoted by \(Z\), whose sample paths belong to
\(C^\eta([0,T])\). In particular,
\[
[Z]_{\rho;[0,T]}<\infty
\qquad\text{almost surely}.
\]

Regard this modification as a centered Gaussian random element of the
separable Banach space
\[
E:=C([0,T])
\]
equipped with the uniform norm. Define the extended-valued functional
\[
N(x):=[x]_{\rho;[0,T]},
\qquad x\in E.
\]
Let \(D\subset[0,T]\) be a countable dense set containing \(0\) and
\(T\). For every continuous \(x\),
\[
N(x)
=
\sup_{\substack{r,q\in D\\r\ne q}}
\frac{|x(r)-x(q)|}{|r-q|^\rho}.
\]
Each term in this countable supremum is continuous on \(E\), and hence
\(N\) is Borel measurable. Moreover,
\[
N(x+y)\le N(x)+N(y),
\qquad
N(\lambda x)=|\lambda|N(x),
\]
and \(\{x\in E:N(x)<\infty\}\) is a vector subspace of \(E\).
Thus \(N\) is a measurable pseudo-seminorm.

Since
\[
\mathbb P\{N(Z)<\infty\}=1,
\]
Fernique's integrability theorem for measurable pseudo-seminorms
\cite[Theorem~1.3.2, p.~11]{Fernique1975} yields an
\(\alpha>0\) such that
\[
\mathbb E\exp\{\alpha N(Z)^2\}<\infty.
\]
This is precisely \eqref{eq:fernique-holder}.
\end{proof}

\subsubsection{Proof of Theorem~\ref{thm:euler-wsfbm}}\label{sec:supp-proof-thm4}
\begin{proof}
We first prove part \textup{(i)}. The exact solutions satisfy the flow
property
\[
\phi_{r,t}\circ\phi_{s,r}=\phi_{s,t},
\qquad 0\le s\le r\le t\le T.
\]
At the grid points, a telescoping decomposition of the exact and Euler
flows gives
\begin{align*}
|X_{t_k}-X_{t_k}^\pi|
&\le
\sum_{j=0}^{k-1}
\operatorname{Lip}(\phi_{t_{j+1},t_k})
\sup_{x\in\mathbb R}
|\phi_{t_j,t_{j+1}}(x)-\Psi_{t_j,t_{j+1}}(x)|.
\end{align*}
By Lemma~\ref{lem:young-flow-defect},
\[
\operatorname{Lip}(\phi_{t_{j+1},t_k})
\le L_T(Z)
\]
and
\[
\sup_{x\in\mathbb R}
|\phi_{t_j,t_{j+1}}(x)-\Psi_{t_j,t_{j+1}}(x)|
\le
CL_T(Z)
\bigl(1+[Z]_{\rho;[0,T]}\bigr)^2
\Delta_j^{2\rho}.
\]
Consequently,
\begin{align*}
|X_{t_k}-X_{t_k}^\pi|
&\le
CL_T(Z)^2
\bigl(1+[Z]_{\rho;[0,T]}\bigr)^2
\sum_{j=0}^{k-1}\Delta_j^{2\rho}\\
&\le
CL_T(Z)^2
\bigl(1+[Z]_{\rho;[0,T]}\bigr)^2
h^{2\rho-1}
\sum_{j=0}^{k-1}\Delta_j\\
&\le
C_TL_T(Z)^2
\bigl(1+[Z]_{\rho;[0,T]}\bigr)^2
h^{2\rho-1}.
\end{align*}

Let \(t\in[t_k,t_{k+1}]\). Since
\[
X_t=\phi_{t_k,t}(X_{t_k}),
\qquad
X_t^\pi=\Psi_{t_k,t}(X_{t_k}^\pi),
\]
we have
\begin{align*}
|X_t-X_t^\pi|
&\le
|\phi_{t_k,t}(X_{t_k})
-\phi_{t_k,t}(X_{t_k}^\pi)|\\
&\quad+
|\phi_{t_k,t}(X_{t_k}^\pi)
-\Psi_{t_k,t}(X_{t_k}^\pi)|\\
&\le
\operatorname{Lip}(\phi_{t_k,t})
|X_{t_k}-X_{t_k}^\pi|
+
\sup_x|\phi_{t_k,t}(x)-\Psi_{t_k,t}(x)|.
\end{align*}
Applying Lemma~\ref{lem:young-flow-defect} once more gives
\begin{align*}
|X_t-X_t^\pi|
&\le
C_TL_T(Z)^3
\bigl(1+[Z]_{\rho;[0,T]}\bigr)^2h^{2\rho-1}\\
&\quad+
CL_T(Z)
\bigl(1+[Z]_{\rho;[0,T]}\bigr)^2h^{2\rho}.
\end{align*}
Since \(h\le T\) and \(L_T(Z)\ge1\), the second term is absorbed into
the first one after changing \(C_T\). By the definition of \(L_T(Z)\),
\[
L_T(Z)^3
\le
C_T\exp\!\left\{
C_T\bigl(1+[Z]_{\rho;[0,T]}\bigr)^{1/\rho}
\right\}.
\]
After enlarging \(C_T\) once more, the resulting estimate is exactly \eqref{eq:euler-young-pathwise-final}.

By \eqref{eq:rough-driver-increment-bound},
\[
\mathbb E|Z_t-Z_s|^2
\le
C|t-s|^{1+b}.
\]
Lemma~\ref{lem:gaussian-holder-fernique}, applied with
\[
H=\frac{1+b}{2},
\]
yields
\[
\mathbb E
\exp\bigl\{\alpha[Z]_{\rho;[0,T]}^2\bigr\}
<\infty
\]
for some \(\alpha>0\). Since \(\rho>1/2\),
\[
\frac1\rho<2.
\]
Moreover, by the definition of \(L_T(Z)\),
\[
L_T(Z)^3
\bigl(1+[Z]_{\rho;[0,T]}\bigr)^2
\le
C
\bigl(1+[Z]_{\rho;[0,T]}\bigr)^2
\exp\left\{
C\bigl(1+[Z]_{\rho;[0,T]}\bigr)^{1/\rho}
\right\}.
\]
For every \(p\ge1\), there exists \(C_p<\infty\) such that
\[
pC(1+x)^{1/\rho}
+
2p\log(1+x)
\le
\frac{\alpha}{2}x^2+C_p,
\qquad x\ge0.
\]
It follows that
\[
\mathbb E\left[
L_T(Z)^3
\bigl(1+[Z]_{\rho;[0,T]}\bigr)^2
\right]^p
<\infty.
\]
Taking \(L^p\)-norms in the pathwise estimate \eqref{eq:euler-young-pathwise-final} proves
\eqref{eq:euler-lp-young}.

We now prove part \textup{(ii)}. On the probability-one event on which
\(\nu\in L^1([0,T])\), write
\[
dA_t
=
\bigl(1+\sqrt b\,|\nu_t|\bigr)\,dt.
\]
For the Carath\'eodory ODE
\[
\dot X_t
=
\mathfrak b(X_t,t)
+
\sqrt b\,\sigma(X_t,t)\nu_t,
\]
the spatial Lipschitz coefficient is bounded by
\[
C\bigl(1+\sqrt b\,|\nu_t|\bigr).
\]
Gronwall's inequality therefore gives
\begin{equation}\label{eq:fv-flow-lipschitz}
\operatorname{Lip}(\phi_{s,t})
\le
\exp\{C(A_t-A_s)\},
\qquad 0\le s\le t\le T.
\end{equation}

Let \(X_r=\phi_{s,r}(x)\). Since the coefficients are bounded,
\[
|X_r-x|
\le
C\int_s^r
\bigl(1+\sqrt b\,|\nu_q|\bigr)\,dq
=
C(A_r-A_s).
\]
The one-step defect of the Euler map is
\[
\begin{split}
\phi_{s,t}(x)-\Psi_{s,t}(x)
&=
\int_s^t
\{\mathfrak b(X_r,r)-\mathfrak b(x,s)\}\,dr\\
&\quad+
\sqrt b\int_s^t
\{\sigma(X_r,r)-\sigma(x,s)\}\nu_r\,dr.
\end{split}
\]
By the bounded first derivatives in \textup{(H)},
\[
|\mathfrak b(X_r,r)-\mathfrak b(x,s)|
+
|\sigma(X_r,r)-\sigma(x,s)|
\le
C\bigl(|X_r-x|+r-s\bigr).
\]
Since
\[
r-s\le A_r-A_s,
\]
we obtain
\begin{align}
\sup_{x\in\mathbb R}
|\phi_{s,t}(x)-\Psi_{s,t}(x)|
&\le
C\int_s^t(A_r-A_s)\,dA_r \notag\\
&=
\frac C2(A_t-A_s)^2
\le
C(A_t-A_s)^2.
\label{eq:fv-one-step-defect}
\end{align}

At a grid point \(t_k\), the same telescoping-flow decomposition used
above, together with
\eqref{eq:fv-flow-lipschitz} and
\eqref{eq:fv-one-step-defect}, gives
\begin{align*}
|X_{t_k}-X_{t_k}^\pi|
&\le
Ce^{CA_T}
\sum_{j=0}^{k-1}
(A_{t_{j+1}}-A_{t_j})^2.
\end{align*}
Because \(A\) is non-decreasing and
\[
t_{j+1}-t_j\le h,
\]
we have
\[
A_{t_{j+1}}-A_{t_j}\le\omega_A(h).
\]
Therefore,
\begin{align*}
\sum_{j=0}^{k-1}
(A_{t_{j+1}}-A_{t_j})^2
&\le
\omega_A(h)
\sum_{j=0}^{k-1}(A_{t_{j+1}}-A_{t_j})\\
&\le
A_T\omega_A(h),
\end{align*}
and hence
\[
\max_k|X_{t_k}-X_{t_k}^\pi|
\le
Ce^{CA_T}A_T\omega_A(h).
\]

For \(t\in[t_k,t_{k+1}]\),
\[
\begin{split}
|X_t-X_t^\pi|
&\le
\operatorname{Lip}(\phi_{t_k,t})
|X_{t_k}-X_{t_k}^\pi|\\
&\quad+
\sup_x|\phi_{t_k,t}(x)-\Psi_{t_k,t}(x)|.
\end{split}
\]
The first term is bounded by
\[
Ce^{CA_T}A_T\omega_A(h).
\]
For the second term,
\[
\sup_x|\phi_{t_k,t}(x)-\Psi_{t_k,t}(x)|
\le
C(A_t-A_{t_k})^2
\le
C\omega_A(h)^2.
\]
Since
\[
\omega_A(h)\le A_T,
\]
this term is bounded by \(CA_T\omega_A(h)\). After changing the
deterministic constant, we obtain
\[
\sup_{0\le t\le T}|X_t-X_t^\pi|
\le
C_Te^{C_TA_T}A_T\omega_A(h),
\]
which proves \eqref{eq:euler-fv-modulus}. Since \(A\) is absolutely
continuous, it is continuous on \([0,T]\), and therefore
\[
\omega_A(h)\longrightarrow0
\qquad\text{as }h\downarrow0.
\]

Finally, assume that \(a\ge0\). By
Proposition~\ref{prop:gaussian-velocity}, \(\nu\) admits a continuous
Gaussian modification. For this modification,
\[
\omega_A(h)
\le
\bigl(1+\sqrt b\,\|\nu\|_{\infty;[0,T]}\bigr)h
\]
and
\[
A_T
\le
T\bigl(1+\sqrt b\,\|\nu\|_{\infty;[0,T]}\bigr).
\]
Consequently,
\[
e^{C_TA_T}A_T\omega_A(h)
\le
C_T
\bigl(1+\|\nu\|_{\infty;[0,T]}\bigr)^2
\exp\bigl\{C_T\|\nu\|_{\infty;[0,T]}\bigr\}h.
\]
The random variable \(\nu\), regarded as a centered Gaussian random
element of \(C([0,T])\), satisfies Fernique's theorem. Hence there
exists \(\beta>0\) such that
\[
\mathbb E
\exp\bigl\{\beta\|\nu\|_{\infty;[0,T]}^2\bigr\}
<\infty.
\]
It follows that the random factor in the preceding estimate belongs to
\(L^p\) for every \(p\ge1\). Taking \(L^p\)-norms in
\eqref{eq:euler-fv-modulus} proves
\eqref{eq:euler-lp-fv}.
\end{proof}

\subsection{Related covariance constructions on Euclidean index sets}\label{sec:euclidean-covariances}
The ordered half-line geometry in \eqref{weighted-covariance}, in particular the term \(s+t-2u\), has no direct rotation-invariant analogue on an unordered Euclidean index set. The constructions below are spatial covariance models based on the same positive-definite weighting mechanism; they are not multidimensional versions of \(R_{f,b}\). Throughout this section, \(d\in\mathbb N\), \(\|\cdot\|\) denotes the Euclidean norm on \(\mathbb R^d\), and \(du\) denotes Lebesgue measure. Let \(E\subseteq\mathbb R^d\) be measurable, let \(A\subseteq E\) be non-empty and closed, and write
\[
d_A(x)=\operatorname{dist}(x,A).
\]
Let \(f:E\to[0,\infty)\) be measurable. For \(u,x\in E\), define
\[
\chi_x(u)=\mathbf 1_{\{0<d_A(u)\le d_A(x)\}}.
\]
Then \(\chi_x(u)\chi_y(u)=\mathbf 1_{\{0<d_A(u)\le d_A(x)\wedge d_A(y)\}}\).

For \(0<H\le1\), set
\[
Q_{H,-}(x,y)=\|x\|^{2H}+\|y\|^{2H}-\|x-y\|^{2H}.
\]
This is twice the covariance of the L\'evy fractional Brownian field and is positive definite on \(\mathbb R^d\). Assume
\begin{equation}\label{eq:spatial-integrability}
\int_E\chi_x(u)f(u)\|x-u\|^{2H}du<\infty,
\qquad x\in E.
\end{equation}
Define
\begin{equation}\label{eq:spatial-levy-weighted}
K_{H,A,f,-}(x,y)=\int_E\chi_x(u)\chi_y(u)f(u)
Q_{H,-}(x-u,y-u)\,du.
\end{equation}

\begin{proposition}\label{prop:spatial-pd}
Under \eqref{eq:spatial-integrability}, \(K_{H,A,f,-}\) is finite and positive definite.
\end{proposition}

\begin{proof}
For fixed \(u\), the translated kernel
\(Q_{H,-}(x-u,y-u)\) is positive definite, and
\(\chi_x(u)\chi_y(u)\) is a rank-one positive definite kernel. Their product is positive definite by the Schur product theorem. Before integrating, note the pointwise covariance inequality
\[
|Q_{H,-}(x-u,y-u)|
\le 2\|x-u\|^H\|y-u\|^H.
\]
Consequently, weighted Cauchy--Schwarz and \eqref{eq:spatial-integrability} imply
\begin{align*}
&\int_E\chi_x(u)\chi_y(u)f(u)|Q_{H,-}(x-u,y-u)|du\\
&\qquad\le2\left(\int_E\chi_x(u)f(u)\|x-u\|^{2H}du\right)^{1/2}
\left(\int_E\chi_y(u)f(u)\|y-u\|^{2H}du\right)^{1/2}<\infty.
\end{align*}
Thus every entry is absolutely integrable, and finite-sum interchange is legitimate. For finite \(x_i\) and coefficients \(c_i\),
\[
\sum_{i,j}c_ic_jK_{H,A,f,-}(x_i,x_j)
=\int_Ef(u)\sum_{i,j}c_ic_j\chi_{x_i}(u)\chi_{x_j}(u)
Q_{H,-}(x_i-u,x_j-u)du\ge0.
\]
The diagonal is
\[
K_{H,A,f,-}(x,x)=2\int_E\chi_x(u)f(u)\|x-u\|^{2H}du<\infty,
\]
so the kernel is finite and positive definite.
\end{proof}

The superficially similar expression
\(Q_{H,+}(x,y)=\|x\|^{2H}+\|y\|^{2H}-\|x+y\|^{2H}\)
is not positive definite on the full Euclidean space. If \(H>1/2\), its diagonal is negative. If \(0<H\le1/2\), the covariance matrix at \(x\) and \(-x\) has diagonal
\((2-2^{2H})\|x\|^{2H}\) and off-diagonal \(2\|x\|^{2H}\), hence negative determinant. This explains why the temporal plus sign does not extend directly to an unordered Euclidean index set.

A second construction starts from any finite positive definite kernel \(K\) on \(E\). Assume
\begin{equation}\label{eq:gram-integrability}
\int_E\chi_x(u)f(u)du<\infty,
\qquad x\in E,
\end{equation}
and define
\begin{align}
C_{A,f}(x,y)&=\int_E\chi_x(u)\chi_y(u)f(u)du,\label{eq:gram-covariance}\\
\widetilde K_{A,f}(x,y)&=K(x,y)C_{A,f}(x,y).\label{eq:schur-spatial-covariance}
\end{align}
The map \(x\mapsto\chi_x\sqrt f\) is a feature map into \(L^2(E,du)\), so \(C_{A,f}\) is positive definite. The Schur product theorem then proves that \(\widetilde K_{A,f}\) is positive definite.

To recover the time-changed Brownian covariance exactly, take the base space \(E=\mathbb R_+\), \(A=\{0\}\), and \(x,y\ge0\). Then \(d_A(u)=u\), and
\[
C_{A,f}(x,y)=\int_0^{x\wedge y}f(u)du.
\]
The restriction to the half-line is essential; integrating over the full real line would include both signs of \(u\).

\begin{example}
Let \(E=\mathbb R^d\), \(A=\overline B(0,1):=\{z\in\mathbb R^d:\|z\|\le1\}\), and \(R_{x,y}:=\|x\|\wedge\|y\|\). Since \(d_A(u)=(\|u\|-1)_+\), the integration region in \eqref{eq:gram-covariance} is \(1<\|u\|\le R_{x,y}\) when \(R_{x,y}>1\), and is empty otherwise. Writing \(\omega_{d-1}=2\pi^{d/2}/\Gamma(d/2)\), we obtain
\[
C_{A,f}(x,y)=\omega_{d-1}\int_1^{R_{x,y}} f_0(r)r^{d-1}dr
\]
for radial weights \(f(u)=f_0(\|u\|)\). In particular,
\[
f_0(r)=r^a:\quad
C_{A,f}(x,y)=\frac{\omega_{d-1}}{a+d}(R_{x,y}^{a+d}-1),
\]
when \(a\ne-d\), with the logarithmic value \(\omega_{d-1}\log R_{x,y}\) when \(a=-d\). For \(f_0(r)=e^{ar}\),
\[
C_{A,f}(x,y)=\omega_{d-1}\int_1^{R_{x,y}}e^{ar}r^{d-1}dr.
\]
These formulas are understood as zero when \(R_{x,y}\le1\).
\end{example}

Typical choices for \(K\) in \eqref{eq:schur-spatial-covariance} include the Mat\'ern, squared-exponential, rational-quadratic, and product-periodic kernels:
\begin{align*}
K_M(x,y)&=\tau^2\frac{2^{1-\mu}}{\Gamma(\mu)}
\left(\frac{\sqrt{2\mu}\|x-y\|}{\ell_M}\right)^\mu
\mathcal K_\mu\left(\frac{\sqrt{2\mu}\|x-y\|}{\ell_M}\right),\\
K_{SE}(x,y)&=\tau^2\exp\left(-\frac{\|x-y\|^2}{2\ell_{SE}^2}\right),\\
K_{RQ}(x,y)&=\tau^2\left(1+\frac{\|x-y\|^2}{2q\ell_{RQ}^2}\right)^{-q},\\
K_P(x,y)&=\tau^2\exp\left[-\sum_{j=1}^d
\frac{2\sin^2(\pi(x_j-y_j)/p_j)}{\ell_j^2}\right].
\end{align*}
Here \(\mathcal K_\mu\) is the modified Bessel function of the second kind; \(\tau,\mu,\ell_M,\ell_{SE},q,\ell_{RQ},p_j,\ell_j\) are positive parameters. At \(x=y\), the Mat\'ern expression is interpreted by its continuous limit \(K_M(x,x)=\tau^2\). Each displayed kernel is positive definite; hence its product with \(C_{A,f}\) is a valid covariance.

\subsection{Numerical evaluation of the covariance matrices}\label{sec:supp-numerical}

We compare four methods for evaluating \eqref{supp-weighted-covariance} and \eqref{supp-log-covariance} for \(\mathcal C_1\) and \(\mathcal C_2\). The implementation uses four strategies. Method~1 is adaptive Gauss--Kronrod quadrature as implemented by \texttt{integrate} in \texttt{R}. Method~2 is the h-adaptive Genz--Malik cubature rule implemented by \texttt{hcubature} \cite{GenzMalik1980}. Method~3 is the p-adaptive Clenshaw--Curtis rule implemented by \texttt{pcubature}. Method~4 evaluates the closed forms below, using standard implementations of the beta, incomplete beta, confluent hypergeometric, Gauss hypergeometric, digamma, and exponential-integral functions. The Gauss--Kronrod methodology is described in \cite{Piessens1983}. We use the standard notation \({}_1F_1\) for the confluent hypergeometric function, \({}_2F_1\) for the Gauss hypergeometric function, \(\psi=\Gamma'/\Gamma\) for the digamma function, and \(\operatorname{Ei}\) for the exponential integral. For \(p,q>0\) and \(0\le x\le1\), we write \(B(p,q)=\int_0^1v^{p-1}(1-v)^{q-1}\,dv\) and \(B_x(p,q)=\int_0^xv^{p-1}(1-v)^{q-1}\,dv\). The symbols \({}_1F_1\) and \({}_2F_1\) refer to their standard real branches on the arguments used below.

Put \(m=s\wedge t\) and \(M=s\vee t\). For \(f(u)=u^a\), \(a>-1\), and \(b>-1\), \(b\ne1\), direct changes of variables give
\begin{equation}\label{supp:power-closed}
\begin{aligned}
R_{a,b}^{(1)}(s,t)=\frac{1}{1-b}\Bigg[&m^{a+b+1}B(a+1,b+1)
+M^{a+b+1}B_{m/M}(a+1,b+1)\\
&-2^{-a-1}(M+m)^{a+b+1}B_{2m/(M+m)}(a+1,b+1)\Bigg],
\end{aligned}
\end{equation}
where \(B_x\) is the incomplete Beta function defined above. Here \(R_{a,b}^{(1)}=R_{f_a^{(1)},b}\). Formula \eqref{supp:power-closed} is interpreted by continuity when \(m=0\).

For the exponential family, define, for \(z\ge0\), \(c\ge m\), and \(b>-1\),
\begin{align}
\mathcal F_b(z;a)&=\frac{z^{b+1}}{b+1}\,{}_1F_1(b+1;b+2;-az),\label{supp:Fdef}\\
\mathcal J_{a,b}(c,m)&=e^{ac}\bigl[\mathcal F_b(c;a)-\mathcal F_b(c-m;a)\bigr].\label{supp:Jdef}
\end{align}
Since \(\frac{d}{dz}\mathcal F_b(z;a)=z^be^{-az}\), one has
\(
\mathcal J_{a,b}(c,m)=\int_0^m e^{ar}(c-r)^b\,dr
\)
for every real \(a\). Consequently, for \(b\ne1\),
\begin{equation}\label{supp:exp-closed}
R_{a,b}^{(2)}(s,t)=\frac{\mathcal J_{a,b}(s,m)+\mathcal J_{a,b}(t,m)-2^b\mathcal J_{a,b}((s+t)/2,m)}{1-b}.
\end{equation}
Here \(R_{a,b}^{(2)}=R_{f_a^{(2)},b}\). Formula~\eqref{supp:exp-closed} is real valued for every \(a\in\mathbb R\), including \(a=0\).

For the logarithmic covariance with power weight, let \(a>-1\) and set
\begin{equation}\label{supp:Hdef}
\begin{split}
H_a(x)&=\int_0^x v^a(1-v)\log(1-v)\,dv\\
&=-\frac{x^{a+1}}{(a+1)(a+2)(a+3)}\Big[(a+1)x^2{}_2F_1(1,a+3;a+4;x)\\
&\hspace{3.5em}-(a+3)x{}_2F_1(1,a+2;a+3;x)\\
&\hspace{3.5em}+(a+3)((a+1)x-a-2)\log(1-x)\Big].
\end{split}
\end{equation}
Formula~\eqref{supp:Hdef} holds for \(0\le x<1\) and is extended continuously to \(x=1\). Moreover,
\[
H_a(1)=B(a+1,2)\{\psi(2)-\psi(a+3)\},
\]
where \(\psi\) is the digamma function. It follows that
\begin{equation}\label{supp:power-log-closed}
\begin{aligned}
K_a^{(1)}(s,t)=&\;2^{-a-1}(M+m)^{a+2}\left[\log(M+m)B_{2m/(M+m)}(a+1,2)+H_a\!\left(\frac{2m}{M+m}\right)\right]\\
&-M^{a+2}\left[\log M\,B_{m/M}(a+1,2)+H_a\!\left(\frac mM\right)\right]\\
&-m^{a+2}\left[\log m\,B(a+1,2)+H_a(1)\right],
\end{aligned}
\end{equation}
Here \(K_a^{(1)}=K_{f_a^{(1)}}\), and every term containing a zero base is interpreted by continuity.

For the logarithmic covariance with exponential weight, define, for \(y\ge0\),
\begin{equation}\label{supp:Phi}
\Phi_a(y)=
\begin{cases}
\displaystyle\frac{\operatorname{Ei}(-ay)-e^{-ay}\{(ay+1)\log y+1\}}{a^2},&a\ne0,\ y>0,\\[1.2ex]
\displaystyle\frac{\gamma+\log|a|-1}{a^2},&a\ne0,\ y=0,\\[1.2ex]
\displaystyle\frac{y^2}{2}\log y-\frac{y^2}{4},&a=0,\ y>0,\\[1.2ex]
0,&a=0,\ y=0,
\end{cases}
\end{equation}
where \(\gamma\) is Euler's constant. For \(c\ge m\ge0\), put
\begin{equation}\label{supp:Ldef}
\mathcal L_a(c,m)=e^{ac}\{\Phi_a(c)-\Phi_a(c-m)\}
=\int_0^m e^{ar}(c-r)\log(c-r)\,dr,
\end{equation}
and therefore
\begin{equation}\label{supp:exp-log-closed}
K_a^{(2)}(s,t)=2\mathcal L_a\!\left(\frac{s+t}{2},m\right)
+2\log 2\,\mathcal J_{a,1}\!\left(\frac{s+t}{2},m\right)
-\mathcal L_a(s,m)-\mathcal L_a(t,m).
\end{equation}
Here \(K_a^{(2)}=K_{f_a^{(2)}}\); the formula is interpreted by continuity when \(m=0\).
For the tested parameter values, the closed forms \eqref{supp:power-closed},
\eqref{supp:exp-closed}, \eqref{supp:power-log-closed}, and
\eqref{supp:exp-log-closed} agreed with direct quadrature.  These identities
are analytic identities, but their literal floating-point evaluation is not
uniformly stable near \(b=1\), and the representation
\(\mathcal L_a(c,m)=e^{ac}\{\Phi_a(c)-\Phi_a(c-m)\}\) also suffers cancellation
when \(a\) is close to zero.  The inference calculations in
Section~\ref{sec:supp-inference} therefore use the corrected evaluator supplied
with the accompanying reproducibility files.  It evaluates the logarithmic power integral through
stable endpoint series, evaluates the exponential logarithmic covariance near
\(a=0\) through its convergent power series in the weight parameter, and
joins the \(b\ne1\) formulas to their exact logarithmic value by a matched
local expansion in \(b-1\).  Direct quadrature of the stabilized kernel
is included as an independent check.  Gaussian likelihoods and conditional
predictions are computed by Cholesky solves; no explicit covariance inverse is
formed.

Let \(T>0\), \(n\in\mathbb N\), and \(t_k=kT/n\), \(k=0,\ldots,n\). In the timing experiments, \(T=10\) and \(n=100,200,\ldots,1500\).

\begin{figure}[htbp]
\centering
\IfFourLegacyFiguresExist{imagen1times1.png}{times2_1.png}{times3_1.png}{timesratio1.pdf}{%
\subfigure[$(f(u):=u^{0.21},b=1.28)$]{\includegraphics[width=60mm,height=4.5cm]{imagen1times1.png}}
\subfigure[$(f(u):=e^{0.21u},b=1.28)$]{\includegraphics[width=60mm,height=4.55cm]{times2_1.png}}
\subfigure[$(f(u):=e^{-0.6u},b=1.28)$]{\includegraphics[width=60mm,height=4.5cm]{times3_1.png}}
\subfigure[]{\includegraphics[width=60mm,height=4.55cm]{timesratio1.pdf}}
}{%
\includegraphics[width=0.92\linewidth]{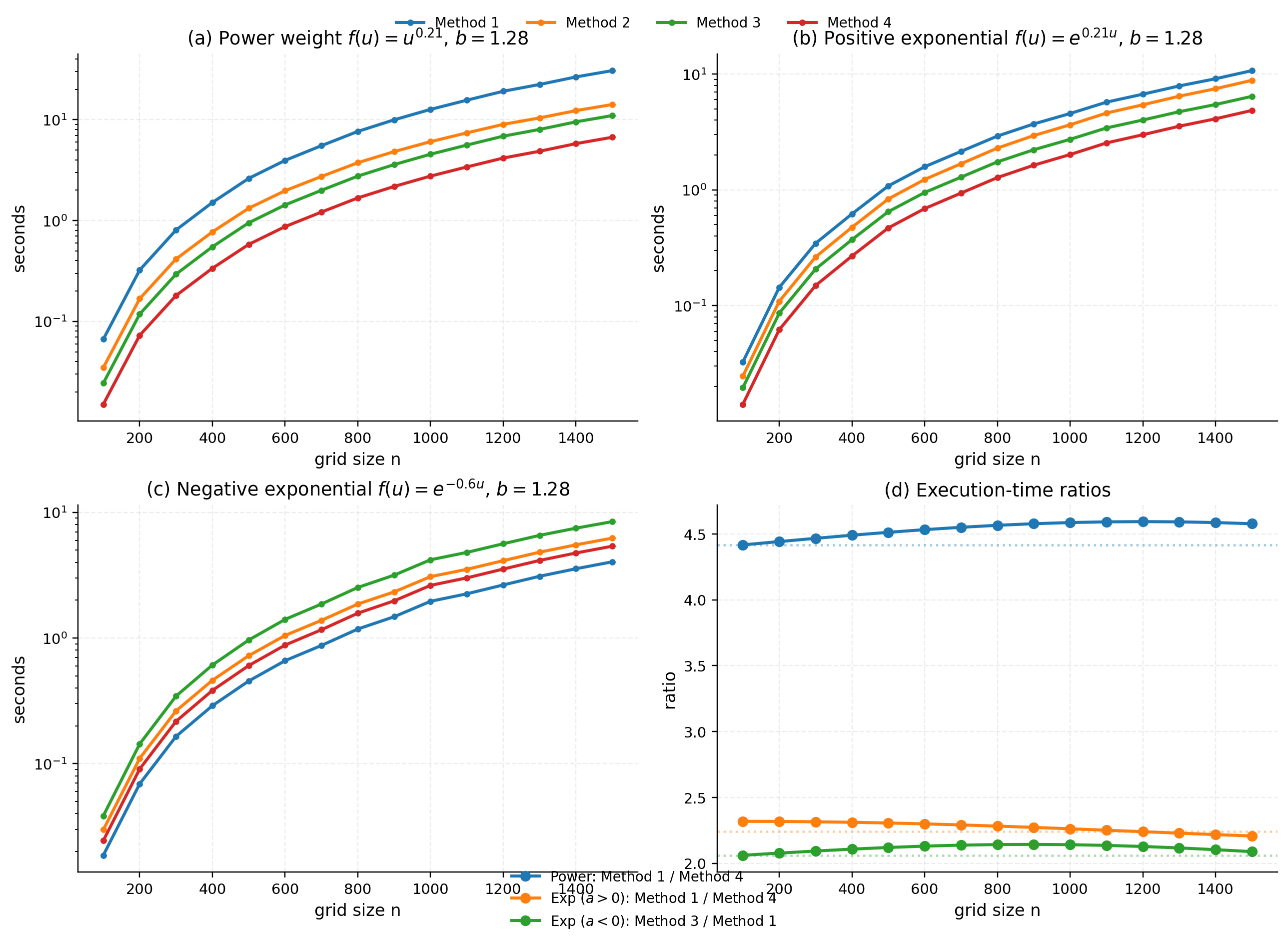}
}
\caption{Execution times for evaluating \(R_{f,b}\). Panels (a)--(c) correspond respectively to a power weight, a positive exponential weight, and a negative exponential weight. Panel (d) reports the execution-time ratios used in the text.}\label{fig:timing-power}
\end{figure}

For the power family in Figure~\ref{fig:timing-power}(a), Method~1 required on average about \(4.417\) times the execution time of Method~4. For the exponential family with \(a>0\), Method~1 required about \(2.240\) times the time of Method~4; see Figure~\ref{fig:timing-power}(b,d). For \(a<0\), Method~3 required about \(2.061\) times the execution time of Method~1, so Method~1 was the fastest implementation in that experiment. Across the four implementations, the maximum coordinate-wise discrepancies between covariance matrices were of order \(10^{-4}\).

\begin{figure}[htbp]
\centering
\IfFourLegacyFiguresExist{times4.png}{times5.png}{times6.png}{timesratio2.pdf}{%
\subfigure[$(f(u):=u^{0.21},\ K_f(s,t))$]{\includegraphics[width=60mm,height=4.2cm]{times4.png}}
\subfigure[$(f(u):=e^{0.21u},\ K_f(s,t))$]{\includegraphics[width=60mm,height=4.2cm]{times5.png}}
\subfigure[$(f(u):=e^{-0.6u},\ K_f(s,t))$]{\includegraphics[width=60mm,height=4.2cm]{times6.png}}
\subfigure[]{\includegraphics[width=60mm,height=4.2cm]{timesratio2.pdf}}
}{%
\includegraphics[width=0.92\linewidth]{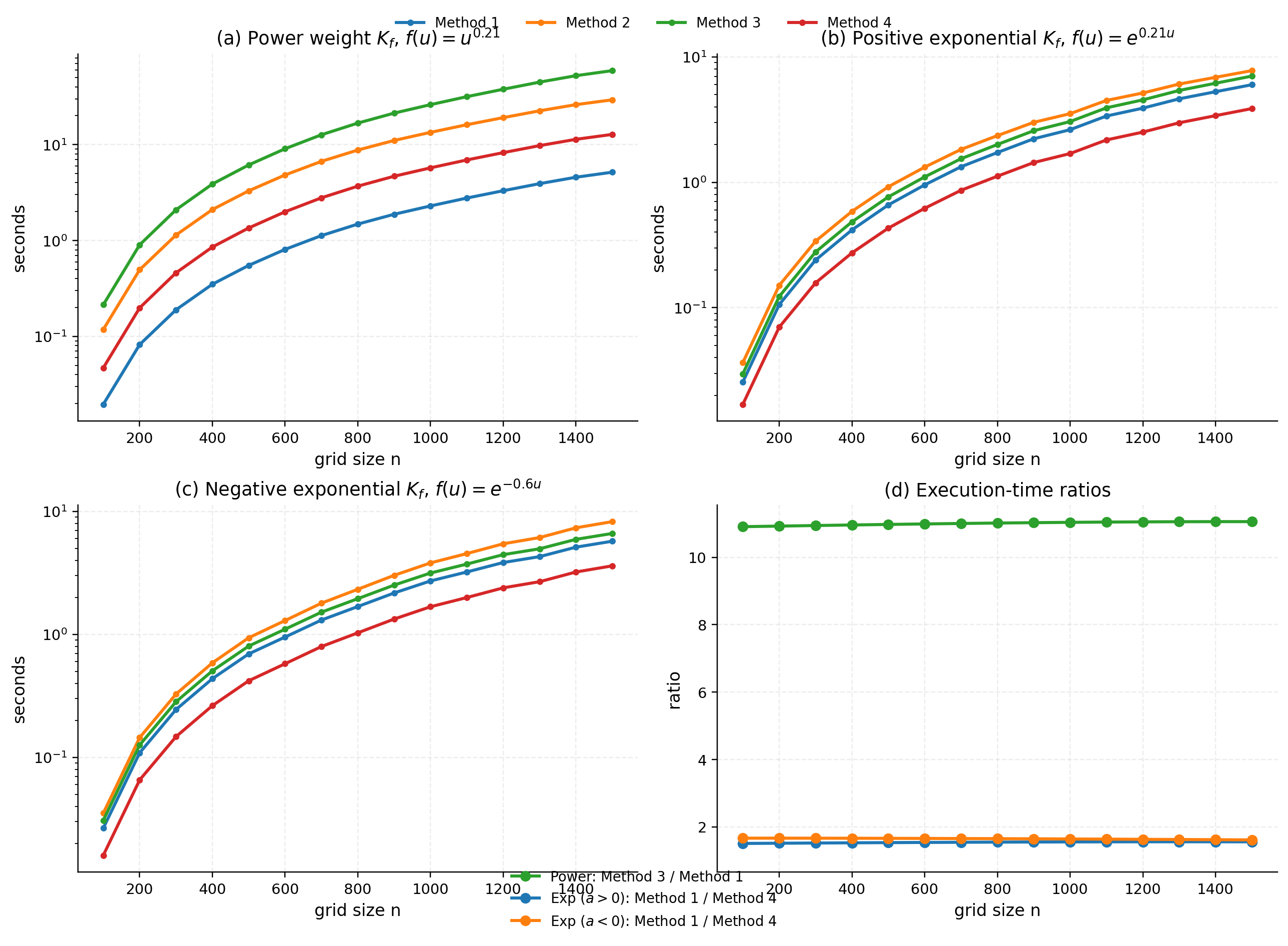}
}
\caption{Execution times for evaluating \(K_f\). Panels (a)--(c) correspond respectively to a power weight, a positive exponential weight, and a negative exponential weight. Panel (d) reports the execution-time ratios for the two exponential cases.}\label{fig:timing-log}
\end{figure}

For the logarithmic covariance, Figure~\ref{fig:timing-log}(a) shows that Method~3 required about \(10.91\) times the execution time of Method~1 for the power family. For positive and negative exponential weights, Method~1 required respectively about \(1.507\) and \(1.614\) times the execution time of Method~4. The maximum coordinate-wise discrepancies remained of order \(10^{-4}\).

Analogous experiments on irregular grids, obtained by sorting independent uniform points on \([0,T]\), produced qualitatively similar rankings.

\subsection{Inference for simulated data}\label{sec:supp-inference}

The calculations use finite-dimensional Gaussian likelihoods.  The corrected
code and the three data files used below are included in the reproducibility
package accompanying this version; the earlier repository is
\url{https://github.com/joseramirezgonzalez/Weighted-sfBm}.  Let
\(t_k=kT/n\), \(k=0,\ldots,n\), and write
\(\bar t_+=(t_1,\ldots,t_n)\).  Since \(\zeta_{0,f,b}=0\) almost surely, the
deterministic value at \(t_0=0\) is excluded from all likelihood covariance
matrices.

For \(i=1,2\), extend the corresponding covariance family to \(b=1\) by
\[
 R_{a,b}^{(i)}(s,t)=
 \begin{cases}
  R_{a,b}^{(i)}(s,t),&b\ne1,\\
  K_a^{(i)}(s,t),&b=1.
 \end{cases}
\]
At \(a=0\), both weight families reduce to \(f\equiv1\); the common
logarithmic covariance is denoted by \(K_1\).  All optimizations were
restricted to \(b\in[0,2]\).  The reported profile intervals use the standard
one-parameter cutoff
\[
 2\{\ell(\widehat\theta)-\ell_p(\theta)\}
 \le \chi^2_{1,0.95}=3.841459.
\]
They describe the individual simulated data sets and are not asserted as
finite-sample coverage results.

\begin{figure}[htbp]
\centering
\includegraphics[width=0.98\linewidth]{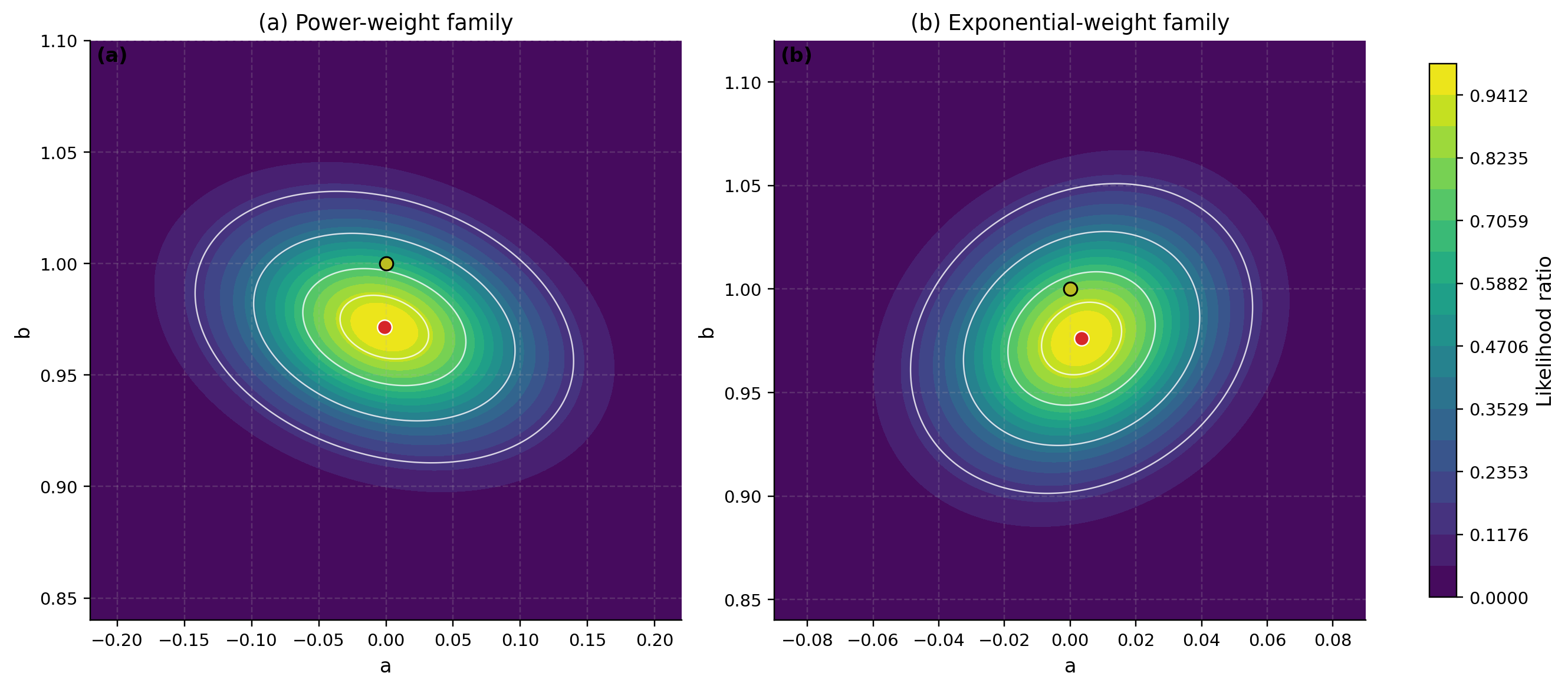}
\caption{Corrected likelihood-ratio surfaces for the two fitted weight
families in the logarithmic-boundary simulation.  The marked points are the
maximizers.}\label{fig:likelihood-surface}
\end{figure}

\begin{figure}[htbp]
\centering
\includegraphics[width=0.98\linewidth]{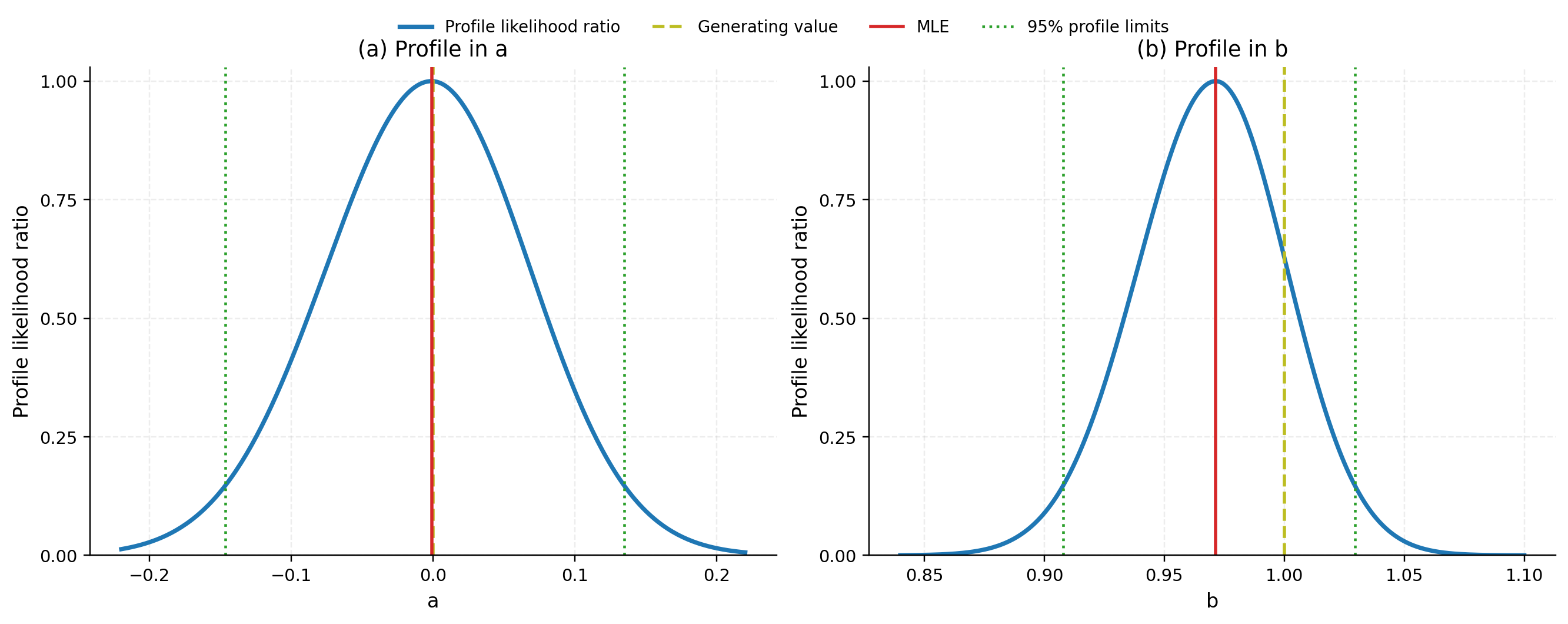}
\caption{Corrected profile likelihood ratios for the power-weight family in
the logarithmic-boundary simulation.  Vertical lines indicate the generating
value, the maximum-likelihood estimate, and the nominal \(95\%\) profile
limits.}\label{fig:power-boundary-profiles}
\end{figure}

In the first experiment, \(T=10\) and \(n=400\).  The supplied centered
Gaussian vector has covariance \(K_1\) on \(\bar t_+\).  The first \(90\%\)
of the grid, corresponding to \(T_{\mathrm{train}}=9\) and
\(n_{\mathrm{train}}=360\), was used for estimation, and the final \(10\%\)
was held out for prediction.  The corrected maximum-likelihood estimates are
\[
 \widehat\theta_{\mathrm{pow}}
 =(-0.001190,\,0.971478),
 \qquad
 \widehat\theta_{\mathrm{exp}}
 =(0.003470,\,0.976083).
\]
For the power family, the nominal \(95\%\) profile intervals are
\[
 a\in(-0.146811,\,0.134797),
 \qquad
 b\in(0.907785,\,1.029506),
\]
and the AIC is \(-1546.872\).  For the exponential family, they are
\[
 a\in(-0.047942,\,0.056165),
 \qquad
 b\in(0.899671,\,1.049246),
\]
and the AIC is \(-1546.889\).  Figure~\ref{fig:likelihood-surface} displays
the corrected likelihood-ratio surfaces, and
Figure~\ref{fig:power-boundary-profiles} displays the power-family profiles.

\begin{figure}[htbp]
\centering
\includegraphics[width=0.98\linewidth]{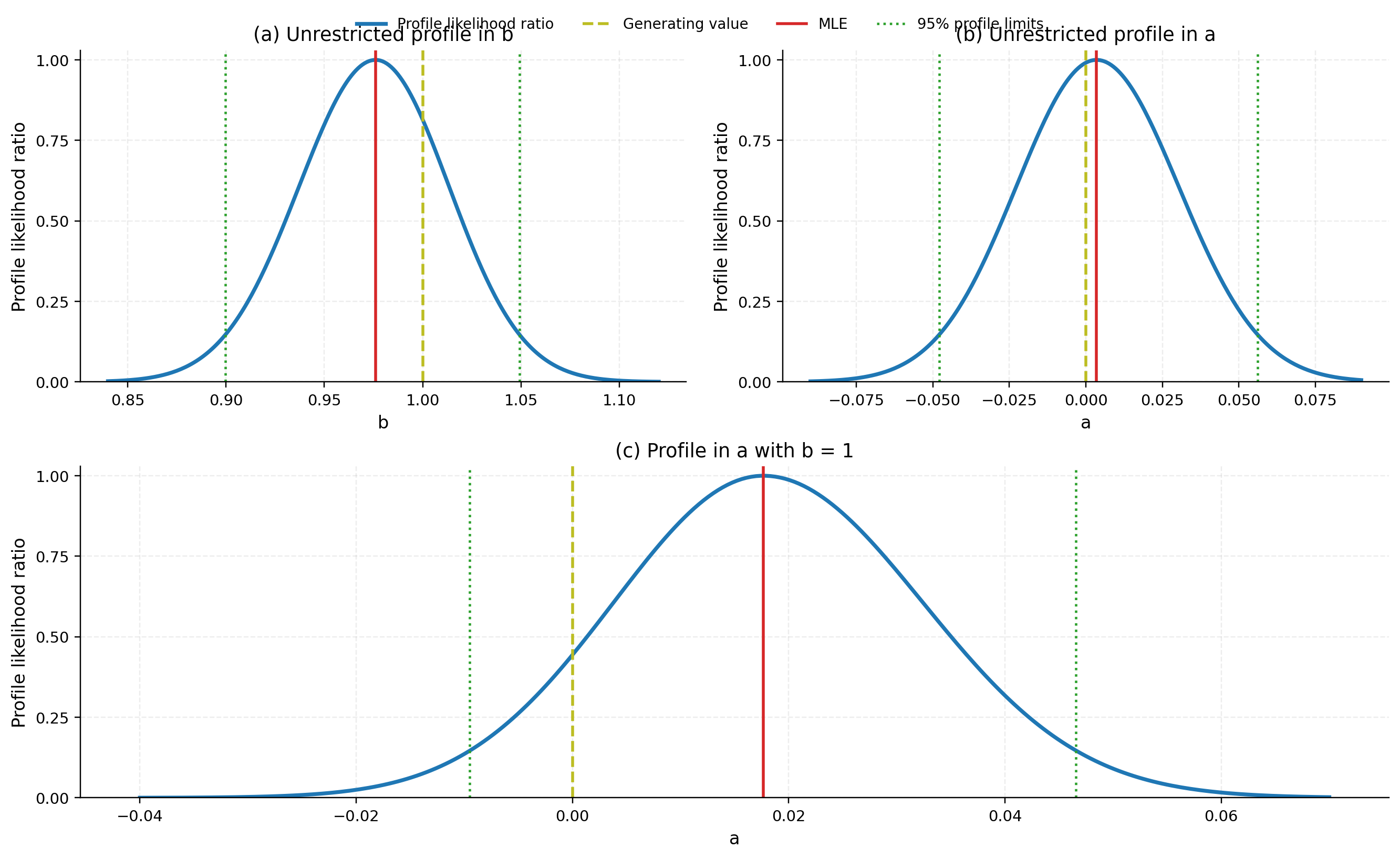}
\caption{Corrected profile likelihood ratios for the exponential-weight
family in the logarithmic-boundary simulation.  Panels (a) and (b) show the
unrestricted profiles in \(b\) and \(a\); panel (c) shows the profile in
\(a\) with \(b=1\).  The discontinuity present in the earlier numerical
implementation is absent.}\label{fig:exponential-boundary-profiles}
\end{figure}

For the exponential family, fixing \(b=1\) gives
\[
 \widehat a_{\,b=1}=0.01765,
 \qquad
 a\in(-0.009481,\,0.046576).
\]
This corrects the previously reported value \(-0.001183\).  In the earlier
code, the fixed-\(b\) optimization call was commented out and \(\widehat a\)
was read from a previously existing optimizer object.  The interval itself
was computed from the fixed-\(b\) likelihood and remains unchanged up to the
reported rounding.  The unrestricted and fixed-boundary profiles are shown
in Figure~\ref{fig:exponential-boundary-profiles}.

The fixed boundary model \(K_1\) has AIC \(-1548.900\).  The training
covariance matrices are ill-conditioned but positive definite: for the fitted
power, fitted exponential, and fixed \(K_1\) models, respectively, the minimum
eigenvalues are \(7.55\times10^{-5}\), \(7.44\times10^{-5}\), and
\(6.66\times10^{-5}\), while the spectral condition numbers are
\(7.79\times10^7\), \(8.03\times10^7\), and \(9.30\times10^7\).  These
diagnostics explain why stable covariance evaluation and Cholesky likelihoods
are necessary.

\begin{figure}[htbp]
\centering
\includegraphics[width=0.98\linewidth]{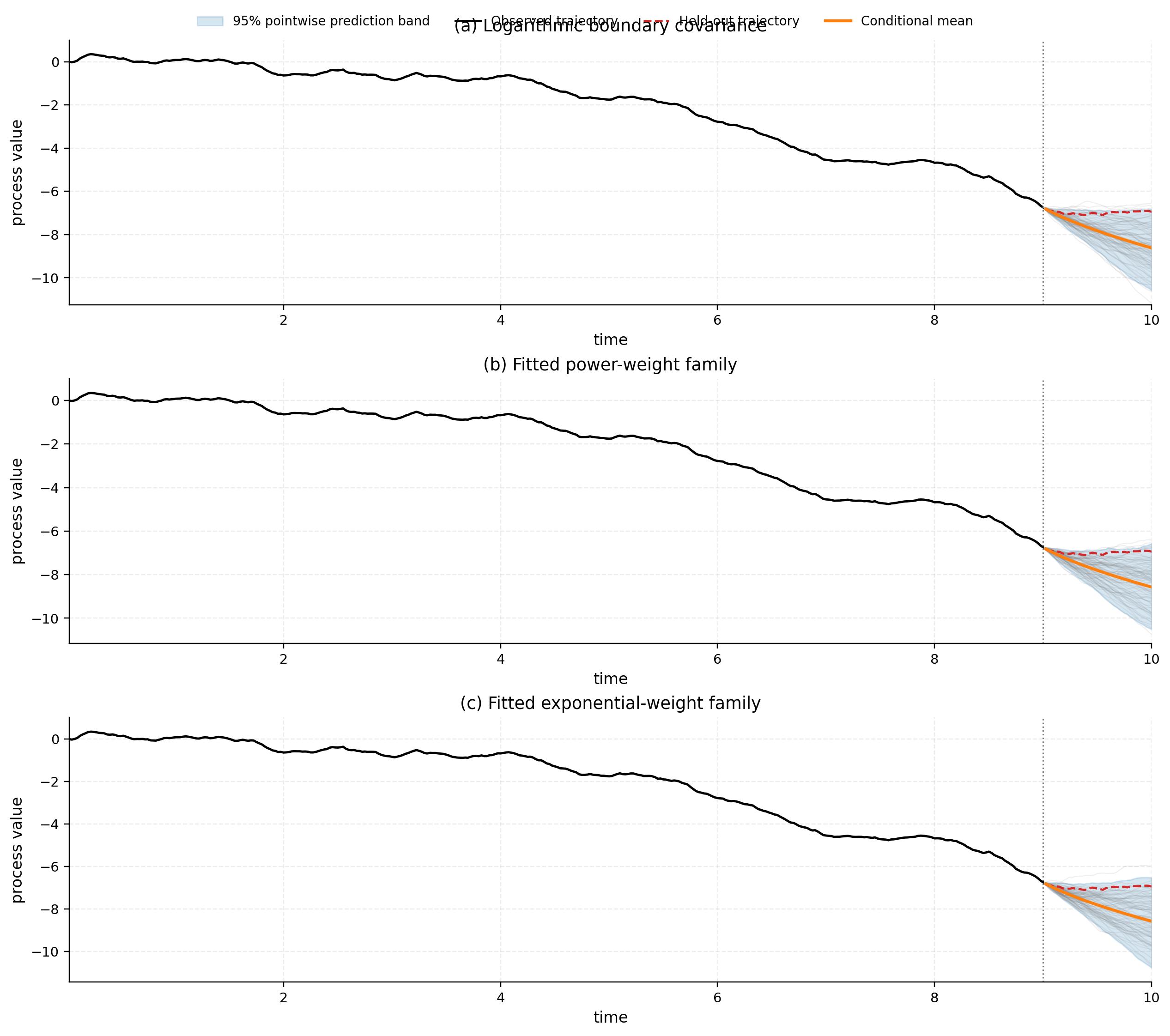}
\caption{Held-out trajectory predictions at the logarithmic boundary and
under the two fitted weight families.  The shaded regions are pointwise
\(95\%\) conditional prediction bands; thin curves are conditional
simulations.}\label{fig:boundary-predictions}
\end{figure}

Predictions for the final \(10\%\) of the grid were computed from the exact
conditional Gaussian mean and covariance.  The mean squared errors of the
conditional means relative to the held-out trajectory are \(0.913125\) for
the fitted power family, \(0.919985\) for the fitted exponential family, and
\(0.966079\) for \(K_1\).  Figure~\ref{fig:boundary-predictions} shows 1000
conditional simulations, their exact conditional mean, and pointwise
\(95\%\) conditional prediction bands.  These are prediction bands, not
confidence bands.

\begin{figure}[htbp]
\centering
\includegraphics[width=0.98\linewidth]{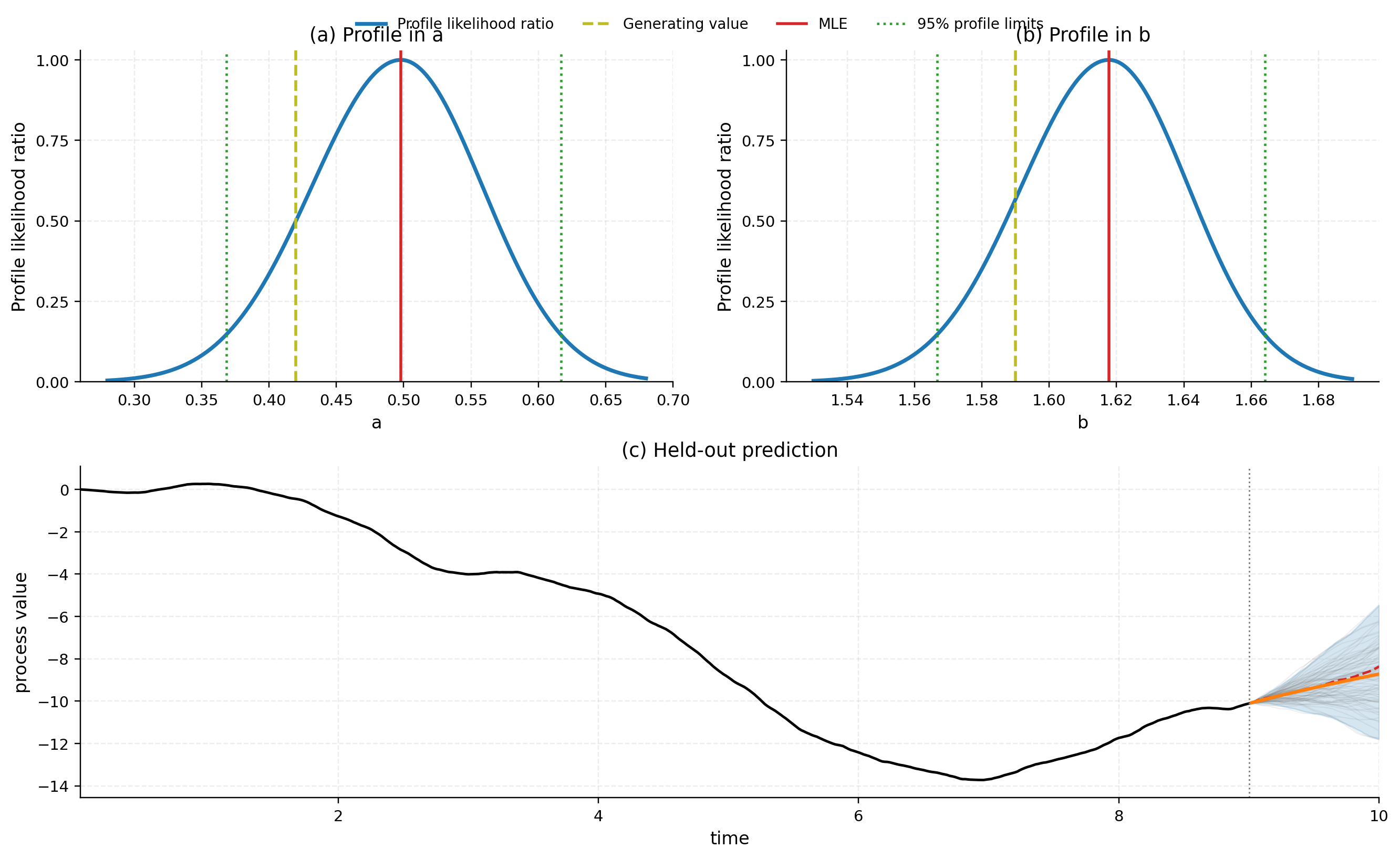}
\caption{Corrected profile likelihood ratios and held-out prediction for the
power-weight simulation with generating parameters \((a,b)=(0.42,1.59)\).
The prediction band is pointwise and conditional.}\label{fig:power-simulation}
\end{figure}

The second experiment uses the power family with generating parameters
\((a,b)=(0.42,1.59)\), \(T=10\), and \(n=500\).  The first \(90\%\) of the
supplied trajectory, corresponding to \(T_{\mathrm{train}}=9\) and
\(n_{\mathrm{train}}=450\), was used for estimation.  The corrected fit is
\[
 (\widehat a,\widehat b)=(0.497986,\,1.617756),
\]
with nominal \(95\%\) profile intervals
\[
 a\in(0.368215,\,0.616754),
 \qquad
 b\in(1.566683,\,1.664032).
\]
The AIC is \(-2950.341\), and the mean squared error of the exact conditional
mean on the held-out \(10\%\) is \(0.012858\).
Figure~\ref{fig:power-simulation} shows the corrected profiles and prediction.

\begin{figure}[htbp]
\centering
\includegraphics[width=0.98\linewidth]{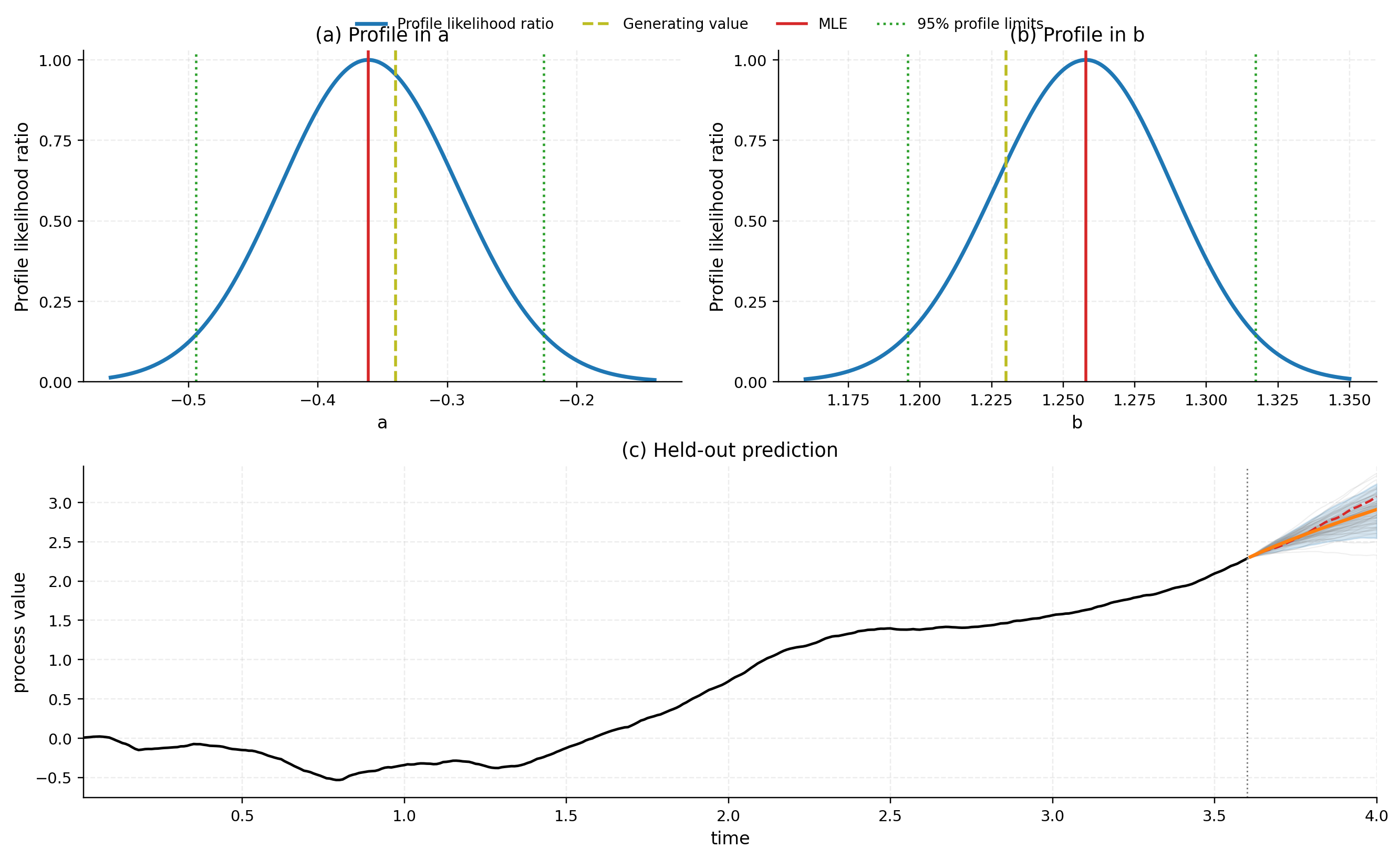}
\caption{Corrected profile likelihood ratios and held-out prediction for the
exponential-weight simulation with generating parameters
\((a,b)=(-0.34,1.23)\).  The prediction band is pointwise and
conditional.}\label{fig:exponential-simulation}
\end{figure}

The third experiment uses the exponential family.  Its generating parameters
are \(a=-0.34\) and \(b=1.23\), with \(T=4\) and \(n=400\).  The first
\(90\%\), consisting of 360 observations through time \(3.6\), was used
for estimation.  The corrected fit is
\[
 (\widehat a,\widehat b)=(-0.361156,\,1.257947),
\]
with nominal \(95\%\) profile intervals
\[
 a\in(-0.494082,\,-0.225511),
 \qquad
 b\in(1.195686,\,1.317231).
\]
The AIC is \(-2872.748\), and the mean squared error of the exact conditional
mean on the held-out \(10\%\) is \(0.005174\).
Figure~\ref{fig:exponential-simulation} shows the corrected profiles and
prediction.

The complete numerical values, covariance diagnostics, profile grids, and
conditional predictions are supplied as CSV files in the reproducibility
package.

\subsection{Sample paths}\label{sec:supp-paths}

\begin{figure}[htbp]
\centering
\iflegacysamplepaths
\subfigure[$f(u)=u^{-0.93},\ b=0.4$]{\includegraphics[width=40mm,height=3cm]{1.5,0.4,a=-0.93.pdf}}
\subfigure[$f(u)=u^{-0.93},\ b=1$]{\includegraphics[width=40mm,height=3cm]{a=-0.93.pdf}}
\subfigure[$f(u)=u^{-0.93},\ b=1.3$]{\includegraphics[width=40mm,height=3cm]{1.5,1.3,a=-0.93.pdf}}
\subfigure[$f(u)\equiv1,\ b=0.4$]{\includegraphics[width=40mm,height=3cm]{1.5,0.4,a=0.pdf}}
\subfigure[$f(u)\equiv1,\ b=1$]{\includegraphics[width=40mm,height=3cm]{a=0.pdf}}
\subfigure[$f(u)\equiv1,\ b=1.3$]{\includegraphics[width=40mm,height=3cm]{1.5,1.3,a=0.pdf}}
\subfigure[$f(u)=e^{-0.8u},\ b=0.4$]{\includegraphics[width=40mm,height=3cm]{1.5,0.4,beta=-0.8.pdf}}
\subfigure[$f(u)=e^{-0.8u},\ b=1$]{\includegraphics[width=40mm,height=3cm]{b=-0.8.pdf}}
\subfigure[$f(u)=e^{-0.8u},\ b=1.3$]{\includegraphics[width=40mm,height=3cm]{1.5,1.3,beta=-0.8.pdf}}
\subfigure[$f(u)=u^{10.6},\ b=0.4$]{\includegraphics[width=40mm,height=3cm]{1.5,0.4,10.6.pdf}}
\subfigure[$f(u)=u^{10.6},\ b=1$]{\includegraphics[width=40mm,height=3cm]{a=10.6.pdf}}
\subfigure[$f(u)=u^{10.6},\ b=1.3$]{\includegraphics[width=40mm,height=3cm]{1.5,1.3,a=10.6.pdf}}
\else
\includegraphics[width=0.92\linewidth]{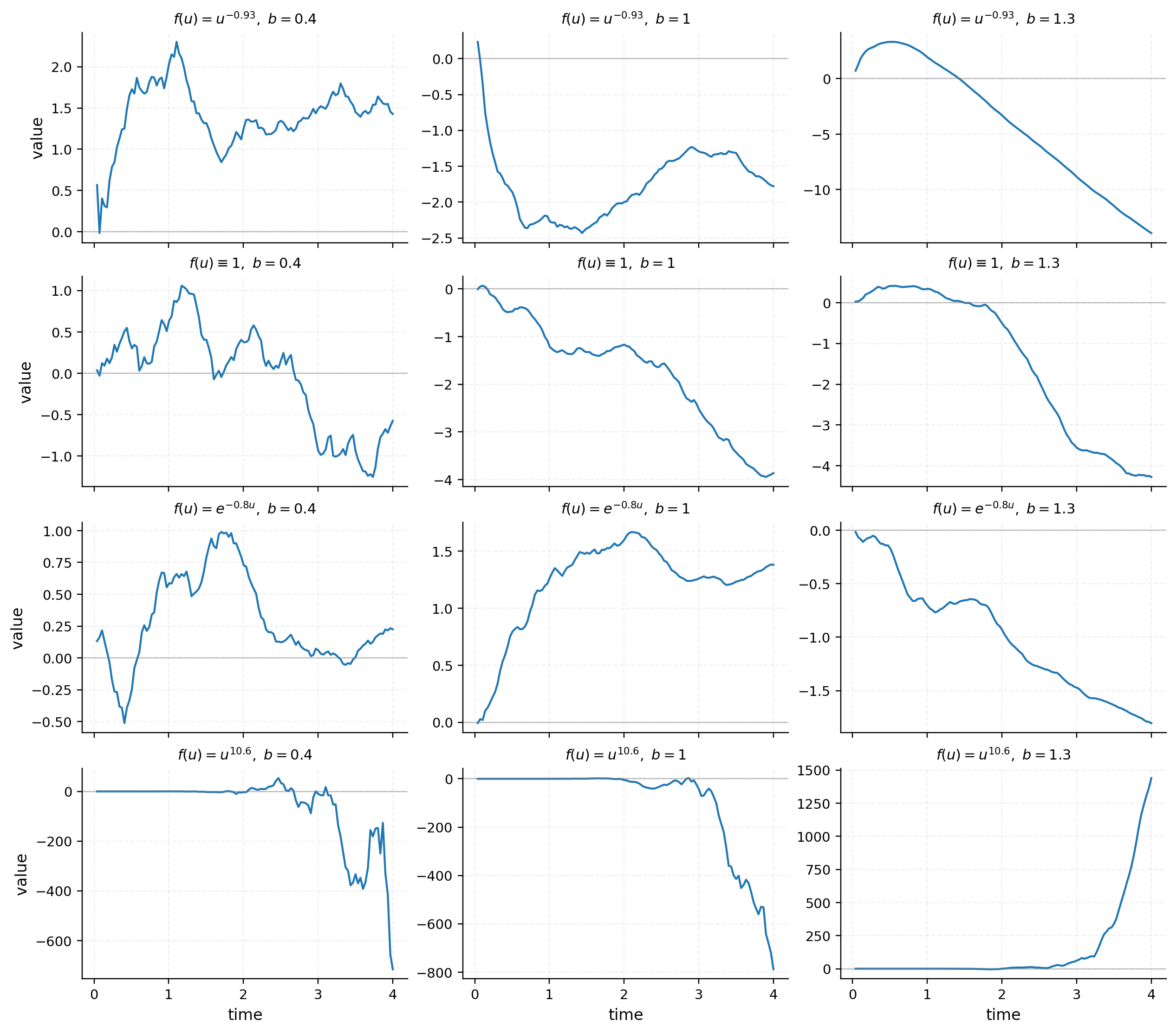}
\fi
\caption{Simulated paths for the parameter values shown in each panel. Panels labeled \(b=1\) use the logarithmic covariance \(K_f\) in \eqref{supp-log-covariance}.}\label{fig:sample-paths}
\end{figure}

\backmatter

\section*{Statements and Declarations}

\noindent\textbf{Funding.}
This work was supported by the S\~ao Paulo Research Foundation (FAPESP), grant 2025/03804-0.

\medskip
\noindent\textbf{Competing interests.}
The author has no relevant financial or non-financial interests to disclose.

\medskip
\noindent\textbf{Data and code availability.}
The code and simulated data used in the numerical illustrations are available at
\url{https://github.com/joseramirezgonzalez/Weighted-sfBm}.

\medskip
\noindent\textbf{Author contributions.}
The sole author carried out all aspects of the work and approved the final manuscript.

\end{document}